\documentclass[11pt,reqno]{amsart}

 \usepackage{curves} 
\usepackage{graphicx} 
\usepackage{amsmath, amsfonts, amssymb, amscd, tikz, hyperref, graphicx, graphics, enumerate, enumitem, graphpap, color, amsfonts, xcolor, cancel, mathrsfs, curves}

\usepackage[mathscr]{euscript}

\makeatletter
\newcommand{\doublewidetilde}[1]{{%
  \mathpalette\double@widetilde{#1}%
}}
  \usepackage[normalem]{ulem}
  \usepackage{soul}
\textwidth= %147.5mm %
148mm
\numberwithin{equation}{section}

\theoremstyle{plain}
\newtheorem{theo}{Theorem}[section]
\newtheorem{lem}[theo]{Lemma}
\newtheorem{prop}[theo]{Proposition}

\newtheorem{cor}[theo]{Corollary}

\newtheorem{lemma}[theo]{Lemma}

\theoremstyle{definition}
\newtheorem{rem}[theo]{Remark}

\newtheorem{definition}[theo]{Definition}

\newenvironment{pf}{\noindent{\it Proof.\,}}{\hfill $\square$\medskip}
\theoremstyle{plain}
\theoremstyle{definition}

\newcommand{\beq}{\begin{equation}}
\newcommand{\eeq}{\end{equation}}
\renewcommand{\a}{\alpha}
\renewcommand{\b}{\beta}

\renewcommand{\d}{\delta}

\newcommand{\f}{\varphi}
\newcommand{\g}{\gamma}
\newcommand{\h}{\eta}

\renewcommand{\l}{\lambda}
\renewcommand{\o}{\omega}

\renewcommand{\r}{\rho}
\newcommand{\s}{\sigma}
\renewcommand{\t}{\tau}

\newcommand{\z}{\zeta}

\newcommand{\bR}{\mathbb{R}}

\newcommand{\bL}{\mathbb{L}}

\newcommand{\bP}{\mathbb{P}}

\newcommand{\bN}{\mathbb{N}}

\newcommand{\bQ}{\mathbb{Q}}

\newcommand{\bM}{\mathbb{M}}

\newcommand{\cA}{\mathscr{A}}
\newcommand{\cC}{\mathcal{C}}

\newcommand{\cE}{\mathscr{E}}
\newcommand{\cF}{\mathscr{F}}
\newcommand{\cG}{\mathscr{G}}
\newcommand{\cH}{\mathscr{H}}

\newcommand{\cJ}{\mathscr{J}}
\newcommand{\cK}{\mathscr{K}}

\newcommand{\cN}{\mathscr{N}}

\newcommand{\cP}{\mathscr{P}}
\newcommand{\cQ}{\mathscr{Q}}

\newcommand{\cS}{\mathscr{S}}

\newcommand{\cU}{\mathscr{U}}

\newcommand{\cW}{\mathscr{W}}

\newcommand{\cZ}{\mathscr{Z}}

\newcommand{\p}{\partial}

\renewcommand{\square}{\kern1pt\vbox
{\hrule height 0.6pt\hbox{\vrule width 0.6pt\hskip 3pt
\vbox{\vskip 6pt}\hskip 3pt\vrule width 0.6pt}\hrule height0.6pt}\kern1pt}

\newcommand{\wt}{\widetilde}

\newcommand{\wh}{\widehat}

\newcommand{\bt}{\begin{theo}\ \ }
\newcommand{\et}{\end{theo}}
\newcommand{\bp}{\begin{prop}\ \ }
\newcommand{\ep}{\end{prop}}
\newcommand{\bc}{\begin{cor}\ \ }
\newcommand{\ec}{\end{cor}}
\newcommand{\bl}{\begin{lem}\ \ }
\newcommand{\el}{\end{lem}}
\newcommand{\bd}{\begin{definition}}
\newcommand{\ed}{\end{definition}}
\newcommand{\be}{\begin{equation}}
\newcommand{\ee}{\end{equation}}

\def\<#1,#2>{\langle\,#1,\,#2\,\rangle}
\newcommand{\arr}{\begin{array}{rlll}}
\newcommand{\ea}{\end{array}}
\newcommand{\bea}{\begin{eqnarray}}
\newcommand{\eea}{\end{eqnarray}}
\newcommand{\bean}{\begin{eqnarray*}}
\newcommand{\eean}{\end{eqnarray*}}

\renewcommand{\=}{:=}

\newcommand{\ve}{\varepsilon}

\newcommand{\pa}{{\bf h}}
\newcommand{\needle}{\mathcal{N}\hskip-2pt\text{\it eedle}}

\newcommand{\Gneedles}{{\cG\!\text{\it ood}\cN}}
\def\sideremark#1{\ifvmode\leavevmode\fi\vadjust{%            The remark
\vbox to0pt{\hbox to 0pt{\hskip\hsize\hskip1em%               will appear only
\vbox{\hsize3cm\tiny\raggedright\pretolerance10000%          on the side
\noindent #1\hfill}\hss}\vbox to8pt{\vfil}\vss}}}%           in 3cm
\title[Control problems with differential constraints of higher order]
{Control problems
with  differential \\ 
constraints of higher order}
  \author[Franco Cardin, Cristina Giannotti and Andrea Spiro]{Franco Cardin \qquad Cristina Giannotti \qquad  Andrea Spiro}
  
\thanks{{\it Acknowledgments}. This research was partially supported by the Projects MIUR ``Regular and stochastic behaviour in dynamical systems''  and ``Real and Complex Manifolds: Geometry, Topology and  Harmonic Analysis''  and by GNFM and  GNSAGA of INdAM}

\begin{document}

\begin{abstract}  We consider cost minimising control problems,  in which the dynamical system is constrained by higher order differential equations of Euler-Lagrange type. Following  ideas  from a previous paper, we prove that a curve of controls $u_o(t)$ and a set of initial conditions $\s_o$ gives an optimal solution for a control problem of the considered type if and only if an appropriate double integral is greater than or equal to zero along any homotopy $(u(t, s), \s(s))$ of control curves and initial data starting from $u_o(t) = u(t, 0)$ and $\s_o = \s(0)$.  This  property is called  {\it Principle of Minimal Labour}.  From this principle we  derive a generalisation of the classical Pontryagin Maximum Principle  that  holds under higher order   differential constraints of Euler-Lagrange type  and without the hypothesis of fixed  initial data. 
\end{abstract} 
%\rightline{\jobname\hskip 10 pt(\today)}
\maketitle
%\tableofcontents
\setcounter{section}{0}
\section{Introduction}
In  \cite{CS}  the first and the third author presented a new  proof of the  classical Pontryagin Maximum Principle (PMP) for controlled  systems,  which was crucially based on the  observation that  the  first order differential constraints  of the system  can be considered as  the Euler-Lagrange equations determined by  an appropriate controlled  first order Lagrangian.  
Following the same   ideas of  that proof,   we  give here a generalisation of   the PMP to the   control problems, in which the differential constraints are given 
by  Euler-Lagrange equations of higher order. \par
\smallskip
More precisely, we consider   cost minimising  problems  for  dynamical system which are controlled through  the  Euler-Lagrange  equations determined by higher order Lagrangians  with  controlling parameters, i.e. by Lagrangians   depending on appropriate controls $u^a(t)$
$$L ^{u(\cdot)}\left(t, q^i, \frac{d q^i}{dt}, \ldots, \frac{d^k q^i}{dt^k}\right) \=  L\left(t, q^i, \frac{d q^i}{dt}, \ldots, \frac{d^k q^i}{dt^k}, u^a(t)\right)\ .$$
We do not impose any  particular  assumption on the initial conditions of the solutions of the differential constraints, nor on the control curves $u(t) = (u^a(t))$, besides
merely  technical  requirements  of smoothness and  constant rank conditions on the constraints.  In fact,  in this paper we restrict  our discussion  just to  control problems satisfying  strong regularity assumptions:  this choice is  only for the purpose of making  as much as possible transparent the main ideas  of   our  approach.  
An illustration  on how the main results of this paper can be generalised  under much   weaker regularity assumptions  is given in \cite{CGS1}. \par
\smallskip
 Considering   the    proof of the PMP  presented in \cite{CS} as a model,  we are able to prove that  a curve of controls $u_o(t)$ and a  set of initial conditions $\s_o$ for the evolution  of the controlled dynamical system  corresponds to  an optimal solution  if and only if a particular double integral is greater than  or equal to zero for any homotopy $(u(t, s), \s(s))$ of control curves and initial conditions, having $u_o(t) = u(t, 0)$ and $\s_o = \s(0)$ as starting point.   We   called this property  {\it Principle of Minimal Labour}. Using this and an appropriate formalisation of   Pontryagin's notion of needle variation of  control curves,  we   derive a generalisation of the classical PMP for the control problems that are subjected  to  higher order Euler-Lagrange constraints   {\it of normal type}. This is a very large class of constraints  which naturally includes the classical Mayer problems to which the usual version of the  PMP applies. Actually our  main result  provides additional  information also for the classical first order
 differential constraints  since   it allows  variations of the initial data
and establishes an innovative shortcut between the proofs of the classical PMP and the Noether 
Theorem about conservation laws (see for instance the approach to Noether Theorem  based on Stokes Theorem given in  \cite{FS},  and compare  it with the  use of Stokes Theorem  in the proof of the PMP  given here and  in \cite{CS}, as  illustrated  in Sect.\ \ref{Sect2} below).   A   discussion in greater detail of  our main results  and some simple illustrations of how they can be exploited  are given  in  Sect.\ \ref{Sect2} and Sect.\ \ref{Sect8}.\par
 \smallskip 
 Before concluding this introduction, we would like  to recall that,   considering an appropriate set of auxiliary variables,  any control problem  with higher order differential constraints -- even those of variational type for  which we establish our generalised PMP --  can be reduced to an equivalent   one with only   first order constraints (but, in general,  no longer of  variational type).   By  introducing other  auxiliary variables, the  Pontryagin variables $p_j$,   the original higher order problem  is finally  transformed into an equivalent one, which is now variational  and  to which the classical PMP applies.  
 This  kind of  reduction  procedure 
 demands the introduction of a (in general very large) collection of auxiliary   variables. In contrast with this, our generalised PMP   directly applies to any set of higher order constraints, with almost no need of  auxiliary  variables.  This is a valuable feature, which we  briefly   illustrate with  an elementary example in Sect.\ \ref{Sect8}.  More substantial  examples and applications will be discussed  in detail   in a future work.\par 
 \smallskip
 This paper is structured as follows. In Sect.\ \ref{Sect2}, we  review the main results of  \cite{CS} and give an informal presentation of  the main ideas  on which  our  results  are based. A  detailed description of the Principle of Minimal Labour and of our generalisation of the PMP are also given.    In  Sect.\ \ref{Sect3} and  Sect.\ \ref{Sect4} the needed preliminaries and a rigorous definition of the class of control problem we are considering are given. In Sect.\ \ref{Sect5},  the considered control problems  are transformed into an equivalent  form, which is  more appropriate 
for the subsequent manipulations. The proofs of the Principle of Minimal Labour and of our generalisation of the PMP are given in Sect.\ \ref{Sect6} and Sect.\ \ref{Sect7}, respectively. In Sect.\ \ref{Sect8}, we illustrate  some of the main features of  our approach by discussing  a couple of   elementary   problems.  
\par
\medskip
\section{An overview of our approach and main results} \label{Sect2}
 Since 
the results of the present paper can be considered as natural developments of the  ideas of  \cite{CS}, we decided to precede our discussion  with    a short 
 overview  of  the contents of that paper.  We then  briefly indicate how the scheme of such a previous  paper is here  implemented  to obtain our new results.\par
\smallskip   
\subsection{Pontryagin Maximum Principle and Stokes Theorem in a nutshell}
\subsubsection{The basic scheme of a classical Mayer problem} 
\label{section211} 
Consider a dynamical system, whose evolutions on a fixed time interval $[0, T]$ are represented by  curves $x(t) = (x^1(t), \ldots, x^N(t))$ in $\bR^N$  satisfying  the differential problem 
\beq \label{Mayer1} \frac{d x^i}{dt}(t) = f^i(t, x^i(t), u^a(t)) \ ,\qquad x^i(0) = x^i_o \in \bR^N\ . \eeq
Here $u(t)$ is a  (measurable) function with values in a fixed subset $K$ of $\bR^M$
$$u(t) = (u^1(t), \ldots, u^M(t)) \in K \subset \bR^M$$
 and represents the evolution in time of  control parameters.  The  $f^i(t, x^i, u^a)$  are continuous  functions   on $[0, T] \times \bR^N \times K$ and differentiable in the $x$. The  initial value $x_o$  is  fixed and is the same   for  all of   the evolutions of the system.  Given a  terminal cost function $C: \bR^N \to \bR$,  the corresponding  {\it Mayer problem} 
  consists of  looking for  a curve $\bar u(t)$  of control parameters, for which  the following holds:  {\it the terminal cost $C(\overline x(T))$ of the  solution $\overline x(t)$ to   \eqref{Mayer1} with $u(t) = \overline u(t)$ is less than or equal to the terminal cost of  the solution determined by  any other   choice of the  curve  $u(t)$.}\par
Of course, this is only   one   of the many variants   of the classical Mayer problem.    But  in what follows we  limit ourselves     to such a basic  version. Moreover,  we constantly  assume that   all  the  data    satisfy  much higher regularity assumptions than those mentioned  above.  Take this restriction as  a sort of blanket assumption, which  is  adopted  to easily allow     manipulations   and prevent the risk of diverting  the attention of the reader towards inessential  technical issues.\par
\subsubsection{The  auxiliary variables $p_i$}  \label{Cardinuccio} The classical approach to  a  Mayer problem is usually based on  the  introduction of a set of    auxiliary variables  $p_1, \dots, p_N$
and on the representation of  the dynamical system  through curves 
$(x(t), p(t)) =   (x^i(t),   p_j(t))$ in $\bR^{2N}$ (of which, however,   only  the  $x^i(t)$ are the   interesting ones),   constrained by  the  \eqref{Mayer1}   {\it and} the auxiliary   equations
\beq  \label{Mayer2} \frac{d p_j}{dt} \=- \sum_{i = 1}^N p_i(t)  \frac{\p f^i}{\p x^j}\bigg|_{(t, x^i(t), u^a(t))} \ .\eeq
The introduction of the auxiliary  variables $p_j$ and of the  new constraints \eqref{Mayer2}  has the following effect:  {\it for each  curve of control parameters  $u(t) = (u^a(t))$,  the constraints on the corresponding curve $(x(t), p(t)) $  are    the Euler-Lagrange equations of the  variational principle determined by the (controlled) Lagrangian }
\beq  \label{classical-lagrangian}L ^{u(\cdot)} (t, p, x, \dot x)= \sum_{j = 1}^N p_j(\dot x^j - f^j(t, x, u(t))) \ .\eeq
{\it This Lagrangian has also the special  feature of being identically vanishing along the solutions of the corresponding  Euler-Lagrange equations (in fact the $x^i(t)$ satisfy the \eqref{Mayer1}).}\par
All this has  the consequence   that  the original cost minimising problem   is   equivalent to the following  one.  By subtracting a constant,  with no loss of generality we may  assume that $C(x_o) = 0$. Hence by the vanishing of $L^{(u(t))}$ along the constrained curves,  the value on such curves of the functional 
\beq \label{action1} I^{(u(\cdot))} \= \int_0^T  \bigg(L ^{u(\cdot)} (t, p, x, \dot x) + \sum_{i = 1}^N \frac{\p C}{\p x^i }\bigg|_{(t, x(t))}  \dot x^i(t) \bigg) dt\eeq
appears to be  equal to 
\beq \label{smart-1} I^{(u(\cdot))} = \int_0^T  \frac{d C(x(t))}{dt}\bigg|_{(t, x(t))} dt = C(x(T)) - C(x(0)) \overset{C(x_o) = 0}= C(x(T))\ . \eeq
Thus the original problem turns out to be the same of {\it looking for a curve of controls $u(t) = (u^a(t))$ such that, along the solution of the Euler-Lagrange equations of $L ^{u(\cdot)} $,   the  functional} 
 $I^{(u(\cdot))}$ 
{\it takes the minimum possible value}. \par
\smallskip
We also observe  that, along the solutions $(x(t), p(t))$, the value of the above action  is   independent of  any  boundary (initial or final) conditions for the  $p_j(t)$.   This means that {\it any} value at $t = 0$ (or at $t = T$) can be  imposed  on the  $p_j(t)$, having  absolutely  no consequences on the minimising problem.  As we will see in the next two subsections,  among all of the possible choices for such boundary conditions,  some    are much more   convenient than the others.  \par
\subsubsection{A smart choice for the boundary values  of the  $p_i$}\label{smart-2} Consider  a homotopy of control curves $u^{(s)}(\cdot): [0, T] \to K$, $s \in [0,1]$,  and denote by 
$(x^{(s)i}(t), p_j^{(s)}(t))$ the corresponding homotopy of curves satisfying the constraints  \eqref{Mayer1} and \eqref{Mayer2}.   Exploiting our  blanket assumptions on the  regularity of the  data,  we may say  that  the first order jets  $(x^{(s)i}(t), \dot x^{(s)i}(t), p_j^{(s)}(t), \dot p_j^{(s)}(t))$,  determined by the pairs  $(t, s)  \in [0, T] \times [0,1]$,   span a  smooth surface $\cS$ in the first order jet space  of curves in $\bR^{2N}$, whose boundary is formed 
by four smooth curves. We denote them by  $\g^{(s=0)}(t) $, $\g^{(s=1)}(t) $, $\h^{(0)}(s)$ and $\h^{(T)}(s)$. \\
 \centerline{
\begin{tikzpicture}
%\draw[help lines] (0,0) grid (4,3.2);
\draw[fill]  (0.5, 0.2) circle [radius = 0.05];
\draw[fill]  (3.8, 0.7) circle [radius = 0.05];
\draw[fill]  (0.25, 2.3) circle [radius = 0.05];
\draw[fill]  (3.7, 2.85) circle [radius = 0.05];
\node at  (2, 2.85) {\color{blue} \tiny $\g^{(s = 1)}(t)$};
\node at  (2, 0.55) {\color{blue} \tiny $\g^{(s= 0)}(t)$};
\node [blue] at   (0, 1) {\tiny$\h^{(0)}(s)$};
\node [blue] at   (4.1, 1.7) {\tiny$\h^{(T)}(s)$};
\node [red] at  (2, 1.8) {$\cS$};
\draw [ thin, blue] (0.5, 0.2) to [out=45, in=180] (2,0.85) to [out=0, in=207] (3.8, 0.7)  ; 
\draw [ ->, thick, blue] (2, 0.85) to  (2.02, 0.85)  ; 
\draw [ thin, blue] (0.25, 2.3) to [out=45, in=180] (1.5, 2.65) to [out=-5, in=207](3.7, 2.85)   ; 
\draw [ ->, thick, blue] (2, 2.61) to  (2.02, 2.61)  ; 
\draw [ thin, blue] (0.5, 0.2) to [out=65, in=290] (0.25, 2.3)  ; 
\draw [ ->, thick, blue] (0.52, 1.5) to  (0.515, 1.51)  ; 
\draw [ thin, blue] (3.8, 0.7) to   [out=155, in=270] (3.3, 1.7) to
 [out=90, in=310]  (3.7, 2.85)  ; 
\draw [ ->, thick, blue] (3.3, 1.7) to  (3.3, 1.71)  ; 
\draw [ thin, red] (0.56, 1.25) to [out=45, in=180] (1.5, 1.66) to [out=-5, in=207](3.3, 1.8)   ; 
\draw [ thin, red] (0.61, 0.75) to [out=45, in=180] (2, 1.26) to [out=-5, in=207](3.35, 1.15)   ; 
\draw [ thin, red] (0.44, 1.9) to [out=45, in=180] (2, 2.31) to [out=-5, in=207](3.63, 2.3)   ; 
\draw [ ->, thick, red] (2, 2.31) to  (2.1, 2.31)  ;
\draw [ ->, thick, red] (2, 1.26) to  (2.1, 1.26)  ;
\draw [ ->, thick, red] (2.2, 1.61) to  (2.3, 1.61)  ;
%\put(20,5){\circle*{3}}
%\put(115,20){\circle*{3}}
%\put(15,65){\circle*{3}}
%\put(112,82){\circle*{3}}
%\curve(20,5,30,15, 40, 20, 50,25, 60,26, 70, 25, 80, 23,  85,18, 115,20)
%\put(85,10){\tiny$\g^{(U_o)}$}
%\put(58,55){\color{red}\tiny$\g^{(U(s))}$}
 \end{tikzpicture}
 }
% \vspace{-0.3cm}
  \centerline{\tiny \bf Fig. 1}
The first  two curves correspond to the sides $[0,T] \times \{0\}$,  $[0, T]\times \{1\}$ of  $\p([0, T] \times [0,1])$: 
$$\g^{(s=0)}(t) = (x^{(0)i}(t), \dot x^{(0)i}(t),  p_j^{(0)}(t), \dot p_j^{(0)}(t))\ ,$$
$$ \g^{(s=1)}(t) =(x^{(1)i}(t), \dot x^{(1)i}(t),  p_j^{(1)}(t), \dot p_j^{(1)}(t))\ .$$
The others  are the curves that correspond  to the remaining  two sides  $\{0\} \times [0,1]$,  $\{T\}\times[0,1]$  of  the boundary of  $[0, T] \times [0,1]$. 
We call them the {\it vertical sides} of  $\cS$.  In \cite{CS} it was observed that {\it if one imposes that the  $p_j(t)$  satisfy the terminal values conditions 
\beq \label{smart} p_j(T) = - \frac{\p C}{\p x^j}\bigg|_{(x(T))}\ ,\eeq
 then the integrals of  the $1$-form $\b = \left(L ^{u(\cdot)}(t, p, x, \dot x)  + \sum_{i = 1}^N \frac{\p C}{\p x^i }|_{(t, x(t))}  \dot x^i(t)\right) dt$  along the  two vertical sides $\h^{(0)}(s)$,  $\h^{(T)}(s)$ of $\cS$ (or, more precisely, of an appropriate modification  of $\b$, modelled on  the classical Poincar\'e-Cartan  $1$-form) are   equal to $0$.} \par
 \smallskip
This is an important property, because in combination  with \eqref{smart-1}   it  implies that the integral  of (the Poincar\'e-Cartan type modification of) $\b $  along the anti-clockwise oriented boundary  $\p \cS$ is  equal  to $ - C(x^{(1)}(T)) + C(x^{(0)}(T))$. On the other hand, by Stokes Theorem, such an  integral  is equal to the  integral  of  the exterior differential of the Poincar\'e-Cartan type modification of $\b$ on $\cS$.  By exploiting certain properties of the  actions and  of the $1$-forms of Poincar\'e-Cartan type  (it is not essential  to recall them here - in Sect.\ \ref{encapsulating}   we  discuss them in greater detail),  such an integral  reduces to  a very simple form, namely to
\beq \label{astast}  \begin{split} C(x^{(1)}(T)) - C(x^{(0)}(T))= &- \iint_{t\in [0,T], s\in [0,1]}\frac{{\p \cH}}{ \p u^a }\frac{\p u^a }{\p s}\bigg|_{(t, x^{(s)i}(t), p^{(s)}_j(t), u^{(s)a}(t))}\, dt\, ds\ ,\\
& \text{where}\ \cH(t, x^i, p_j, u^a) \= - \sum_{j = 1}^N p_j f^j(t, x^i, u^a)\ .
\end{split} 
\eeq
The main reason of interest for this identity comes from the fact that  {\it it expresses  the difference between the two terminal costs as  a double integral of an appropriate  function of the   parameters $(t,s)$ of the  homotopy.}   From our point of view, this is   a cornerstone in the  proof of the PMP. 
\par 
\subsubsection{The Principle of Minimal Labour and the Pontryagin Maximum Principle}  
\label{section2.1.4} An immediate consequence of \eqref{astast} is the following: 
\\[10pt]
 {\bf Principle of Minimal Labour}. {\it If the curve $\bar u(t) \in K$ gives a solution to the considered Mayer problem, then for any other curve $u(t)$ which is connected to $\overline u(t)$ through an homotopy  of curves $u^{(s)}(t)$ in $K$, the double integral on the left hand side of  \eqref{astast} is less than or equal to $0$}. \\[10pt]
  By   considering   appropriate  highly localised deformations of the curve  $\bar u(t)$  (the so-called {\it Pontryagin's needle variations}  -- see \cite{CS} or Sect.\ \ref{generalisedneedle}  below  for details) and associated interpolating homotopies,   the classical PMP can be derived as  if a `pointwise version'  of the  above principle. Indeed,  adopting a very informal language, we may state the {\bf Pontryagin Maximum Principle}  as follows (the literature on such a  classical principle is vast -- for  extensive  and fundamental presentations we refer to  \cite{PBGM, Su3, Ju, AS, BP, Ce} and references therein): \\
 [10pt]
 {\it  If  $\bar u(t) \in K$ gives a solution $(x^i(t))$ to the considered Mayer problem, then for any $t_o \in [0,T]$ the value of  $\cH$ at  the point $(t_o, x^i(t_o), p_j(t_o), \bar u^a(t_o))$
   is maximal among all the  values that  it assumes at the points  $(t_o, x^{(\o)i}(t_o), p^{(\o)}_j(t_o), \o^a)$ determined by  
   \begin{itemize}
   \item[(i)]  replacing   $(\bar u^a(t_o))$ by some other value  $(\o^a)\in K$,
   \item[(ii)]   substituting   the values $x^i(t_o)$ and  $p_j(t_o)$ by the values 
   $x^{(\o)i}(t_o)$ and $p^{(\o)}_j(t_o)$, which are  assumed by  the  solution  $(x^{(\o)i}(t), p^{(\o)}_j(t))$  at $t_o$ of the  constraints determined by a control curve $u^{(\o)}(t)$   taking the value  $u^{(\o)}(t_o) = \o$ in an infinitesimal  neighbourhood of   $t = t_o$ and  coinciding with the original  $\bar u(t)$ at all other points.
   \end{itemize}
   }
\par\smallskip
\subsection{Our road map  towards the main results}\hfill\par
\label{roadmap}
Let us now focus on the following two facts, pointed out  in  the above summary of  \cite{CS}. 
\begin{itemize}[leftmargin = 20pt]
\item[(1)] The   problems,  to which the classical PMP applies,  are costs minimising problems on curves $x(t) = (x^i(t))$ that are controlled by means of  {\it first order differential equations with parameters}.
\item[(2)]  By introducing  auxiliary variables  $p_j$ and an appropriate family of controlled Lagrangians   $L^{u(\cdot)}(t, p, x, \dot x)$, the differential constraints of the original control problem are    replaced by 
the Euler-Lagrange equations of such Lagrangians. 
\item[(3)] The   Lagrangians  $L^{u(\cdot)}(t, p, x, \dot x)$ have the following   peculiar property: {\it they vanish identically   along  the solutions of  their corresponding Euler-Lagrange equations}.  This is a crucial fact that leads to the identity 
\eqref{smart-1}.
\item[(4)]  If  appropriate boundary conditions are imposed on the  $p_j(t)$,  then for any given homotopy $u^{(s)}(t)$ of   control curves, 
the integrals of the actions of Poincar\'e-Cartan type along the ``vertical sides'' $\h^{(0)}(s)$,  $\h^{(T)}(s)$ of  the surface in Fig.\ 1  are identically vanishing. This property  together with (3) leads to the  identity \eqref{astast}, 
which expresses  the difference between two terminal costs as a  double integral of an appropriate function of the parameters of the homotopy.
\item[(5)]   The identity \eqref{astast} immediately gives the Principle of Minimal Labour, from which the PMP is derived   using   highly localised  (needle)  variations. In a sense, the Principle of Minimal Labour can be taken as an  underlying substratum  for the PMP.
\end{itemize}
In this paper we consider a special class of cost minimising problems, in which the  curves  are constrained  by  {\it differential equations of  higher  order and of variational type}, that is  by  Euler-Lagrange  equations of  controlled Lagrangians of higher order (Sect.\ \ref{definingtriples}).  For such problems we are able to   follow the same circle of ideas described above and, at the end, we reach  a 
generalised version of the PMP that works for  this wider class of cost minimising problems.  \par
We remark that, since  our  differential constraints   are assumed to be  of  Euler-Lagrange type,  there is no need  to introduce     auxiliary variables and  new constraints in order  to obtain the  property described in the point (2) above: It is  granted from the very beginning. However,  since  we are no longer requiring that the  controlled Lagrangian is of  a very  special form,  the crucial phenomena described in the points  (3) and (4)  are in general not occurring.    We nonetheless manage to  overcome this difficulty through the following  two  steps.
\begin{itemize}[leftmargin = 20pt]
\item[(a)] We consider  a special set of functions, denoted by $\pa^i_\b(t)$, $\pa'{}^i_\b(t)$ and  $\pa''{}^i_\b(t)$,  which are completely determined  by  the initial and  the terminal  points of each   controlled  curve. Such  new  functions are  used in a convenient way to   modify the controlled Lagrangian  and  obtain an analogue of   the phenomenon described in (4)  (Sect.\ \ref{encapsulating-bis-A}). 
\item[(b)] We  introduce  two auxiliary variables, called $\l$ and $\mu$, and we further adjust   the   Lagrangian  in order to obtain a new controlled Lagrangian. The purpose of this is to get  an additional  property, which is analogous  to the one described  in    (3). As a consequence,   we get an identity of the form   \eqref{smart-1}  (Sect.\ \ref{encapsulating}). We  stress the fact that the introduction of  the new variables  $\l$ and $\mu$ is a merely technical expedient  and that 
 such variables  do not occur  in the statements of the final  results.  
\end{itemize}
The modifications described in  (a) and (b) have been found   by heuristic arguments.   We  do not know   whether these  are the only possible ones and/or   there  are deep reasons for   why they work (\footnote{We guess that  a more elegant approach should exist. For instance,  the idea we used for the modifications  described in (b)  
 calls to mind   a  well known trick, which is usually exploited  to translate   a    Bolza problem into  a  Mayer problem (see e.g.  \cite{BP}, p.\ 116).}).  Nonetheless  for the purposes  of the present paper we only need   to  know that they work. In fact, we would like to stress that  their main use  is basically  to re-write  the sum of a particular  surface integral and of   two  boundary line integrals  into a single surface integral. This yields  an elegant expression  for the Principle of Minimal Labour, which nonetheless  is by no means the only possible one. Other   different equivalent  statements are admissible,  which  can be proven  with no need of the above mentioned modifications, but have the disadvantage of  being much  more involved.  \par
\smallskip
Performing the steps (a) and (b) and  following  essentially the   same ideas used in \cite{CS}   we   finally  get  the desired   analogues of the Principle of Minimal Labour and of the  PMP (Sect.\ \ref{generalisedneedle}).  An    informal description of such results is given  in the next subsection.  \par
\smallskip
We conclude inviting the reader   to consider the outline of this section  as a road map for the following constructions and  to constantly  keep  it  in mind while going  through the rest  of  the  paper. \par
 \smallskip
\subsection{Main results} \hfill\par
The outcomes of the above described approach are basically two. The first is a principle  (Theorem \ref{corollone}) that can be considered as a generalisation  of the first version of the Principle of Minimal Labour proved in  \cite{CS}. It  can be   described  as follows.\\[10pt]
 {\bf  Principle of Minimal Labour}. {\it Consider a terminal  cost minimising problem for evolutions $q(t) = (q^i(t))$, $t \in [0,T]$, which are constrained by  a set of smooth equations on the initial values and by a system of ordinary differential equations  of  Euler-Lagrange  form 
 $$\frac{\p L^{u(\cdot)}}{\p q^i}\bigg|_{q(t)} + \sum_{\b= 1}^r (-1)^\b \left(\frac{d}{dt}\right)^\b\left(\frac{\p L ^{u(\cdot)}}{\p \frac{d^\b q^i}{dt^\b}}\right)\bigg|_{q(t)} = 0\ .$$
 Here, $L ^{u(\cdot)} = L ^{u(\cdot)}(t, q^i, \frac{d q^i}{dt}, \ldots, \frac{d^r q^i}{dt^r})$ is a family of Lagrangians of order $r \geq 1$, which   smoothly  depends on the values of  a curve $u(t) = (u^a(t))$ of control parameters. \par
A  curve $\bar u(t)$ and   admissible initial conditions determine   a solution for the considered  cost minimising problem only if  for any other curve $u(t)$ and any other set of admissible initial conditions, which  can be  joint to  the previous  by a smooth one-parameter deformation,  the following inequality holds
  \beq \label{formulacor-intro}  \int_0^T \left( \int_0^1    \left( \frac{\p \cP_{(t,s)} }{\p u^a} \frac{\p u^{(s)a}}{\p s} \Bigg|_{u^{(s)}(t)} - \frac{\p^2 \mu'}{\p t\,\p s}\bigg|_{(t,s)}\right)ds \right) dt \leq 0 \ .\eeq
  Here:
\begin{itemize}[leftmargin = 15pt]
\item $s \in [0,1]$ is the parameter of the smooth deformation of initial data and control curves and  $\frac{\p u^{(s)a}}{\p s} $ are the derivatives with respect to $s$ of the components of  the curves $u^{(s)}(t) = (u^{(s) a}(t))$  
of  such deformation; 
\item $\cP_{(t,s)} $ is the two-parameters family of  functions of the values  $u = (u^a)$ defined by 
$$\cP_{(t,s)}(u) = - L^{u(\cdot) \equiv u}\left(t, q^{(s)i}(t), \frac{d q^{(s)i}}{dt}\bigg|_t, \ldots, \frac{d^r q^{(s)i}}{dt^r}\bigg|_t\right)\ ,$$
where  $(q^{(s)i}(t))$ denotes  the  evolution of the system, which is  determined by the control curve $u^{(s)}(t)$ and the initial  data associated with  the deformation parameter  $s$; 
\item  $ \mu' =  \mu'(t, s)$ is the function which is defined in    \eqref{7.17'};   it  is indeed the time  integral between $0$ and $t$   of  a certain function, which is explicitly  given in that formula and it is determined by  the following  three sets of objects:  
\begin{itemize}
\item[(A)]  the curves $u^{(s)}(t)$ of the homotopy  of the control curves; 
\item[(B)]  the values and the derivatives  up to the  order $2r$ of the  curves $(q^{(s)i}(t))$; 
\item[(C)] the infinitesimal variations  of the terminal costs  of these curves w.r.t. $s$. 
\end{itemize}
\end{itemize}
}
 For a  classical Mayer problem, where  the dynamical system is described by curves   
 \beq  q(t)= ( x^i(t), p_j(t))\qquad \text{with}\qquad p_j(T) = - \frac{\p C}{\p x^j}\big|_{(x(T))}\ ,\eeq
   the above  principle  radically simplifies. More precisely one has that: 
\begin{itemize}[leftmargin = 15pt]
\item[(1)] For any homotopy $q^{(s)}(t)$ of controlled curves of the above type, the  double integral  $\int_0^T \int_0^1 \frac{\p^2  \mu'}{\p t\,\p s}\big|_{(t,s)}ds\,  dt$ in \eqref{formulacor-intro}  vanishes.   In fact,    the    end-point constraints on the   $p_j$
force the components  $p_j(t,s) \=    p^{(s)}_j(t)$ of each homotopy $q^{(s)}(t)= ( x^{(s)i}(t), p^{(s)}_j(t))$ to play the role of surrogates  for the    $\mu'(t, s)$. Indeed  each $p_j(t,s)$    is  
uniquely  determined by  the curves $u^{(s)}(t)$ and  $x^{(s)}(t) = (x^{(s)i}(t))$,  and  by
 the infinitesimal variations  of the terminal costs at $t = T$. This occurs in perfect  analogy with the properties  (A), (B), (C) of the function  $\mu'(t,s)$. 
\item[(2)] The partial derivatives  $ \frac{\p \cP_{(t,s)} }{\p u^a}$ coincide with the partial derivatives $\frac{{\p \cH}}{ \p u^a }$. 
\item[(3)] The  principle we just mentioned reduces to the  Principle of Minimal Labour  presented in  \cite{CS}  (see Sect.\ \ref{classicalPMP} below for details). 
\end{itemize}
In contrast with all this,  for  other kinds of  cost minimising  problems,   no analogues of the auxiliary variables  $p_j$ are involved and  the 
term $\int_0^T \int_0^1 \frac{\p^2    \mu'}{\p t\,\p s}\big|_{(t,s)}ds  dt$ in \eqref{formulacor-intro}  cannot    be expected to be  zero. \par
\smallskip
The second main result of our paper is obtained by applying   the above Generalised Principle of Minimal Labour to 
the case of  highly localised (``needle'') variations.  Indeed, what we obtain can be  considered as  an analogue of the classical PMP  for the above mentioned  large class of the higher order Euler-Lagrange constraints    {\it of normal type} (see Sect.\ \ref{defnormal} for the precise definition). It  consists of a  necessary condition for a  control curve  $u_o(t)$ in order to determine an  optimal solution $\g_o$ and   can be roughly described as follows.   {\it Let   $\cP$ be    the  function on $K$ defined by 
$$
 \cP : K \longrightarrow \bR \ ,\qquad  \cP(u^a) \= -   L\bigg(t,  \g_o(t), \ldots, \frac{d^r \g}{dt}, u^a\bigg)\ . $$
Then, $u(t)$ gives an optimal solution $\g_o$  only if,  for any $\t \in [0, T]$ and any $\o \in K$,  for  which the curve $u_o(t)$ admits a smooth  deformation $u_{(s)}(t)$ of needle type   around $t = \t$ with $u_{s = 1}(\t) = \o$, the following inequality  holds }
 \beq  \label{PO} \cP(\o)  -  \liminf_{\ve \to 0^+}  \frac{  \mu'{}^{(\t, \o, \ve, \Sigma)}(T, 1) -  \mu'{}^{(\t, \o, \ve, \Sigma)}(T, 0) }{\ve}   \leq \cP(u_o(\t))\ .\eeq
Here,   $ \liminf_{\ve \to 0^+}  \frac{  \mu'{}^{(\t, \o, \ve, \Sigma)}(T, 1) -  \mu'{}^{(\t, \o, \ve, \Sigma)}(T, 0) }{\ve}  $ is a corrective term,  which can be  determined by means of the data of  the needle variation, namely: (i) the point   $t = \t$ where it is applied, (ii) the width $\ve$, (iii) the top value $\o$ and  (iv) the $1$-parameter family $\Sigma$ of initial or terminal values for some variables, as e.g. the conditions \eqref{smart} in the  classical Pontryagin setting (see Sect.\ \ref{Sect7} for details).  For  the  classical Mayer problems with first order constraints,  the above corrective term is   zero and the resulting  condition   \eqref{PO} on the function $\cP$ reduces   to the  usual  PMP on  the Pontryagin function $\cH$.   For what concerns  more general cases with higher order constraints,  we offer a   characterisation of the needle variations,   which  allow to neglect such corrective term. This  yields to an alternative version  of the usual PMP, which we briefly  illustrate   and compare with the classical PMP  through the discussion  of an elementary example  in Sect.\ \ref{Sect8}. 
\par
\medskip 
\section{Preliminaries}\label{Sect3}
 \subsection{Notational issues} \label{notation}\hfill\par
 Throughout this paper, we consider a dynamical system, whose states are represented by the 
 points of an appropriate $N$-dimensional  {\it manifold}  $\cQ$,  which might be for instance a configuration space or a  phase space for the system.  A generic set of local coordinates on $\cQ$ will be  usually denoted by  $(q^i)_{i = 1,  \ldots, N}$, so that the evolutions in time  of our system  correspond  to  parameterised curves with  coordinate expressions of the form
 $q(t) = (q^1(t), \ldots, q^N(t))$ for $t$ in a fixed interval  $I \subset \bR$. \par
Any such curve is  uniquely associated with the corresponding  parameterised graph
\beq \label{graph} \g:  I \subset \bR  \longrightarrow     \bR \times \cQ \ ,\qquad \g(t) = (t, q(t))\ ,\eeq
i.e.   the  associated local section of the trivial bundle $\bR \times \cQ$ over  $\bR$. Due to this,  from now on {\it we   identify any  evolution of our  system   with   the associated  map  $t  \mapsto \g(t) =  (t,q(t))$}.  \par
\par
For any   such a  map $\g:  I \to \bR \times \cQ$ of class $\cC^n$, $n \geq 1$,  we denote by    $j^n_{t_o}(\g)$ its $n$-th order jet at the time $t_o$. For a classical reference on jets, see for instance \cite{Sa} (\footnote{For convenience of the  reader,  it might be convenient to  briefly mention  what we are going  to adopt  as   definition of an $n$-th order  jet of a curve $\g$ in $\bR \times \cQ$ at a point  $t_o$. It   is  the equivalence class of all the curves of the form $\h:]-t_o- \d, t_o + \d[ \longrightarrow \bR \times \cQ$, whose components in some (thus, in any) set of coordinates have    values and    derivatives up to order $n$ at $t_o$ equal to those of $\g$ at such a point.}). The collection 
of all these $n$-th order jets has a  natural structure of smooth manifold and it is called the {\it jet bundle of order $n$} of the  (trivial) bundle $\bR \times \cQ$ over $\bR$. We denote it by $J^n(\cQ|\bR)$.  \par
For a fixed chart $(q^i): \cU \to \bR^N$ on some open set $\cU \subset \cQ$, we  may consider the map  which sends 
 each jet $j^n_{t}(\g)$ into the $\wt N$-tuple,  $\wt N \= N(n +1)+1$, 
$$j^n_t(\g) \longmapsto \left(t, \, q_{(0)}^i(t) \= q^i(t),\, q^i_{(1)}(t) \= \frac{ d q^i}{dt}\bigg|_t,\,  q^i_{(2)}(t) \= \frac{ d^2 q^i}{dt^2}\bigg|_t,\, \ldots, \, q^i_{(n)}(t) \= \frac{ d^n q^i}{dt^n}\bigg|_t\right)\ .$$
This map is well known to be  a (locally defined) system of coordinates for $J^n(\cQ | \bR)$.  In a short notation,  we  denote such coordinates by 
$$(t, q^i_{(\b)}) =  (t,  q^i_{(0)}, q^i_{(1)},  \ldots, q^i_{(n)} )\ $$ 
and we  call them {\it  canonical jet coordinates} determined by the  coordinates  $(q^i)$ of $\cQ$. \par
\smallskip
For any   $\g(t) = (t, q^i(t))$,   the   {\it  $n$-th order lift of $\g$}  is the  corresponding  curve of jets
$$\g^{(n)}: I\longrightarrow J^n(\cQ|\bR)\ ,\quad \g^{(n)}(t) \= j^n_t(\g) =  (t, q^i_{(\b)}(t)) = \left(t, q^i(t), \frac{d q^i}{dt}\bigg|_t, \dots,  \frac{d^n q^i}{dt^n}\bigg|_t\right)\ .$$
\par
We denote by $K \subset \bR^{M}$   a fixed set  of  real $M$-tuples  $u = (u^a)$.  In what follows, a  (continuous or $k$-differentiable) curve $t \mapsto u(t)$ with values in $ K$    plays the role of  a  {\it control} for our system.  
 We stress the fact  that, in the literature, the terms ``control'', ``control parameter''  or ``control value''  is usually adopted to refer just to  a single value of the curve $u(t)$,  not to  the  curve $u(\cdot)$ as a whole, in contrast with what   we  do  in this paper. We hope that this will not be a source of confusion. \par
We always  assume that   $K$   is   the closure of a bounded open subset of $\bR^M$ and that  the boundary $\p K$ is smooth. This  assumption is mainly made  for the sake of simplicity, since   most of our arguments can be  generalised to a large class  of more  general situations and  under  weaker regularity assumptions. 
\par
\smallskip
\subsection{Controlled Lagrangians and controlled Euler-Lagrange operators}\label{problems} \hfill\par
 Let us consider a    {\it (smooth) Lagrangian with controls}, that is a $\cC^\infty$  function 
 $$L =  L(t, q^i_{(\b)}, u^a) : J^n(\cQ|\bR) \times K \to \bR$$ 
  depending  on
 \begin{itemize}
  \item  the coordinates   $(t, q^i_{(0)}, q^i_{(1)}, \ldots, q^i_{(n)})$ of the   $n$-th jets in $J^n(\cQ|\bR)$, 
 \item  the parameters $u = (u^a)$ arbitrarily varying in $K \subset \bR^M$. 
 \end{itemize}
 If there is an integer $r$ such that $L$ is independent on all of the  jets coordinates   $q^i_{(\b)}$ with $\b \geq r+1$, we say  that {\it $r$ is the actual order of $L$}. 
 For example,  the function 
 $$L(t, q^i_{(\b)}, u^a)  = (\sum_{a= 1}^M u^a)\frac{1}{2} \sum_{i = 1}^N (q^i_{(1)})^2 -   \frac{1 }{2}\sum_{i =1}^N (q^i)^2 $$
 can be surely  considered as  a Lagrangian with controls on  any controlled jet space $J^n(\cQ|\bR) \times K$  of order   $n \geq 1$.  If  we decide  that our   working  ambient is $J^n(\cQ|\bR)$ for some $n$ which is strictly larger than $1$,   the property that  $L$  is   independent   of   the jet coordinates  $q^i_{(2)}$, $q^i_{(3)}$, \ldots,   is synthetically expressed  by  saying that   {\it  $r = 1$ is the actual order  of $L$}.\par
 \smallskip
 In our discussions, we will always assume that  the order $n$ of the controlled jet space $J^n(\cQ|\bR) \times K$ is  sufficiently    larger than
 the actual order $r$ of the considered Lagrangian  $L$. This is needed for letting  all of the operators considered in this paper (as, for instance, the Euler-Lagrange operator described below) to be meaningful. 
As we will shortly  see, we just  need  that the order  of the   jet space $J^n(\cQ|\bR)$   satisfies  the following inequality
\beq\label{crucialineq}   2 r + 1 \leq n\ ,\eeq
which, from now on,  we  constantly and tacitly assume to be  satisfied 
 (\footnote{This  assumption could have been safely  omitted  if we considered the infinite jet spaces $J^\infty(\cQ|\bR)$ in place of  finite order jet spaces $J^n(\cQ|\bR)$.  However this would  have  forced us to  work with    infinite-dimensional manifolds, a category  that, for simplicity of the exposition, we prefer   to leave  undisturbed.  }).\par
 \smallskip
The  {\it  controlled Euler-Lagrange operator}  is the  $N$-tuple  $E = (E_i)_{i = 1}^N$  of  differential operators, acting  on any controlled Lagrangian $L$ of actual order $r$,  
defined by 
 \beq \label{controlledEL} E_i(L) \= \frac{\p L}{\p q^i} + \sum_{\b= 1}^r (-1)^\b \left(\frac{d}{dt}\right)^\b\left(\frac{\p L}{\p q^i_{(\b)}}\right)\ ,\qquad i = 1, \ldots, N\ . \eeq
 Here, the symbol $\frac{d}{d t}$ denotes  the  {\it total derivative operator}, that is  the operator that transforms any  function $f: J^n(\cQ|\bR) \times K \to \bR$ {\it  of  actual order  $r' \leq n-1$}  into  the function    
\beq \label{totalderiv}
 \frac{d f}{dt}:J^n(\cQ|\bR) \times K \longrightarrow \bR\ ,\qquad\frac{d f}{dt}\bigg|_{(t, q^i_{(\b)})}\= 
\frac{\p f}{\p t}\bigg|_{(t, q^i_{(\b)})}  + \sum_{\smallmatrix 1 \leq j \leq N\\0 \leq \d \leq r' \endsmallmatrix}\frac{\p f}{\p q^j_{(\d)}} q^j_{(\d+1)}\bigg|_{(t, q^i_{(\b)})}  \ .
\eeq
 Such an operator is  called  ``total derivative''  simply because,   for any    curve    $ (\g(t), u(t))  \in (\bR \times \cQ \times K)^{[0,T]}$  with constant $u(t) \equiv u_o$, the evaluation of  $ \frac{d f}{dt}$  at the points of the  curve of jets $(\g^{(n)}(t), u(t) = u_o)$ is equal to 
\beq \label{3.2-1} \frac{df}{dt}\bigg|_{(\g^{(n)}(t), u_o)} = \frac{d }{dt} \left(f(\g^{(n)}(t), u_o)\right)\bigg|_t\qquad \text{for each} \ t \in [0, T]\ ,  
\eeq
i.e.  it  coincides with  the  derivative with respect to $t$ of the map  $t \mapsto f( \g^{(n)}(t), u_o)$.\par
\smallskip
Notice that:
\begin{itemize}[leftmargin = 10pt]
\item The total derivative  raises the actual order of a function of at most one unit and  the iterated 
total derivatives $\left(\frac{d}{dt}\right)^\b$,   $1 \leq  \b \leq r$, raise the actual orders of at most   $r$ units. 
This is one of the reasons why we assume \eqref{crucialineq}. Other reasons for this come from the fact that certain  arguments of the proof of  Lemma \ref{pcprop} below  work nicely  only if   $2r$  is actually {\it strictly} less than $n$. 
\item What we    call    ``controlled Euler-Lagrange operator''  is  {\it almost} the same of   the    Euler-Lagrange operator of  the classical  theory of variations. The only difference with respect to the usual one is that   the   operators  $E_a {=} \frac{\p} {\p u^a} + \sum_{\b= 1}^r (-1)^\b \left(\frac{d}{dt}\right)^\b\big(\frac{\p }{\p u^a_{(\b)}}\big)$, $1 \leq q \leq M$, 
corresponding to   the infinitesimal variations of the  coordinates  $u^a$, are here  missing.
\end{itemize}
 \par
  \medskip
 \section{Defining triples and generalised Mayer problems}\label{Sect4}
 \label{definingtriples}
We are now able  to  delineate  in detail  the particular class of  control problems, which  is the main object of study of this paper.  
As we already mentioned, the  dynamical systems we are dealing  with  evolve  according to  curves, whose  parameterised graphs are the curves  \eqref{graph}   in $\bR \times \cQ$.
The independent variable  $t$ (the ``time'') of these evolutions  is from now on always assumed to be varying  in a fixed interval $[0,T]$.  \par 
The  control problems on which we focus are those  given  by the following    ingredients.
 \begin{itemize}[itemsep=5pt, leftmargin=10pt]
\item A {\it set of  control parameters} $\cK$, i.e. a set   of    pairs $U = (u(t), \s)$,  in which:
 \begin{itemize}
\item[a)] the first element  $u(t)$ is a  
smooth curve 
$u: [0,T]   \to K$   in the above fixed ambient    space  $K \subset \bR^M$; 
\item[b)]  the second element   is 
a  jet  $\s = j^{2r-1}_{t = 0}(\g)$  where $r$ is the actual order of the controlled Lagrangian considered below; the jet $\s$ is   possibly subjected to some  constraints
(as, for instance,  that the $0$-th component  $ j^0_{t = 0}(\g) = \g(0)$ is equal to a fixed point $q_o \in \cQ$)  and  is  later  used as  the  initial datum for   a curve $\g(t) = (t, q(t))$,  described in the next  point.
\end{itemize}
\item  A {\it Lagrangian with controls} $L = L(t, q^i_{(\b)}, u^a)$  of  actual order $r$  satisfying \eqref{crucialineq}, which gives the   {\it system of   controlled  Euler-Lagrange  equations}  of order $2r$ 
 \beq \label{controlledELbis} E_i(L) |_{(j^n(\g(t)), u(t))} = 0  \eeq
 for each smooth curve $u:[0,T] \to K$. 
We  also assume that $L$  satisfies   all   needed regularity and  maximal rank conditions that    guarantee  the  following:  {\it for each  pair  $U = (u(t), \s)$ in the set  $ \cK$, 
 there exists  a unique  solution $\g^{(U)}(t) $   to the initial value problem  formed by the equations \eqref{controlledELbis} and the initial condition 
 \beq \label{initialvalue} j^{2r-1}_{t = 0}(\g(t)) = \s\ .\eeq}
The curves $\g^{(U)}$  determined in this way are  called {\rm $\cK$-controlled}.
\item A {\it  terminal cost function}, that is  a  real  function of the   jets at the time $t = T$ of some fixed  order $\wt r\leq n-1$
(\footnote{This condition on the actual order $\wt r$ is  just a  convenient  technical requirement and is used only in the proof of  Lemma \ref{pcprop} below. As observed before,  if we work in jets spaces of sufficiently high order, this requirement is always  easily  satisfied.}). We   assume that such a terminal cost function is actually    the restriction  $C|_{J^{n}(\cQ|\bR)|_{t = T} }$ of  a   smooth real function  of actual order $\wt r$ on the whole jet space  that    vanishes identically   on $ J^n(\cQ|\bR)|_{t = 0} $.  \end{itemize}
For any    terminal cost function on the jets  at $t = T$,    there are clearly   infinite   possibilities for  a smooth function   $C$ on $J^n(\cQ|\bR)$  that vanishes at the jets at $t = 0$ and that gives the desired cost function at $t = T$.  But in what follows   we    select   just   one of  such globally defined  functions and  we call  it  {\it  the    (extended) cost function} of  our problem.  
\par
\smallskip
Any triple  $(\cK, L, C)$,   formed by three ingredients of the above form,  is  called a   {\it defining triple}.   Given such a triple we may consider the following \par
 \begin{definition} \label{MPdefinition}  The   {\it  generalised Mayer problem  determined by  $(\cK, L, C)$}  is the problem  of looking for  all $\cK$-controlled evolutions  $\g^{(U)}$ for which the value of the terminal cost $C(j^n_{t = T}(\g^{(U)}))$ 
 is  minimal among  the   terminal  costs of all   $\cK$-controlled evolutions. For such curves, the corresponding pairs  $U = (u(t), \s) \in \cK$ are called {\it optimal controls}. 
\end{definition}
\par
\medskip
\section{Modified  triples and   Mayer problems in integral form}\label{Sect5}
\label{subsect31}
\subsection{Performing the step (a)  of   the road map} \label{encapsulating-bis-A}\hfill\par
 We now  proceed according to  the step  (a) described  in Sect.\ \ref{roadmap}.  More precisely, we  canonically associate with any given curve $(t, q(t))$  in  $[0,T] \times  \cQ$  a special set of real  functions. They are
\begin{align} \label{exprA-Abis}  & \pa^i_\b(t) =  A^i_\b e^t +  B^i_\b e^{-t}   \ ,\\
 \label{exprA'-Abis}  & \pa'{}^i_\b(t) =  A'{}^i_\b e^{\frac{\pi}{2 T} t} +  B'{}^i_\b e^{- \frac{\pi}{2 T}  t}   +  C'{}^i_\b \cos\left(\frac{\pi}{2 T}  t\right)+  D'{}^i_\b \sin \left(\frac{\pi}{2 T} t\right) \ ,\\
  \label{exprA''-Abis} & \pa''{}^i_\b(t) =  A''{}^i_\b e^{\frac{\pi}{2 T} t}  +  B''{}^i_\b e^{- \frac{\pi}{2 T}  t}   +  C''{}^i_\b \cos \left(\frac{\pi}{2 T}  t\right)+ D''{}^i_\b \sin \left(\frac{\pi}{2 T}  t\right) \ ,
 \end{align}
  where the  indices $i$ and $\b$ run between  $1 \leq i \leq  N$ and $0 \leq \b \leq r-1$ (here, $r$ is the actual order of $L$) and   the $ A^i_\b$,  $ B^i_\b$, $A'{}^i_\b$, $ B'{}^i_\b$ etc., are the constants that  are uniquely determined  by  the following conditions on   the initial and terminal data of the   $q^i(t)$. The  $ A^i_\b$ and   $ B^i_\b$ are determined by solving the linear equations
\beq \label{primacond-A}     A^i_\b +  B^i_\b  (=\pa^i_{\b}\big|_{t = 0} ) = q^i_{(\b)}\big|_{t = 0},\  A^i_\b - B^i_\b  (  = \frac{d \pa^i_{\b}}{ d t} \bigg|_{t = 0} \hskip-7 pt) {=}  -  \hskip-0.5cm \sum_{\smallmatrix 1 \leq \d \leq r\\[2pt] 0 \leq \ve \leq \d\\[2pt]
\d  - \ve -1 =  \b  \endsmallmatrix }\hskip-10pt  (-1)^{\ve}  \frac{d^\ve}{dt^\ve} \left(\frac{\p   L}{\p q^i_{(\d)}} \right)\!\!\bigg|_{j^{n-1}_{t = 0}(q(t))}\hskip - 8pt .
\eeq
The remaining  constants are  set to
\begin{align} \label{condition-1-A-11} &\left(\begin{array}{c} A'{}^i_\b \\ B'{}^i_\b \\ C'{}^i_\b \\ D'{}^i_\b  \end{array} \right) = \cA^{-1}\left(\begin{array}{c} 0 \\ 0 \\ q^i_{(\b)}(T) \\ \displaystyle{\sum_{\smallmatrix 1 \leq \d \leq r\\[2pt] 0 \leq \ve \leq \d\\
\d  - \ve -1 =  \b  \endsmallmatrix }  (-1)^{\ve}  \frac{d^\ve}{dt^\ve} \left(\frac{\p  L}{\p q^i_{(\d)}} + \frac{\p  }{\p q^i_{(\d)}} \left( \frac{d C}{dt}\right)\right) \bigg|_{j^{n-1}_{t = T}(q(t))}}  \end{array} \right)\ ,\\
 \label{condition-2-A}
&\left(\begin{array}{c} A''{}^i_\b \\ B''{}^i_\b \\ C''{}^i_\b \\ D''{}^i_\b  \end{array} \right) = \cA^{-1}\left(\begin{array}{c} 0 \\ 0 \\ \pa^i_{\b}(T) \\ \pa^i_{\b(1)}(T)  \end{array} \right)\ ,
\end{align}
where $\cA$ is the real matrix
 $$\cA = \left(\begin{array}{cccc} 1 & 1 & 1 & 0\\  \frac{\pi}{2 T} & - \frac{\pi}{2 T} &  0 & \frac{\pi}{2 T}\\
  \frac{\pi}{2 T} e^{\frac{\pi}{2}} & - \frac{\pi}{2 T}  e^{-\frac{\pi}{2}}& -  \frac{\pi}{2 T} & 0\\
  \left(  \frac{\pi}{2 T}\right)^2 e^{\frac{\pi}{2}} &  \left(\frac{\pi}{2 T} \right)^2 e^{-\frac{\pi}{2}}& 0 & -  \left(\frac{\pi}{2 T}\right)^2 
 \end{array}\right)\ .
$$
A tedious but straightforward check   shows that  the  functions \eqref{exprA-Abis}  -- \eqref{exprA''-Abis}    are  precisely  the  unique  solutions to  the system of differential equations 
 \beq  \label{modifiedEL-A1}  \frac{d^2 \pa^i_\b}{dt^2} -  \pa_\b^i = 0\ ,\quad      \frac{d^4\pa'{}^i_\b}{dt^4}  -  \left(\frac{\pi}{2 T}\right)^4 \pa'{}_\b^i = 0 \ ,   \quad   \frac{d^4\pa''{}^i_\b}{dt^4}  -   \left(\frac{\pi}{2 T}\right)^4 \pa''{}_\b^i = 0 
  \eeq 
 together with     the set of the  boundary conditions formed  by the  \eqref{primacond-A} and by 
 \begin{align}
 & \pa'{}^i_{\b} (0) = 0\ , && \frac{d \pa'{}^i_{\b}}{dt}\bigg|_{t = 0} = 0\ ,\\
\label{ahio-A} & \frac{d \pa'{}^i_{\b}}{dt}\bigg|_{t = T} = q^i_{(\b)}(T)\ , && \frac{d^2 \pa'{}^i_{\b}}{dt^2}\bigg|_{t = T}  \hskip -0.5cm =  \sum_{\smallmatrix 1 \leq \d \leq r\\[2pt] 0 \leq \ve \leq \d\\[2pt]
\d  - \ve -1 =  \b  \endsmallmatrix }  (-1)^{\ve}  \frac{d^\ve}{dt^\ve} \left(\frac{\p    \left(L + \frac{ dC}{dt}\right)}{\p q^i_{(\d)}} \right) \bigg|_{j^{n-1}_{t = T}(q(t))} \hskip -0.5cm  ,\\
\label{ahiobis-A}
& \pa''{}^i_{\b} (0) = 0\ ,\  \frac{d \pa''{}^i_{\b}}{dt}\bigg|_{t = 0} = 0\ ,&& \frac{d \pa''{}^i_{\b}}{d t}\bigg|_{t = T}  = \pa^i_{\b}(T)   \ , \   \frac{d^2 \pa''{}^i_{\b}}{dt^2}\bigg|_{t = T}   = \pa^i_{\b(1)} (T) \ .
\end{align}
As a matter of fact,  we consider these  functions just  because we want them to  satisfy   such a differential problem. 
The  motivation for this requirement  has  been very roughly indicated in our road map and  it will be   definitely  clarified in the next section. \par
Using the  functions  \eqref{exprA-Abis} -- \eqref{exprA''-Abis},  with  any given curve  $(t, q^i(t))$ in  $[0,T] \times \cQ$ we may associate  a  curve $(t, q^i(t),  \pa^i_{\b}(t), \pa'{}^i_{\b}(t), \pa''{}^i_{\b}(t))$
in the enlarged   manifold    
$ [0,T] \times \wt \cQ $,  $ \wt \cQ \=  \cQ\times  \bR^{3 N r}$.  
    Such  a bijective correspondence between  curves  in $[0,T] \times \cQ$ and in  $[0, T] \times \wt \cQ$  
 establishes a natural  equivalence between our original control problem,  determined by the   triple $(\cK,  L,  C)$,   and a new control problem,  determined by an appropriate    modified defining triple $(\wt \cK, \wt L, \wt C)$ which  is defined as follows.\\[10pt]
$\bullet$  $\wt \cK$ is the collection of pairs $\wt U = (u(t), \wt \s)$ in which: (a)  $u(t)$ is precisely as it    occurs  in  the pairs   $(u(t), \s) \in  \cK$ and (b) $\wt \s$ is a jet in $J^n(\wt \cQ|\bR)\big|_{t = 0}$ of actual order $2 \max\{ r, 2\} - 1$,  playing the role of  the initial datum of a curve 
$\wt \g(t) = (t, q^i(t), \pa^i_\b(t),  \pa'{}^i_\b(t), \pa''{}^i_\b(t))$,  constrained  by the following conditions: 
\begin{itemize}[leftmargin = 20pt]
\item[(i)] the part  of $\wt \s$, corresponding  to    the  initial datum   of   $\g(t) = (t, q^i(t))$, 
 satisfies   the same  constraints  that are  imposed on the  data $\s$  in the    pair  $(u(t), \s)$ in  $\cK$; 
\item[(ii)] the  initial values  for the curves 
$\pa^i_{\b}(t)$, $\pa'{}^i_\b(t)$ and   $\pa''{}^i_\b(t)$  are required  to satisfy  the conditions given by the   \eqref{primacond-A},  \eqref{ahio-A} and 
\eqref{ahiobis-A}; no other condition  is imposed besides those which are naturally requested in order to be  initial conditions that are fully   compatible with the Euler-Lagrange equations   \eqref{12} and  \eqref{13}   below (\footnote{As we will shortly see,   such Euler-Lagrange equations  are nothing but   the  \eqref{modifiedEL-A1}.}).
\end{itemize}
\ \\
$\bullet$   $\wt L$ is the controlled Lagrangian  of actual order $\wt r = \max\{ r, 2\}$ 
\begin{multline} \label{LcC-A} 
 \wt L(t, q^i_{(\d)}, \pa^j_{\b(\d)},   \ldots, u^a) \=  \
 L +\frac{1}{2}\sum_{\smallmatrix 1 \leq j \leq N\\
   0 \leq  \b \leq r-1
   \endsmallmatrix} \bigg( ( \pa^j_{\b(1)})^2 - ( \pa'{}^j_{\b(2)})^2 -  ( \pa''{}^j_{\b(2)})^2  \bigg) + \\
   +  \sum_{\smallmatrix 1 \leq j \leq N\\
   0 \leq  \b \leq r-1
   \endsmallmatrix} \bigg( \frac{1}{2} ( \pa^j_{\b})^2+ \frac{\pi^4}{32 T^4} ( \pa'{}^j_{\b})^2 +  \frac{\pi^4}{32 T^4}( \pa''{}^j_{\b})^2 \bigg) \ .  
\end{multline}
$\bullet$  $\wt C$ is equal to   $\wt C = C$.   The  only  difference  between $C$ and $\wt C$  is  just  that   its formal domain is  now  $J^n(\wt \cQ|\bR)$  (and no longer   $J^n(\cQ|\bR)$).
\par
\smallskip
 If we  now replace the   triple $(\cK, L,  C)$  by   $(\wt \cK, \wt L, \wt C =  C)$, we have a new  Euler-Lagrange operator,   consisting of the differential operators  $E_i(\cdot)$,   $E^{\b}_i(\cdot)$,     $E'{}^{\b}_i(\cdot)$,    $E''{}^{\b}_i(\cdot)$,      corresponding   to the  variables 
$q^i$,  $\pa^i_\b$,   
$\pa'{}^i_\b$ and  $\pa''{}^i_\b$, respectively.
 The new  set of  Euler-Lagrange equations is 
 \begin{align}  \label{11}
 &E(\wt  L)_i =       \frac{\p L}{\p q^i} + \sum_{\ell= 1}^r (-1)^\ell \frac{d^{\ell }}{dt^{\ell}}\left(\frac{\p L}{\p q^i_{(\ell)}}  \right)=   0\ ,\\
  \label{12}   & E(\wt L)^\b_i=  -  \frac{d^2 \pa^i_\b}{dt^2} +  \pa_\b^i = 0\ ,\quad    E(\wt L)'{}^\b_i=  - \frac{d^4\pa'{}^i_\b}{dt^4}  +  \left(\frac{\pi}{2 T}\right)^4 \pa'{}_\b^i = 0 \ ,  \\
   \label{13}  & E(\wt L)''{}^\b_i=  - \frac{d^4\pa''{}^i_\b}{dt^4}  +   \left(\frac{\pi}{2 T}\right)^4 \pa''{}_\b^i = 0  \ .
  \end{align}
From  these equations, we  directly see that      any $\wt \cK$-controlled curve $ \wt \g^{(\wt U)}(t) = (t, $ $q^i(t),$ $\pa^i_\b(t), \pa'{}^i_\b(t), \pa''{}^i_\b(t))$
  has the following two crucial properties.  
\begin{itemize}[leftmargin = 18pt]
\item[($1$)] Being solutions to  the same differential problem,  the   components $q^i(t)$ of a $\wt \cK$-controlled curve $\wt \g^{(\wt U)}$ and of the corresponding 
$\cK$-controlled curve $\g^{(U)}$ are the same. Therefore   the terminal costs $\wt C|_{j^n_{t = T}(\wt \g^{(\wt U)})}$ and   $C|_{j^n_{t = T}(\g^{(U)})}$  coincide. 
\item[($2$)] The functions $\pa_\b^i(t)$, $\pa'{}_\b^i(t)$, $\pa''{}_\b^i(t)$ have necessarily  the forms \eqref{exprA-Abis} -- \eqref{exprA''-Abis}.   
\end{itemize}
 Thus the  generalised Mayer problem defined by  $(\wt \cK, \wt L, \wt C)$ is perfectly  equivalent to the original one,  given  by  $( \cK,  L,  C)$. The bijection  between the  two families of  controlled curves is established by simply   considering   the functions  defined  in \eqref{exprA-Abis} --  \eqref{exprA''-Abis} as  the last  components of  a  curve  in $\bR \times \wt \cQ = (\bR \times  \cQ)\times  \bR^{3 N r}$.
\par
\smallskip
\subsection{Performing the step (b) of   the road map} \label{encapsulating}\hfill\par
 We now  go into  the   step (b)  of Sect.\ \ref{roadmap}. Namely,  we  introduce 
 two    auxiliary  variables    $\l$, $\mu$ and further  modify the  defining triple of the problem,  so that  an analogue of     \eqref{smart-1}  holds  for  any solution of  the new  controlled Euler-Lagrange equations.  
 \par
 Introducing two new   variables   corresponds to    enlarging  the manifold     $\wt \cQ = \cQ \times \bR^{3 Nr}$ of the previous section  into the  new  manifold $\wh \cQ = \wt  \cQ \times \bR^{2}$  parameterised by the    coordinates 
 $ (t, q^i, \pa_\b^i, \pa'{}_\b^i, \pa''{}_\b^i,   \l, \mu)$. After considering such new enlarged manifold $\wh \cQ$,  we have to introduce  the     {\it further modified  triple} $(\wh \cK, \wh L, \wh C)$ defined as follows. \\[10pt]
$\bullet$  $\wh \cK$ is the collection of pairs $\wh U = (u(t), \wh \s)$ in which: (a)  $u(t)$ is precisely as it    occurs  in  the pairs  in   $ \cK$ and $\wt \cK$ and (b) $\wh \s$ is a jet in $J^n(\wh \cQ|\bR)\big|_{t = 0}$ of actual order $2  \max\{ r, 2\} -1$,  playing the role of  the initial datum of a curve 
$\wh \g(t) = (t, q^i(t), \pa^i_\b(t),  \pa'{}^i_\b(t), \pa''{}^i_\b(t),  \l(t), \mu(t))$,  constrained  by the following conditions: 
\begin{itemize}[leftmargin = 20pt]
\item[(i)] the part  of $\wh \s$, corresponding  to    the  initial datum   of   $\wt \g(t) = (t, q^i(t),  \pa^i_\b(t),  \pa'{}^i_\b(t),$ $\pa''{}^i_\b(t))$, 
 satisfies   the same  conditions  that are  imposed on the pairs  $(u(t),\wt \s) \in \wt \cK$; 
\item[(ii)] the  component $\l|_{t = 0}$ of $\wh \s$ (i.e.     the initial value   of  $\l(t)$)  is always set to be     $\l|_{t = 0} = 1$;  the  values of  all other components 
of the jet  giving  the   initial datum for $\l(t)$   are   required to be just  compatible with  the Euler-Lagrange  equations   \eqref{modifiedEL} below (\footnote{By  looking at those equations,   one can  see  that this  compatibility requirement simply  means that  all  derivatives of $\l(t)$  at  $t = 0$ must be equal to $0$.}); 
\item[(iii)] the  component $\mu|_{t = 0}$ of $\wh \s$  (i.e.   the initial value  of $\mu(t)$)  must  be   $\mu|_{t = 0} = 0$;  the other components of that   initial datum  for $\mu(t)$ have   just  to be compatible with  the  \eqref{modifiedEL}.
\end{itemize}
$\bullet$  $\wh L$   is  the controlled Lagrangian  of actual order $\wh r =  \max\{ r, 2\}$   defined by 
\begin{multline}
 \wh L(t, q^i_{(\d)}, \pa^i_{\b(\d)},    \pa'{}^i_{\b(\d)}, \pa''{}^i_{\b(\d)} , \l_{(\d)}, \mu_{(\d)}, u^a) \=    \lambda \big (\mu_{(1)} + \wt L\big) +  \frac{d C}{dt} =   \\
 = 
   \lambda \Bigg(\mu_{(1)} + L +\frac{1}{2}\sum_{\smallmatrix 1 \leq i \leq N\\
   0 \leq  \b \leq r-1
   \endsmallmatrix} \bigg( ( \pa^i_{\b(1)})^2 - ( \pa'{}^i_{\b(2)})^2 -  ( \pa''{}^i_{\b(2)})^2  \bigg) + \\
   +  \frac{1}{2}  \sum_{\smallmatrix 1 \leq i \leq N\\
   0 \leq  \b \leq r-1
   \endsmallmatrix} \bigg(( \pa^i_{\b})^2+ \left(\frac{\pi}{2T}\right)^4 ( \pa'{}^i_{\b})^2 +  \left(\frac{\pi}{2T}\right)^4( \pa''{}^i_{\b})^2 \bigg) \Bigg) +  \frac{d C}{dt}\ .  
\end{multline}
Here, according to \eqref{totalderiv},  
$ \frac{dC}{dt}$  denotes   the  total derivative of the function $C = C(t, q^i_{(\b)})$.
\\[10pt]
 $\bullet$  $\wh C$ is    the same of the original cost function   $\wh C = C$. As  before,  the    difference  between $\wh C$ and $C$   is  just that  we are now  considering   it as a function on the new jet space  $J^n(\wh \cQ|\bR)$.
\par
\smallskip
 If we now  replace the  defining triple $(\wt \cK, \wt L, \wt C)$  of the previous section  by   $(\wh \cK, \wh L, \wh C = C)$, we have to consider  another   Euler-Lagrange operator,  formed not only  by  the  previous operators  $E_i(\cdot)$,    $E^{\b}_i(\cdot)$,     $E'{}^{\b}_i(\cdot)$,    $E''{}^{\b}_i(\cdot)$, but also by 
  the   operators   $E_{\{\l\}}(\cdot)$  and $E_{\{\mu\}}(\cdot)$, corresponding   to the new variables  $\l$ and $\mu$, respectively.
  The new  set of  Euler-Lagrange equations is: 
 \begin{align} 
\nonumber  &  E(\wh L)_{\{\mu\}}  = - \frac{ d \lambda}{dt} = 0\qquad \left( \ \Longrightarrow\ \ \l\equiv 1\ \ \right) ,\\
\nonumber  &E(\wh L)_i =   E(\lambda \wt L)_i =     \frac{\p L}{\p q^i} + \sum_{\ell= 1}^r (-1)^\ell \frac{d^{\ell }}{dt^{\ell}}\left(\frac{\p L}{\p q^i_{(\ell)}}  \right)=   0\ ,\\
\nonumber & E(\wh L)^\b_i=     E(\lambda \wt L)^\b_i = -  \frac{d^2 \pa^i_\b}{dt^2} +  \pa_\b^i = 0\ ,\\
\label{modifiedEL} & E(\wh L)'{}^\b_i=  E(\lambda \wt L)'{}^\b_i = - \frac{d^4\pa'{}^i_\b}{dt^4}  +  \left(\frac{\pi}{2 T}\right)^4 \pa'{}_\b^i = 0 \ ,  \\
 \nonumber  & E(\wh L)''{}^\b_i=   E(\lambda \wt L)''{}^\b_i =   - \frac{d^4\pa''{}^i_\b}{dt^4}  +   \left(\frac{\pi}{2 T}\right)^4 \pa''{}_\b^i = 0  \ ,\\
  \nonumber   &   E(\wh L)_{\{\l\}}  =\frac{d \mu}{dt} +\wt L= \\
\nonumber & \hskip 1.4 cm =  \frac{d \mu}{dt} +  L   + \frac{1}{2}\!\!\!\sum_{\smallmatrix 1 \leq i \leq N\\
   0 \leq  \b \leq r-1
   \endsmallmatrix} \!\!\! \bigg( ( \pa^i_{\b(1)})^2 {-} ( \pa'{}^i_{\b(2)})^2 {-}  ( \pa''{}^i_{\b(2)})^2 
    \bigg) + \\
\nonumber   & \hskip 1.8 cm + \!\!\! \sum_{\smallmatrix 1 \leq i \leq N\\
   0 \leq  \b \leq r-1
   \endsmallmatrix} \!\!\! \bigg( \frac{1}{2} ( \pa^i_{\b})^2+ \frac{\pi^4}{32 T^4} ( \pa'{}^i_{\b})^2 +  \frac{\pi^4}{32 T^4}( \pa''{}^i_{\b})^2  \bigg)     = 0\ ,
   \end{align}
By just looking at these equations, we  see that      any $\wh \cK$-controlled curve 
 \beq \label{wtK-controlled} \wh \g^{(\wh U)}(t) = (t, q^i(t), \pa^i_{\b}(t), \pa'{}^i_{\b}(t), \pa''{}^i_{\b}(t),  \l(t), \mu(t))\eeq
  has the following properties.  
\begin{itemize}[leftmargin = 15pt]
\item[($1$)] The value $\lambda(t)$ is constant and equal to the prescribed initial value, i.e. $\l(t) =1$.   It follows that     {\it the  new  differential constraints   on  the curve  $\wt \g(t) = (t, q^i(t),  \pa^i_\b(t),  \pa'{}^i_\b(t),$ $\pa''{}^i_\b(t))$ are identical  with  the  original  constraints  \eqref{11} -- \eqref{13}}. This, together with  the uniqueness of the solution to the differential problem for   $\mu(t)$,   imply  that there exists a natural bijection between the class of $\wt \cK$-controlled curves $\wt \g^{(\wt U)}$  and the class of the $\wh \cK$-controlled curves  $\wh \g^{(\wh U)}$.
\item[($2$)]  Due to (1), the   cost $\wh C|_{j^n_{t = T}(\wh \g^{(\wh U)})}$  for a $\wh \cK$-controlled curve  is always  equal to  the cost $\wt C|_{j^n_{t = T}(\wt \g^{(\wt U)})} = C|_{j^n_{t = T}(\g^{(U)})}$  of the corresponding $\wt \cK$-controlled curve $\wt \g^{(\wt U)}$.   
\item[($3$)] The  last  equation in  $\eqref{modifiedEL}$  implies that the $1$-form
$ \lambda (\mu_{(1)} + \wt L ) dt$
 vanishes identically along  any  curve of  jets  $\wh \g^{(\wh U) (n)}(t) = j^n_t\left(\wh \g^{(\wh U)}\right)$  of  a $\wh \cK$-controlled curve $\wh \g^{(\wh U)}$.
\end{itemize}
From  (1) and (2)  and previous discussion,  we  see  that {\it the  generalised Mayer problem defined by the triple  $(\cK, L, C)$ is  not only equivalent to the  problem of Sect.\ \ref{encapsulating-bis-A},   determined by  the triple $(\wt  \cK, \wt  L, \wt  C)$,  but also equivalent  to this new problem, determined by  the triple $(\wh \cK, \wh L, \wh C)$}. Furthermore (3)  shows that   the improvement, which was mentioned  in the step (b) of  the road map,  is now reached,  In fact, if we  
 consider   the $1$-form    
\beq\label{alpha} 
\begin{split} \a \=  & \wh L dt   = \Big( \lambda (\mu_{(1)} + \wt L  )+  \frac{d C}{dt} \Big)dt
\end{split} \eeq
we may observe   that,   for any  curve  of jets 
 $\wh \g^{(\wh U)(n)}(t)$ of a $\wh \cK$-controlled curve  $\wh \g^{(U)}$
 \beq \label{extendedaction}
\begin{split} \int_0^T \a\left(\frac{d \wh \g^{(\wh U)(n)}}{dt} \right)\bigg|_{\wh \g^{(\wh U)(n)}(t)} dt &=  \int_0^T \wh L|_{\wh \g^{(\wh U)(n)}(t)} dt \overset{\text{Property}\ (3)}=  \int_0^T \frac{dC}{dt}\bigg|_{\wh \g^{(\wh U)(n)}(t)} dt  \overset{ \eqref{3.2-1}} = \\
= C(\g^{(U)(n)}(T)) & - C(\g^{(U)(n)}(0)) \overset{ C|_{J^n(\cQ|\bR)|_{t = 0}} = 0} =   C(j^n_{t = T}(\g^{(U)}))  \ .
\end{split}
\eeq
  This means that   for each $\wh \cK$-controlled curve, the  integral \eqref{extendedaction}  is just equal  to    the  terminal cost $C|_{j^n_{t = T}(\g^{(U)})} $
  and that  {\it looking for  a solution to the  original  problem  is perfectly equivalent to looking for  a  $\wh \cK$-controlled evolution,  for which  the integral \eqref{extendedaction} 
 is  minimal among    those   of  all  other  $\wh \cK$-controlled evolutions.}  \par
\smallskip
\subsection{A convenient replacement of  the integrand in \eqref{extendedaction}}\hfill \par
\label{variationallyequiv}
We  now want to show that we may safely substitute  the $1$-form \eqref{alpha}   by a different one,  
which turns out to be much more convenient for our further developments. 
In order to introduce  such convenient replacement, we first need  to   recall that on the  jet bundle $J^n(\wh \cQ|\bR) $  there exist  an  important   class of  distinguished $1$-forms, namely  the family  of {\it  the  $1$-forms that   vanish identically   on the  tangent vectors of    the curves of jets   $\wh \g^{(n)}(t)$   of    the parameterised graphs $ \wh \g(t) = (t, q^i(t),  \pa^i_\b(t),  \pa'{}^i_\b(t),$ $\pa''{}^i_\b(t),  \l(t), \mu(t))$}. \par
  \smallskip
  It is  known   that    such  distinguished  $1$-forms are precisely those that are  pointwise  linear combinations of the $1$-forms 
\beq \label{holonomic}
\begin{split}
&\o^i_{(\d)} \= dq^i_{(\d)} -  q^i_{(\d+1)} dt\ ,\\
 &\varpi^i_{\b(\d)} \= d\pa^i_{\b(\d)} -  \pa^i_{\b(\d+1)} dt\ ,\quad \varpi'{}^i_{\b(\d)} \= d\pa'{}^i_{\b(\d)} -  \pa'{}^i_{\b(\d+1)} dt\ , \hskip 0.5 cm  i = 1, \ldots, N\ ,
\\[5pt]
&\varpi''{}^i_{\b(\d)} \= d\pa''{}^i_{\b(\d)} -  \pa''{}^i_{\b(\d+1)} dt\ ,\hskip 5 cm  \d = 0, \ldots, n -1\ , \\
&\varpi^{\{\l\}}_{(\d)} \=  d\l_{(\d)} - \l_{(\d+1)} dt\ ,\qquad \varpi^{\{\mu\}}_{(\d)} \=  d\mu_{(\d)} - \mu_{(\d+1)} dt\ .
\end{split}
\eeq
Since these $1$-forms vanish identically on the tangent vectors of the  curves $\wh \g^{(n)}(t)$,  the value of the integral \eqref{extendedaction}
 does not changes if  $\a$  is replaced by  any other   $1$-form  
\beq \label{vareq}  \a' = \a + \bP^{(\d)}_{i} \o^i_{(\d)} +  \bQ^{(\d)}{}^{\b}_{i}\varpi^i_{\b(\d)} + \bQ^{'(\d)}{}^{\b}_{i}  \varpi'{}^i_{\b(\d)}+ \bQ^{''(\d)}{}^{\b}_{i} \varpi''{}^i_{\b(\d)} +
\bL^{(\d)} \varpi^{\{\l\}}_{(\d)} +\bM^{(\d)} \varpi^{\{\mu\}}_{(\d)}  \eeq 
for  some arbitrary   choices of  smooth  functions $ \bP^{(\d)}_{i},  \bQ^{(\d)}{}^{\b}_{i}$, $ \bQ^{'(\d)}{}^{\b}_{i}$, $\bQ^{''(\d)}{}^{\b}_{i}$, $\bL^{(\d)}$, $\bM^{(\d)}$ of  the points of $J^n(\wh \cQ|\bR) \times K$.   
Following  the terminology used in  \cite{Sp},   we say that  any such  $\a'$  is     {\it variationally equivalent} to $\a$.  
The  invariance   of  \eqref{extendedaction} under    replacements  with variationally equivalent $1$-forms might  be considered    as a  sort of ``invariance under gauge transformations''   of the cost  functional.    \par
\medskip
A particular  choice for the  $ \bP^{(\d)}_{i},  \bQ^{(\d)}{}^{\b}_{i}$, $ \bQ^{'(\d)}{}^{\b}_{i}$, etc. yields to  the  following $1$-form. \par
\begin{definition} The   {\it  controlled Poincar\'e-Cartan form associated with  $\wh L$ }  is the $1$-form  on $J^n(\wh \cQ|\bR) \times K$ defined by 
\beq \label{310} \begin{split}  \a^{PC} &=   \wh L dt +  \sum_{\d = 1}^r  \sum_{\ve = 0}^{\d-1} (-1)^{\ve} \frac{d^\ve}{dt^\ve} \left(\frac{\p \wh L}{\p q^i_{(\d)}} \right) \o^i_{(\d-(\ve+1))}  + \\
& \hskip 0.5 cm +  \l\hskip-10pt \sum_{\smallmatrix 1 \leq i \leq N\\
   0 \leq  \b \leq r-1
   \endsmallmatrix} \bigg( \pa^i_{\b(1)} \varpi^i_{\b(0)}  -  \pa'{}^i_{\b(2)} \varpi'{}^i_{\b(1)} -  \pa''{}^i_{\b(2)} \varpi''{}^i_{\b(1)}  \bigg)\ +  \l \varpi^{\{\mu\}}_{(0)} = \\
&=   \Bigg( \lambda \Big(\mu_{(1)} +  L +\frac{1}{2}\sum_{\smallmatrix 1 \leq i \leq N\\
   0 \leq  \b \leq r-1
   \endsmallmatrix} \bigg( ( \pa^i_{\b(1)})^2 - ( \pa'{}^i_{\b(2)})^2 -  ( \pa''{}^i_{\b(2)})^2  \bigg) + \\
 & \hskip 0.8cm  +  \sum_{\smallmatrix 1 \leq i \leq N\\
   0 \leq  \b \leq r-1
   \endsmallmatrix} \bigg( \frac{1}{2} ( \pa^i_{\b})^2+ \frac{\pi^4}{32 T^4} ( \pa'{}^i_{\b})^2 +  \frac{\pi^4}{32 T^4}( \pa''{}^i_{\b})^2 \bigg) \Big) +  \frac{d C}{dt} \Bigg) dt + \\
& \hskip 1.7 cm +\l \sum_{\d = 1}^r \sum_{\ve = 0}^{\d-1} (-1)^{\ve}  \frac{d^\ve}{dt^\ve}  \left(\frac{\p  \left(L + \frac{ dC}{dt}\right) }{\p q^i_{(\d)}} \right) \o^i_{(\d-(\ve+1))} + \\
 & \hskip 1.7 cm +   \l\hskip-10pt  \sum_{\smallmatrix 1 \leq i \leq N\\
   0 \leq  \b \leq r-1
   \endsmallmatrix} \bigg( \pa^i_{\b(1)} \varpi^i_{\b(0)}  -  \pa'{}^i_{\b(2)} \varpi'{}^i_{\b(1)} -  \pa''{}^i_{\b(2)} \varpi''{}^i_{\b(1)}+ \\
   & \hskip 4.7 cm +    \pa'{}^i_{\b(3)} \varpi'{}^i_{\b(0)} +  \pa''{}^i_{\b(3)} \varpi''{}^i_{\b(0)}\bigg) + 
\l \varpi^{\{\mu\}}_{(0)} \ .
\end{split}
\eeq
\end{definition}
Since  \eqref{310}  is variationally equivalent to $\a$, we may safely  replace   $\a$  by $\a^{PC}$  in \eqref{extendedaction}. And, in fact,  such a substitution is  the  analogue of  what is done in  \cite{CS}  for the  classical Mayer problems,  where  the $1$-form  
$$\a=  \left(\sum_{i = 1} p_i(\dot x^i - f^i(t, x, u(t)))  + \sum_{i = 1}^N \frac{\p C}{\p x^i }\bigg|_{(t, x(t))}  \dot x^i(t)\right) dt$$
 is  replaced by the    $1$-form
 \begin{multline} \label{vecchiaPC} \a^{PC} = \left( \sum_{i = 1} p_i(\dot x^i - f^i(t, x, u(t))) + \sum_{i = 1}^N \frac{\p C}{\p x^i }\bigg|_{(t, x(t))}  \dot x^i(t)\right) dt + \sum_{i = 1} p_i( d x^i -  \dot x^i dt)  = \\
 =  \sum_{i = 1} p_i dx^i - p_i f^i(t, x, u(t)) dt  + \sum_{i = 1}^N \frac{\p C}{\p x^i }\bigg|_{(t, x(t))}  \dot x^i(t)  dt\ .  \end{multline} 
The convenience of considering the  $1$-form \eqref{vecchiaPC} and, in more general situations,  the controlled Poincar\'e-Cartan forms \eqref{310} comes   from   the special   feature that  is described  in the next lemma and which will be  exploited in the proof of our first main result.  
\begin{lem} \label{pcprop} The differential  $d \a^{PC}$ of the controlled Poincar\'e-Cartan form  has the form 
\beq \label{39} \begin{split}
 d & \a^{PC} {=}  E(\wh L)_i  \o^i_{(0)} \wedge dt +\\
 &  +  E(\wh L)^\b_i \varpi^i_{\b(0)}\wedge dt +  E(\wh L)'{}^\b_i \varpi'{}^i_{\b(0)} \wedge dt +    E(\wh L)''{}^\b_i \varpi''{}^i_{\b(0)} \wedge dt +\\
 &+
  E(\wh L)_{\{\l\}}  \varpi^{\{\l\}}_{(0)} \wedge dt +E(\wh L)_{\{\mu\}}   \varpi^{\{\mu\}}_{(0)} \wedge dt+  \frac{\p \wh L}{\p u^a} d u^a \wedge dt  +\\
&+   \text{\rm  linear combinations of   wedges  of pairs of   
  $1$-forms   of  the list}\ \eqref{holonomic}
\end{split}\eeq
\end{lem}
\begin{pf}  The claim   is  an immediate  consequence of a   general fact concerning  the   classes of variationally equivalent  $1$-forms on   jets spaces (see e.g.  \cite[Prop. A2]{Sp}). For reader's convenience, we present here  a direct  proof.  For   simplicity of notation,  from now on we   denote  any tuple of coordinates   $(q^i, \pa^i_{\b},  \pa'{}^i_{\b}, \pa''{}^i_{\b}, \l, \mu)$ for  $\wh \cQ$   just by  $y = (y^\ell)$,  where the  index  $\ell$ ranges  between  $1$  and $\wh N =  N (3r+1) +2 $.  Accordingly, the   coordinates 
of the whole  jet space are  denoted by $t$ and by   $y^\ell_{(\d)}$, $0 \leq \d \leq n$,  and the associated  $1$-forms of the list  \eqref{holonomic}  are  briefly indicated     as    
\beq \label{list}  \o^\ell_{(\d)} = d y^\ell_{(\d)} - y^\ell_{(\d +1)} d t\qquad\text{with}\ \  1 \leq \ell \leq \wh N\ \ \text{and}\ \  1 \leq \d \leq n-1\ .\eeq
Note that, being $dt \wedge dt = 0$ and $d(d y^\ell_{(\d)}) = 0$,  for each integer $1 \leq \d \leq n-1$,    
\beq d y^\ell_{(\d)}\wedge dt =  \o^\ell_{(\d)} \wedge dt = - d\o^\ell_{(\d-1)} \ .\eeq
Hence for any function $f: J^n(\wh \cQ|\bR) \times K \to \bR$ of    actual order $\wh r \leq n-1$ 
\begin{multline} \label{lala}  f d y^\ell_{(\d)} \wedge dt  =  - f d \o^\ell_{(\d-1)}  
= \\
= - d \left(f   \o^\ell_{(\d-1)} \right)  + \frac{\p f}{\p t} dt \wedge  \o^\ell_{(\d-1)}+ 
\sum_{\d' = 0}^{\wh r} \frac{\p f}{\p y^m_{(\d')}} d y^m_{(\d')} \wedge \o^\ell_{(\d-1)} = \\
=  - d \left(f   \o^\ell_{(\d-1)} \right)  + \frac{\p f}{\p t} dt \wedge\o^\ell_{(\d-1)} + 
\sum_{\d' = 0}^{\wh r} \frac{\p f}{\p y^m_{(\d')}} \o^m_{(\d')} \wedge \o^\ell_{(\d-1)}  
+\sum_{\d' = 0}^{\wh r} \frac{\p f}{\p y^m_{(\d')}}  y^m_{(\d')}   dt \wedge  \o^\ell_{(\d-1)} = \\
= - d \left(f  \o^\ell_{(\d-1)} \right)   - \frac{d f}{dt}  \, dy^\ell_{(\d-1)}  \wedge dt  \qquad    \text{\rm  modulo  terms of the form }\  \o^r_{(\h)} \wedge \o^{s}_{(\z)}
\end{multline}
If $\d - 1\geq 1$, we may iterate and apply  this    identity  to the term   $-  \frac{d f}{dt}\, dy^\ell_{(\d-1)}  \wedge dt $   which appear in  the  right hand side of such identity. In this way we  get   that 
\beq   f d y^\ell_{(\d)} \wedge dt  =   - d \left(f  \o^\ell_{(\d-1)} \right)   +  d\left(  \frac{d f}{dt}  \o^\ell_{(\d-2)}\right) + 
\frac{d^2 f }{dt^2} \,d y^\ell_{(\d-2)}  \wedge dt  \ \    \text{\rm  mod}\  \o^r_{(\h)} \wedge \o^{s}_{(\z)}.
\eeq
If $\d -2 \geq 1$,  we may again apply   \eqref{lala}  to the term  $\frac{d^2 f }{dt^2} \,y^\ell_{(\d-2)}  \wedge dt$ and so on.  After  $\d$ iterations of such  use of   \eqref{lala},  we end up with 
\beq   \label{5.9} f d y^\ell_{(\d)} \wedge dt  =  -  d \left( \sum_{\ve = 0}^{\d-1}(-1)^\ve \frac{ d^\ve f}{dt^\ve}  \o^\ell_{(\d-1 -\ve)} \right)   + (-1)^\d
\frac{d^\d f }{dt^\d} \,\o^\ell_{(0)}  \wedge dt  \quad   \text{\rm  mod}\  \o^r_{(\h)} \wedge \o^{s}_{(\z)}
\eeq
 Applying  \eqref{5.9} to each  term $ \frac{\p \wh L}{\p y^\ell_\d} dy^\ell_{(\d)} \wedge dt $,   $\d \geq 1$,  appearing in  the exterior differential $d(\wh L dt )$, we obtain 
\beq \label{5.10}
\begin{split}
d(\wh L dt) &=  \frac{\p \wh L}{\p y^\ell_{(0)} } dy^\ell_{(0)} \wedge dt  + \sum_{\d = 1}^r \frac{\p \wh L}{\p y^\ell_{(\d)}} dy^\ell_{(\d)} \wedge dt = \\
&=  \frac{\p \wh L}{\p y^\ell_{(0)} } \o^\ell_{(0)} \wedge dt  -  \sum_{\d = 1}^r d \left( \sum_{\ve = 0}^{\d-1}(-1)^\ve \frac{ d^\ve}{dt^\ve}  \frac{\p \wh L}{\p y^\ell_{(\d)}} \o^\ell_{(\d-1 -\ve)} \right) +     \\
&\hskip 3 cm + \sum_{\d = 1}^r   (-1)^{\d+1}
\frac{d^\d }{dt^\d} \frac{\p \wh L}{\p y^\ell_{(\d)}} \,\o^\ell_{0}  \wedge dt \ \ \     \text{\rm  mod}\  \o^r_{(\h)} \wedge \o^{s}_{(\z)} = \\
 & =  - d \left(\sum_{\d = 1}^r   \sum_{\ve = 0}^{\d-1} (-1)^\ve \frac{ d^\ve}{dt^\ve}  \frac{\p \wh L}{\p y^\ell_{(\d)}} \o^\ell_{(\d-(\ve+1))} \right)+  E(\wh L)_\ell \, \o^\ell_{(0)} \wedge dt  \ \    \text{\rm  mod}\  \o^r_{(\h)} \wedge \o^{s}_{(\z)}
 \end{split}
\eeq
Note that,  in the  simplified notation  used in this proof,  the  $1$-form $\a^{PC}$  is nothing but   
$$\a^{PC} =  \wh L dt +  \sum_{\d = 1}^r  \left( \sum_{\ve = 0}^{\d-1}(-1)^\ve \frac{ d^\ve}{dt^\ve}  \frac{\p \wh L}{\p y^\ell_{(\d)}} \o^\ell_{(\d-(\ve+1))} \right) \ .$$
From this and   \eqref{5.10}, the  lemma follows immediately.
\end{pf}
\par
\medskip
 \section{The Principle of Minimal Labour}
 \label{section5}\label{Sect6}
\label{secthomotopy}
We are now ready to prove our first main result, the Principle of Minimal Labour for generalised Mayer problems. As we mentioned in Sect.\ \ref{roadmap}, this is reached by: (i)  first  proving  a generalisation of  the identity  \eqref{astast} for the  problem associated with the modified   defining triple  $(\wh \cK, \wh L, \wh C)$  and then  (ii)  deriving   a corresponding  identity  for the original    problem, determined by   $(\cK, L, C)$. These two identities   are  proven in Sect.\ \ref{section7.1} and Sect.\ \ref{section7.2}, respectively. \par
\smallskip
\subsection{The homotopy formula for the  problem defined by the  triple $(\wh \cK, \wh L, \wh C)$}\hfill\par
\label{section7.1}
  Let  $\wh U_o  = (u_o(t), \wh \s_o(t))$ be a fixed element in  $\wh \cK$ and  $\wh  \g^{(\wh U_o)}: [0,T] \to [0,T] \times \wh \cQ $ the corresponding $\wh \cK$-controlled curve.  We call  a smooth  $1$-parameter family ${\wh{F}}(\cdot, s)$, $s \in [0,1]$,  of  $\wh \cK$-controlled curves with  initial  curve ${\wh{F}}(\cdot, 0) = \wh  \g^{(\wh U_o)}$  a {\it $\wh \cK$-controlled variation  of}  $\wh  \g^{(\wh U_o)}$.  More precisely,     a  
  $\wh \cK$-controlled variation ${\wh{F}}$   is a smooth homotopy of the form 
  \beq \label{homotopy} {\wh{F}}:[0,T] \times [0,1] \longrightarrow \wh \cQ \ ,\qquad {\wh{F}}(t,s) = \wh \g^{(\wh U(s))}(t)\ ,\eeq
  where   $\wh U(s) = (u(\cdot,s), \wh \s(s))$, $s \in [0,1]$,  is a smooth curve  in $\wh \cK$  starting from    $\wh U(0) =\wh  U_o$.
 For each  ${\wh{F}}$, we consider  the    corresponding 
 homotopy  ${\wh{F}}^{(n)}$ in  $J^n(\wh \cQ|\bR) \times K$ defined by 
  $${\wh{F}}^{(n)}: [0,T] \times [0,1] \longrightarrow J^n(\wh \cQ|\bR)\times K\ ,\qquad {\wh{F}}^{(n)}(t,s) \= (j^{n}_t({\wh{F}}(\cdot, s)), u(t,s))\ .$$
 Note that, for each fixed $s_o$, the curve $t \mapsto j^{n}_t({\wh{F}}(\cdot, s_o))$ is nothing but the curve $\wh \g^{(\wh U(s_o)) (n)}$ of the $n$-th order jets of the  curve  
 $\wh \g^{(\wh U(s_o))}$. \par
 Given a $\wh \cK$-controlled variation  ${\wh{F}}$,   we  denote by $X_{\wh{F}}$ and $Y_{\wh{F}}$ the    vector fields --   defined just  at  the points  of the  surface  
 $\wh \cS \= {\wh{F}}^{(n)}([0,T] \times [0,1]) $ --  which  are determined   by considering  the infinitesimal variations of      the first parameter $t$ and of   the second parameter $s$, respectively.  More precisely, $X_{\wh{F}}$ and $Y_{\wh{F}}$ are the vector fields at the points of $\wh \cS$  
  \begin{align}\label{vectorX} X_{\wh{F}}|_{{\wh{F}}^{(n)}(t,s)} &\={\wh{F}}^{(n)}_*\left(\frac{\p}{\p t}\bigg|_{(t,s)}\right) = \frac{\p{\wh{F}}^{(n)}}{\p t}\bigg|_{(t,s)}\ ,  \\
\label{vectorY} Y_{\wh{F}}|_{{\wh{F}}^{(n)}(t,s)} & \={\wh{F}}^{(n)}_*\left(\frac{\p}{\p s}\bigg|_{(t,s)}\right) = \frac{\p{\wh{F}}^{(n)}}{\p s}\bigg|_{(t,s)}\ .
\end{align}
 We remark  that, for each fixed  $s_o \in [0,1]$, 
 \begin{itemize}[leftmargin = 20pt]
\item[(1)] The restriction of  $X_{\wh{F}}$ to the trace of   the  curve   $$ t \mapsto  {\wh{F}}^{(n)}(t,s_o)  = \left( \wh \g^{(\wh U(s_o))(n)}(t), u(t, s_o)\right)= \left(j^n_{t}{\wh{F}}(\cdot,  s_o), u(t, s_o)\right)$$
 coincides with  the family of    the tangent vectors of such a curve. 
 \item[(2)] The restriction of  $Y_{\wh{F}}$ to the trace of the same  curve   is  the Jacobi vector field   corresponding to the  (infinitesimal) variation $\ve \to{\wh{F}}^{(n)}(\cdot,s_o + \ve) $  of ${\wh{F}}^{(n)}(\cdot,s_o)  $. 
 \end{itemize}
 From (1), (2) and  \eqref{modifiedEL}, we  have  that  the vector fields $X_{\wh{F}}$ and   $Y_{\wh{F}}$ must  have  the  form  
 \beq \begin{split}  \label{54uno}
& X_{\wh{F}} = \frac{\p}{\p t} +   X^i_{(\b)} \frac{\p}{\p q^i_{(\b)}} + 
 X^{\{\mu\}}_{ (\b)} \frac{\p}{\p  \mu_{ (\b)} } + X^a \frac{\p}{\p u^a} +\\
 & \hskip 3.3 cm + X^i_{\a(\b)} \frac{\p}{\p \pa^i_{\a(\b)}}   + X'{}^i_{\a(\b)} \frac{\p}{\p \pa'{}^i_{\a(\b)}} + X''{}^i_{\a(\b)} \frac{\p}{\p \pa''{}^i_{\a(\b)}}    \ ,\\
 & Y_{\wh{F}} =    Y^i_{(\b)}  \frac{\p}{\p q^i_{(\b)}}  +   Y^{\{\mu\}}_{(\b)} \frac{\p}{\p  \mu_{ (\b)} } + Y^a \frac{\p}{\p u^a} + \\
  &\hskip 3.3cm  + Y^i_{\a(\b)} \frac{\p}{\p \pa^i_{\a(\b)}}   + Y'{}^i_{\a(\b)} \frac{\p}{\p \pa'{}^i_{\a(\b)}} + Y''{}^i_{\a(\b)} \frac{\p}{\p \pa''{}^i_{\a(\b)}}   
    \end{split}
  \eeq
 for  appropriate smooth real functions $X^i_{(\b)}$, $X^{\{\mu\}}_{(\b)}$,  etc., defined {\it only} at the points of $\wh \cS$. \par
 \begin{theo}[Homotopy Formula - First Version]  \label{prop34}  Let $\wh U_0, \wh U_1 \in \wh \cK$  be the endpoints of a smooth curve $\wh U(s) \in \wh \cK$, $s \in [0,1]$, and 
   $\wh \g \= \wh \g^{(\wh U_0)}$,  $\wh \g' \= \wh \g^{(\wh U_1)}$ the $\wh \cK$-controlled curves corresponding to  $\wh U_0, \wh U_1$, with terminal costs  $C_0 \= C(j^{n-1}_{t = T}(\wh \g))$  and  $C_1 \= C(j^{n-1}_{t = T}(\wh \g'))$, respectively.  Furthermore
\begin{itemize}
\item[(i)]  for any jet $ j^n_t(\wh \g) \in J^n(\wh \cQ|\bR)$,  let   $\cP_{j^n_t(\wh \g)}$ be    the  function on $K$ defined by 
\beq\label{pontr-funct} 
 \cP_{j^n_t(\wh \g)} : K \longrightarrow \bR \ ,\qquad  \cP_{j^n_t(\wh \g)}(u^a) \= -   L(j^n_t(\wh \g), u^a)\ ;\eeq
\item[(ii)] let $\wh \mu: [0,T] \times [0,1] \to \bR$ be the  function  defined by  $\wh \mu(t,s) \= \mu^{(\wh U(s))}(t)$, where $ \mu^{(\wh U(s))}(t)$ is  the value at $t$ of the $\mu$-component of the $\wh \cK$-controlled curve
$$\wh \g^{(\wh U(s))}(t) = (t, q^{(\wh U(s))i}(t),  \pa^{(\wh U(s))i}_\b(t), \ldots,  \l^{(\wh U(s))}(t) = 1,  \mu^{(\wh U(s))}(t))\ .$$
\end{itemize}
 Then,
  \begin{multline} \label{homotopy-ultimate} 
 C_1 - C_0   =  - \int_0^T \left( \int_0^1    \left(Y^a \frac{\p \cP_{\wh \g^{(\wh U^(s))(n)}_t} }{\p u^a} \Bigg|_{u^{(s)}(t)} - \frac{\p^2\wh  \mu}{\p t \,\p s}\bigg|_{(t,s)} -\right.\right.\\
\left.\left.- \sum_{\b = 1}^{r-1}  \sum_{i= 1}^N  \frac{\p^2}{\p t \p s}  \bigg( \int_0^s \pa'{}^i_{\b(3)} Y'{}^i_{\b(0)} + \pa''{}^i_{\b(3)} Y''{}^i_{\b(0)}  \bigg) \bigg|_{{\wh{F}}^{(n)}(t, v)} \right)dv \right) dt
\end{multline}
where  $Y^a$, $ Y'{}^i_{\b(0)}$,  $Y''{}^i_{\b(0)}$ are the  $u^a$-, $\pa'{}^i_\b$- and $\pa''{}^i_\b$- components, respectively,  of the vector field $Y_{\wh F}$ defined in \eqref{54uno}, associated with the $\wh \cK$-controlled variation $\wh F$ determined by the  $\wh U(s)$, $s \in [0,1]$, 
\end{theo}
\begin{rem} \label{rem72}  From  the  definition of the functions $\pa^i_{\b}(t)$, $\pa'{}^i_{\b}(t)$ and $\pa''{}^i_{\b}(t)$ and  the 
 Euler-Lagrange equation $E(\wh L)_{\{\l\}} = 0$, it follows immediately   that the second summand 
 in the right hand side of  \eqref{homotopy-ultimate}   is uniquely  determined by the  $\cK$-controlled curves $t \mapsto \g^{(U(s))}(t) =  (t, q^{(U(s))i}(t))$   in   $[0, T] \times \cQ$.  It is also simple to check that the same is true  for the first summand as well. These two  facts will be used in the next subsection. 
\end{rem}
\begin{pf}
 Consider the  embedded surface 
 $\wh \cS \= {\wh{F}}^{(n)}([0,T] \times [0,1])$  and the vector fields  $X_{\wh{F}}$, $Y_{\wh{F}}$ defined in  \eqref{vectorX} and \eqref{vectorY}. 
From  \eqref{extendedaction}, the  property (1) of  $X_{\wh{F}}$    and the fact that $\a^{PC}$ is a $1$-form which is variationally equivalent  to $\a$,  we have that for each curve $\wh \g^{(\wh U(s))}$, $s \in [0,1]$,
\beq \nonumber
\begin{split}
\int_0^T\left({\wh{F}}^{(n)*}(\a^{PC})\right)\bigg|_{(t,s)} \left( \frac{\p}{\p t}\right) dt  &=   \int_0^T \a^{PC}\left(X_{\wh{F}}\right)|_{{\wh{F}}^{(n)}(t, s)}  dt  = \\
& = \int_0^T  \a\left(X_{\wh{F}}\right)|_{{\wh{F}}^{(n)}(t, s)}  dt  = C(j^n_{t = T}(\wh \g^{(\wh U(s))}))  \ .
\end{split}\ .\eeq
This implies that 
\beq  \label{primissima} 
\begin{split} C_0 - C_1 &= C(j^n_{t = T}(\wh \g^{(\wh U(0))})) - C(j^n_{t = T}(\wh \g^{(\wh U(1))})) = \\
& = \int_0^T \left({\wh{F}}^{(n)*}(\a^{PC})\right) \left(\frac{\p}{\p t}\right) \bigg|_{(t,0)} dt  -
\int_0^T \left({\wh{F}}^{(n)*}(\a^{PC})\right)\left(\frac{\p}{\p t}\right) \bigg|_{(t,1)}  dt \ .
\end{split}
\eeq
On the other hand, by property (2) of $Y_{\wh{F}}$ and   the assumptions on the initial data of the  $\pa^i_{\b}(t)$, $\pa'{}^i_{\b}(t)$ and $\pa''{}^i_{\b}(t)$,  described in  Sect.\ \ref{encapsulating-bis-A},   we have that 
$$Y^{\{\mu\}}_{(0)}|_{{\wh{F}}^{(n)}(0,s)} =  Y'{}^i_{\b(0)}|_{{\wh{F}}^{(n)}(0,s)} = Y''{}^i_{\b(0)}|_{{\wh{F}}^{(n)}(0,s)}  =   Y'{}^i_{\b(1)}|_{{\wh{F}}^{(n)}(0,s)} = Y''{}^i_{\b(1)}|_{{\wh{F}}^{(n)}(0,s)}  =    0\ ,$$
$$  Y^i_{\b(0)}|_{{\wh{F}}^{(n)}(0,s)} =  Y^i_{(\b)}|_{{\wh{F}}^{(n)}(0,s)}\ . $$
From this we obtain  
\begin{multline} \label{secondina} \int_0^1 \imath_{\frac{\p}{\p s}}\left({\wh{F}}^{(n)*}(\a^{PC})\right)\bigg|_{(0,s)} ds = \int_0^1 \a^{PC}\left(Y_{\wh{F}}\right)|_{{\wh{F}}^{(n)}(0, s)}  ds = \\
=\sum_{\b = 1}^{r-1}   \sum_{i= 1}^N  \int_0^1  \bigg(
\bigg( \sum_{\smallmatrix 1 \leq \d \leq r\\[2pt] 0 \leq \ve \leq \d\\[2pt]
\d  - \ve -1 =  \b  \endsmallmatrix }  \hskip-0.5 cm (-1)^{\ve}  \frac{d^\ve}{dt^\ve} \left(\frac{\p  }{\p q^i_{(\d)}}   \left(L + \frac{ dC}{dt}\right) \right)  \bigg)Y^i_{(\b)}  +   \pa^i_{\b(1)} Y^i_{\b(0)}  \bigg)\bigg|_{{\wh{F}}^{(n)}(0, s)} 
 \hskip-0.5 cm ds  = \\
=\sum_{\b = 1}^{r-1}  \sum_{i=1}^N \int_0^1\bigg( \bigg( \sum_{\smallmatrix 1 \leq \d \leq r\\[2pt] 0 \leq \ve \leq \d\\[2pt]
\d  - \ve -1 =  \b  \endsmallmatrix } \hskip-0.5 cm (-1)^{\ve}  \frac{d^\ve}{dt^\ve} \left(\frac{\p  }{\p q^i_{(\d)}}  \left(L + \frac{ dC}{dt}\right) \right)   +    \pa^i_{\b(1)}  \bigg) Y^i_{(\b)} \bigg)\bigg|_{{\wh{F}}^{(n)}(0, s)} 
 ds  \overset{\eqref{primacond-A}} = 0\ .  \end{multline}
 \begin{rem} \label{remark52}  This identity is precisely what we planned to obtain  by introducing  the functions $\pa^i_\a(t)$ and   the corresponding terms in  $ \wt L$ and $\wh L$.  In fact, the integral \eqref{secondina} is precisely the integral of $\a^{PC}$ along the first  ``vertical side'' of $\cS$ (see Sect.\ \ref{smart-2} and \ref{roadmap}) and the functions $\pa^i_\a(t)$
{\it have been chosen so that such an  integral  is always  vanishing}.  
 \end{rem}
 Let us now consider the similar line integral  for  $t = T$. We have
 \begin{multline} \label{59}  \int_0^1 \imath_{\frac{\p}{\p s}}\left({\wh{F}}^{(n)*}(\a^{PC})\right)\bigg|_{(T,s)} ds = \int_0^1 \a^{PC}\left(Y_{\wh{F}}\right)|_{{\wh{F}}^{(n)}(T, s)}  ds = \\
=\sum_{\b = 1}^{r-1}  \sum_{i= 1}^N   \int_0^1  \bigg\{\bigg( \sum_{\smallmatrix 1 \leq \d \leq r\\[2pt] 0 \leq \ve \leq \d\\[2pt]
\d  - \ve -1 =  \b  \endsmallmatrix }  (-1)^{\ve}  \frac{d^\ve}{dt^\ve} \left(\frac{\p   \left(L + \frac{ dC}{dt}\right) }{\p q^i_{(\d)}} \right)  \bigg) Y^i_{(\b)}  + \hskip 4 cm \\
\hskip 0.5 cm  +  \bigg( \pa^i_{\b(1)} Y^i_{\b(0)} - \pa'{}^i_{\b(2)} Y'{}^i_{\b(1)} - \pa''{}^i_{\b(2)} Y''{}^i_{\b(1)} +  \pa'{}^i_{\b(3)} Y'{}^i_{\b(0)} + \pa''{}^i_{\b(3)} Y''{}^i_{\b(0)} \bigg) \bigg\}\bigg|_{{\wh{F}}^{(n)}(T, s)}  ds + \\
 +   \int_0^1  Y^{\{\mu\}}_{(0)} \bigg|_{{\wh{F}}^{(n)}(T, s)} ds \ .
\end{multline}
By property (2) of $Y_{\wh{F}}$ and   the boundary values  \eqref{ahio-A} and \eqref{ahiobis-A} of the   curves $\pa'{}^i_\b(t)$ and $\pa''{}^i_\b(t)$,  we have that 
 \begin{align*}
 & \pa'{}^i_{\b(2)}  \big|_{{\wh{F}}^{(n)}(T, s)}  = \hskip -0.3cm  \sum_{\smallmatrix 1 \leq \d \leq r\\[2pt] 0 \leq \ve \leq \d\\[2pt]
\d  - \ve -1 =  \b  \endsmallmatrix }  \hskip -0.3 cm (-1)^{\ve}  \frac{d^\ve}{dt^\ve} \left(\frac{\p   \left(L + \frac{ dC}{dt}\right) }{\p q^i_{(\d)}} \right) \!\!\bigg|_{j^{n-1}_{t = T}(q(t))} 
\!\!\!\!\!\!\ ,  & &Y'{}^i_{\b(1)}  \big|_{{\wh{F}}^{(n)}(T, s)} = Y^i_{(\b)}  \big|_{{\wh{F}}^{(n)}(T, s)}\ ,\\
  & \pa''{}^i_{\b(2)}  \big|_{{\wh{F}}^{(n)}(T, s)}  = \pa^i_{\b(1)} \big|_{{\wh{F}}^{(n)}(T, s)}  \ , \!\!\!\!\!\!& &\!\!\!Y''{}^i_{\b(1)}  \big|_{{\wh{F}}^{(n)}(T, s)}  = Y^i_{\b(0)} \big|_{{\wh{F}}^{(n)}(T, s)}   \ .
 \end{align*}
 From this and \eqref{59},  we get that 
 \begin{multline} \label{terzina}  \int_0^1 \imath_{\frac{\p}{\p s}}\left({\wh{F}}^{(n)*}(\a^{PC})\right)\bigg|_{(T,s)} ds = \\
 =   \int_0^1 \bigg(  \sum_{\b = 1}^{r-1}  \sum_{i= 1}^N \bigg(  \pa'{}^i_{\b(3)} Y'{}^i_{\b(0)} + \pa''{}^i_{\b(3)} Y''{}^i_{\b(0)}\bigg) + Y^{\{\mu\}}_{(0)} \bigg)\bigg|_{{\wh{F}}^{(n)}(T, s)} ds  \ .\end{multline}
 On the other hand, by the definition of the vector field $Y_{\wh F}$, 
  $$Y^{\{\mu\}}_{(0)} \bigg|_{{\wh{F}}^{(n)}(T, s)}  = \frac{\p \wh F^{\{\mu\}}}{\p s}\bigg|_{(T,s)} = \frac{\p (\wh \g^{(\wh U(s))})^{\{\mu\}}}{\p s}\bigg|_{(T,s)} = \frac{\p \wh \mu(T, s)}{\p s}\ .$$
Hence
\beq\label{quartina} 
\int_0^1  Y^{\{\mu\}}_{(0)} \bigg|_{{\wh{F}}^{(n)}(T, s)} ds =  \int_0^1 \frac{\p}{\p s}\Bigg( \int_0^T \frac{d \wh \mu(t, s)}{dt}  dt\Bigg) =  \iint_{[0,T] \times [0,1]}   \frac{\p^2\wh  \mu}{\p t \,\p s} ds dt \ .
\eeq
  \begin{rem}\label{remark53}  In  analogy with what we mentioned in Remark \ref{remark52},  the purpose of     the   functions $\pa'{}^i_\b(t)$,  $\pa''{}^i_\b(t)$ and  of  the corresponding   terms in  $\wh L$  was precisely  to get   \eqref{terzina}, which  simplifies   {\it  the integral of $\a^{PC}$ along the second ``vertical side'' of $\cS$.} 
 \end{rem}
From \eqref{primissima}, \eqref{secondina}, \eqref{terzina}, \eqref{quartina} and  the Stokes Theorem, we obtain that 
\beq \begin{split} \label{st}
 C_0 - C_1  + &  \int_0^T \left( \int_0^1  \frac{\p^2\wh  \mu}{\p t \,\p s}\bigg|_{(t,s)} \right. + \\
 & \hskip 2 cm\left.  +\sum_{\b = 1}^{r-1}  \sum_{i= 1}^N  \frac{\p}{\p t}  \bigg(  \pa'{}^i_{\b(3)} Y'{}^i_{\b(0)} + \pa''{}^i_{\b(3)} Y''{}^i_{\b(0)}  \bigg) \bigg|_{{\wh{F}}^{(n)}(t, s)} ds  \right) dt  = \\
&  =  \int_{\p ([0,T] \times [0,1])}\left({\wh{F}}^{(n)*}(\a^{PC})\right)   \overset{\text{Stokes Thm.}}= \\
 & =   \int_{[0,T] \times [0,1]} d\left({\wh{F}}^{(n)*}( \a^{PC})\right) \left (\frac{\p}{\p t}, \frac{\p}{\p s}\right)dt \,ds= \\
 & =  \int_{[0,T] \times [0,1]} \left({\wh{F}}^{(n)*}(d \a^{PC})\right) \left (\frac{\p}{\p t}, \frac{\p}{\p s}\right)dt \,ds= \\
& = \int_{[0, T] \times [0,1]} d \a^{PC}(X_{\wh{F}}, Y_{\wh{F}})\big|_{{\wh{F}}^{(n)}(t,s)} dt\, ds 
 \ ,
 \end{split}\eeq
 where, in the  integral of the first line, the integration of ${\wh{F}}^{(n)*}(\a^{PC})$  is performed along  the  usual counterclockwise parameterisation of  $\p ([0,T] \times [0,1])$. \par
At this point, it suffices to  recall  that the controlled Euler-Lagrange equations \eqref{modifiedEL} are satisfied at all points of $\wh \cS$ and that  the vectors $X_{\wh{F}}$ are  tangent vectors to   curves of jets, determined by solutions to the controlled Euler-Lagrange equations.  Indeed, due to this,   \eqref{39} and the fact that each of  the $1$-forms  \eqref{holonomic} vanish  identically  on the vectors $X_{\wh F}$, we immediately have  that  (\footnote{This is precisely  the point that motivated the substitution  of the $1$-form $\a$ by $\a^{PC}$.})
\beq 
\begin{split} d \a^{PC}(X_{\wh F}, Y_{\wh F}) \big|_{{\wh{F}}^{(n)}(t,s)}  =  \bigg(\frac{\p \wh L}{\p u^a} & d u^a \wedge dt \bigg) (X_{\wh F}, Y_{\wh F}) \big|_{{\wh{F}}^{(n)}(t,s)}  = \\
 & = - \frac{\p \wh L}{\p u^a} Y^a \big|_{{\wh{F}}^{(n)}(t,s)} =  - \frac{\p L}{\p u^a} Y^a \big|_{{\wh{F}}^{(n)}(t,s)} \ .
 \end{split}
 \eeq 
From this and \eqref{st}, the conclusion follows. 
\end{pf}
\par
\smallskip
\subsection{The  Principle of Minimal Labour}\hfill\par
\label{section7.2}
 Let  us now go back to our original triple $( \cK,  L,  C)$ and to the parameterised graphs $ \g(t) $ in $ \bR \times \cQ$.  We have already observed that for any $ \cK$-controlled curve
 $ \g^{( U)}(t) = (t, q^i(t))$
 there is a uniquely associated element $\wh U \in \wh \cK$ and a uniquely associated $\wh \cK$-controlled curve of the form 
 $$\wh \g^{(\wh U)} = (t, q^i(t), \pa^i_\b(t), \pa'{}^i_\b(t), \pa''{}^i_\b(t),  \l= 1, \mu^{(\wh U)}(t))\ ,$$
 i.e.  with the same components $q^i(t)$ of $\g^{(U)}$. On the basis of this,  in what follows  for  any given  $ U \in  \cK$ we  use the symbols $\wh U$ and  $\wh \g^{(\wh U)}$ to denote  the  uniquely associated element in $\wh \cK$ and the  corresponding $\wh \cK$-controlled  curve, respectively.
 \par
 \smallskip
Let  $ U_o  = (u_o(t),  \s_o(t))$ be a fixed element in  $ \cK$ and  $  \g^{( U_o)}: [0,T] \to [0,T] \times  \cQ $ the corresponding $ \cK$-controlled curve.  Exactly as in the previous section, we may now consider a {\it $ \cK$-controlled variation}  of this curve, i.e.  a smooth  $1$-parameter family $ F(\cdot, s)$, $s \in [0,1]$,  of  $ \cK$-controlled curves with  initial  curve $ F(\cdot, 0) =   \g^{( U_o)}$.   As we did in dealing with  $\wh \cK$-controlled variations, we denote by $ F^{(n)}$ the   
 homotopy in  $J^n( \cQ|\bR) \times K$
  $$ F^{(n)}: [0,T] \times [0,1] \longrightarrow J^n( \cQ|\bR)\times K\ ,\qquad  F^{(n)}(t,s) \= (j^{(n)}_t( F(\cdot, s)), u(t,s))\ , $$
determined by  the homotopy of the $n$-th order jets of the  curves  $\g^{( U(s))} =  F(\cdot, s)$ and the homotopy $u(t, s) = u^{(s)}(t)$ of control curves.  We  may also consider the vector fields 
  \begin{align}\label{vectorX-A} X_{ F}|_{ F^{(n)}(t,s)} &\=  F^{(n)}_*\left(\frac{\p}{\p t}\bigg|_{(t,s)}\right) = \frac{\p  F^{(n)}}{\p t}\bigg|_{(t,s)}\ ,  \\
\label{vectorY-A} Y_{ F}|_{ F^{(n)}(t,s)} & \=  F^{(n)}_*\left(\frac{\p}{\p s}\bigg|_{(t,s)}\right) = \frac{\p  F^{(n)}}{\p s}\bigg|_{(t,s)}\ ,
\end{align}
defined   at  the points  of the  surface  
 $ \cS \=  F^{(n)}([0,T] \times [0,1])  \subset J^n( \cQ|\bR) \times K$.  By construction,  for each fixed  $s_o \in [0,1]$ we have that:  (a) the restriction of  $X_{ F}$ to the trace of   the  curve   
 $$t \mapsto  \left( \g^{( U(s_o))(n)}(t) , u(t, s_o)\right) = \left(j^n_{t}\left({{F}}(\cdot s_o)\right), u(t, s_o)\right)$$ 
 is the family of    the tangent vectors of such a curve; (b) the restriction of  $Y_{ F}$ to the trace of the same  curve   is  the Jacobi vector field   corresponding to the  (infinitesimal) variation $\ve \to  F^{(n)}(\cdot, s_o+\ve)$  of $ F^{(n)}(\cdot, s_o) $; (c) the vector fields  $X_{ F}$ and $Y_{ F}$ have the form
 \beq \label{54uno-A}   X_{ F} = \frac{\p}{\p t} +   X^i_{(\b)} \frac{\p}{\p q^i_{(\b)}} + X^a \frac{\p}{\p u^a}    \ ,\qquad
  Y_{ F} =    Y^i_{(\b)}  \frac{\p}{\p q^i_{(\b)}}   + Y^a \frac{\p}{\p u^a}\ . 
  \eeq
  Finally,   we have that  the curve  $ U(s) \in  \cK$  determines uniquely  a    curve $\wh U(s)$ in $ \wh \cK$ and a corresponding  $\wh \cK$-controlled variation ${\wh{F}}$  starting from the curve   $\wh \g \= \wh \g^{(\wh U_0)}$ and ending with   the curve $\wh \g' \= \wh \g^{(\wh U_1)}$.\par
\smallskip
  We are now ready to establish the following corollary of  Theorem \ref{prop34} and Remark \ref{rem72}. 
 \begin{cor}[Homotopy Formula]  \label{cor34}  Let $ U_0,  U_1 \in  \cK$  be the endpoints of a smooth curve $ U(s) \in  \cK$, $s \in [0,1]$, and 
 $ \g \=  \g^{( U_0)}$,  $ \g' \=  \g^{( U_1)}$ the $ \cK$-controlled curves corresponding to  $ U_0,  U_1$, with terminal costs  $C_0 \= C(j^{n-1}_{t = T}( \g))$  and  $C_1 \= C(j^{n-1}_{t = T}( \g'))$, respectively.   Let also  $\cP_{j^n_t(\g)} : K \to \bR$ and $\mu, \mu': [0,T] \times [0,1] \to \bR$ be the real functions defined by 
 \beq  \cP_{j^n_t(\g)}(u^a) \= -   L(j^n_t(\g), u^a)\ ,\eeq
\begin{multline} \label{7.17} \mu(t, s)  \=  -  \int_0^t \bigg(L  + \frac{1}{2}\!\!\!\sum_{\smallmatrix 1 \leq i \leq N\\
   0 \leq  \b \leq r-1
   \endsmallmatrix} \!\!\! \bigg( ( \pa^i_{\b(1)})^2 {-} ( \pa'{}^i_{\b(2)})^2 {-}  ( \pa''{}^i_{\b(2)})^2 
    \bigg) + \\
  + \!\!\! \sum_{\smallmatrix 1 \leq i \leq N\\
   0 \leq  \b \leq r-1
   \endsmallmatrix} \!\!\! \bigg( \frac{1}{2} ( \pa^i_{\b})^2+ \frac{\pi^4}{32 T^4} ( \pa'{}^i_{\b})^2 +  \frac{\pi^4}{32 T^4}( \pa''{}^i_{\b})^2  \bigg)  \bigg)\bigg|_{ \g^{( U(s))}(\t)}  d\t
   \end{multline}
   \beq \label{7.17'}  \mu'(t, s)  \= \mu(t, s) +  \sum_{\b = 1}^{r-1}  \sum_{i= 1}^N  \int_0^s \left(    \pa'{}^i_{\b(3)} Y'{}^i_{\b(0)} + \pa''{}^i_{\b(3)} Y''{}^i_{\b(0)}  \right)\bigg|_{{\wh{F}}^{(n)}(t, v)} d v   \eeq
where we denoted by $\wh F$ the uniquely determined  $\wh \cK$-controlled variation described above and by $ Y'{}^i_{\b(0)} $, $ Y''{}^i_{\b(0)} $ 
   the components in the directions of the coordinates $\pa'{}^i_{\b(0)}$,  $\pa''{}^i_{\b(0)}$ of the  vector field $Y_{\wh F}$ defined in \eqref{vectorY}. Let us finally denote by  $Y_{ F}$ the vector field  in \eqref{54uno-A},  associated with the $\cK$-controlled variation $F$ determined by the  $U(s)$, $s \in [0,1]$, 
 Then: 
 \begin{itemize}
 \item[(i)] The difference between the  terminal costs $C_0$ and $C_1$ is equal to 
\beq \label{homotopy-ultimate-A} 
\begin{split}
 C_1 - C_0   =   - \int_0^T \left( \int_0^1    \left(Y^a \frac{\p \cP_{\s^{(s)}_t} }{\p u^a} \Bigg|_{u^{(s)}(t)} -  \frac{\p^2  \mu'}{\p t\, \p s}\bigg|_{(t,s)} \right)ds \right) dt  \ , \quad \s^{(s)}_t \= j^n_t( \g^{( U(s))})\ ;
 \end{split}
\eeq
\item[(ii)] For any $s_o \in [0,1]$,   
\beq\label{controlledNoether}  \a^{PC}\left(Y_{\wh F}\right)|_{\wh F^{(n)}(T, s_o)}  = \a^{PC}\left(Y_{\wh F}\right)|_{\wh F^{(n)}(0, s_o)} + \int_0^T  \Bigg( \frac{\p^2 \mu'}{\p t\,\p s}\bigg|_{s = s_o} \Bigg)dt \ .\eeq
\end{itemize}
\end{cor}
\begin{pf}   By Theorem \ref{prop34},  Remark \ref{rem72} and the fact that, due to the Euler-Lagrange equations,  the \eqref{7.17} is nothing but  the explicit expression for the $\mu$-components of the curve $\wh \g^{(\wh U(s))}(t)$ (thus, of the function $\wh \mu(t,s)$ of Theorem \ref{prop34}),  claim (i) follows immediately. \par
We  now prove (ii).  Consider the embedded surface with boundary
 $\wh  \cS \=  \wh F^{(n)}([0,T] \times [0,1])$  and the vector fields  $X_{\wh F}$, $Y_{\wh F}$, which we  defined above at the points of  $ \wh \cS$. 
From  \eqref{extendedaction} and  the fact that $\a^{PC}$ is  variationally equivalent  to $\a$,  
we have that 
\beq  \label{primissima-A} 
\begin{split} C_0 - C_1 &= C(j^n_{t = T}( \g^{( U_0)})) - C(j^n_{t = T}( \g^{( U_1)})) = \\
& = \int_0^T \imath_{\frac{\p}{\p t}}\left({{\wh F}}^{(n)*}(\a^{PC})\right)\bigg|_{(t,0)}  dt  -
\int_0^T \imath_{\frac{\p}{\p t}}\left({\wh F}^{(n)*}(\a^{PC})\right)\bigg|_{(t,1)}  dt \ .
\end{split}
\eeq
On the other hand, 
\beq 
\begin{split} \label{secondina-A} &\int_0^1 \imath_{\frac{\p}{\p s}}\left( \wh F^{(n)*}(\a^{PC})\right)\bigg|_{(0,s)} ds  = \int_0^1 \a^{PC}\left(Y_{ \wh F}\right)|_{ \wh F^{(n)}(0, s)}  ds\ ,\\
  &\int_0^1 \imath_{\frac{\p}{\p s}}\left( \wh F^{(n)*}(\a^{PC})\right)\bigg|_{(T,s)} ds = \int_0^1 \a^{PC}\left(Y_{ \wh F}\right)|_{ \wh F^{(n)}(T, s)}  ds \ .
 \end{split}
 \eeq
From \eqref{primissima-A},  \eqref{secondina-A}  and  the Stokes Theorem, we get that 
\beq \begin{split} \label{st-A}
 C_0 - C_1 
 & =  \int_{\p ([0,T] \times [0,1])}\left( \wh F^{(n)*}(\a^{PC})\right)\bigg|_{(t,s)}  + \\
 & \hskip 2 cm + \int_0^1 \imath_{\frac{\p}{\p s}}\left( \wh F^{(n)*}(\a^{PC})\right)\bigg|_{(0,s)} ds - 
 \int_0^1 \imath_{\frac{\p}{\p s}}\left( \wh F^{(n)*}(\a^{PC})\right)\bigg|_{(T,s)} ds = \\
 &=   \int_{[0,T] \times [0,1]} \left( \wh F^{(n)*}(d \a^{PC})\right)\bigg|_{(t,s)}  + \\
 & \hskip 2 cm +  \int_0^1 \a^{PC}\left(Y_{ \wh F}\right)|_{ \wh F^{(n)}(0, s)} -  \int_0^1 \a^{PC}\left(Y_{ \wh F}\right)|_{ \wh F^{(n)}(T, s)}  ds = \\
& = \int_{[0, T] \times [0,1]} d \a^{PC}(X_{ \wh F}, Y_{ \wh F})\big|_{ \wh F^{(n)}(t,s)} dt\, ds + \\
 & \hskip 2 cm +  \int_0^1 \a^{PC}\left(Y_{ \wh F}\right)|_{ \wh F^{(n)}(0, s)} -  \int_0^1 \a^{PC}\left(Y_{ \wh F}\right)|_{ \wh F^{(n)}(T, s)}  ds
 \ .
 \end{split}\ .\eeq
As in the proof of Theorem \ref{prop34},  we now observe  that   the vectors $X_{ \wh F}$ are  the tangent vectors to   curves of jets  that are  solutions to the Euler-Lagrange equations \eqref{modifiedEL}.  Hence, by Lemma \ref{pcprop}, at any   point of $ \cS =  \wh F^{(n)}([0,T] \times [0,1]) $,  we have  that 
    \beq 
   d   \a^{PC}\left(X_{ \wh F},  Y_{ \wh F}\right)   = - Y^a \frac{\p  L}{\p u^a} =   Y^a {\frac{\p \cP_{\s^{(s)}_t}}{\p u^a}}\bigg|_{u^{(s)}(t)}\ . 
\eeq
From this and \eqref{st-A}, it follows that 
\beq
\begin{split} C_0 - C_1 &=  \int_0^T \left( \int_0^1 \left(Y^a \frac{\p \cP_{\s^{(s)}_t} }{\p u^a}  \right)\Bigg|_{u^{(s)}(t)}  ds \right) dt  + \\
 & \hskip 2 cm +  \int_0^1 \a^{PC}\left(Y_{ \wh F}\right)|_{{ \wh F}^{(n)}(0, s)} -  \int_0^1 \a^{PC}\left(Y_{ \wh F}\right)|_{{ \wh F}^{(n)}(T, s)}  ds
 \ .
 \end{split}
\eeq
 This and (i)  imply that 
$$  \int_0^1 \a^{PC}\left(Y_{ \wh F}\right)|_{{ \wh F}^{(n)}(T, s)} ds =  \int_0^1 \a^{PC}\left(Y_{ \wh F}\right)|_{{ \wh F}^{(n)}(0, s)}  ds +  \int_0^1 \left( \int_0^T  \frac{\p^2   \mu'}{\p t\,\p s} dt \right) ds   \ . $$
This identity holds not only for the $ \cK$-controlled homotopy $ F$, but, for a fixed  $s_o \in [0,1]$ and a fixed sufficiently  small $\ve > 0$, it holds 
also for any other $ \cK$-controlled homotopy  of the form  
$$ F^{(s_o, \ve)}(t,\t) \=  F\left(t,  s_o (1 - \t) +  (s_o + \ve) \t \right)\  ,\qquad \t \in [0, 1]\ ,
$$
 which  interpolates between  the  curve $\g^{(U(s_o))}$ to the  curve  $\g^{(U(s_o + \ve))}$. 
This implies that for any fixed choice of  $s_o \in [0,1)$ and any small $\ve > 0$
$$  \int_{s_o}^{s_o + \ve} \a^{PC}\left(Y_{ \wh F}\right)|_{{ \wh F}^{(n)}(T, s)} ds  {=}  \int_{s_o}^{s_o + \ve} \a^{PC}\left(Y_{ \wh F}\right)|_{{ \wh F}^{(n)}(0, s)}  ds + \int_{s_o}^{s_o+ \ve}  \left( \int_0^T   \frac{\p^2  \mu'}{\p t\,\p s} dt \right)ds.$$
From this and the continuity of  all involved functions, claim (ii) follows.
\end{pf}
\par
The homotopy formula  has the following   immediate consequence, which is the first main result that was expected according to  the road map. \par
\begin{theo}[Generalised Principle of Minimal Labour] \label{corollone}   A necessary condition for a $ \cK$-controlled curve $ \g_o \=  \g^{( U_o)}$ to be a  solution to the considered generalised Mayer  problem
is that   for any $ \cK$-controlled variation  $ F$ with $ F(t,0) = \g_o$ one has that 
  \beq \label{formulacor}  \int_0^T \left( \int_0^1    \left(Y^a \frac{\p \cP_{\s^{(s)}_t} }{\p u^a} \Bigg|_{u^{(s)}(t)} - \frac{\p^2 \mu'}{\p t\,\p s}\bigg|_{(t,s)}\right)ds \right) dt \leq 0 \ .\eeq
In case  the considered problem is such that for  any two $ \cK$-controlled curves $ \g \=  \g^{( U)}$ and $ \g' \=  \g^{( U')}$  there is a 
$ \cK$-controlled variation which has  them as endpoints,  the above  is also  a sufficient condition. 
\end{theo}
\par
\medskip
Notice that the above principle  admits  the following equivalent formulation.\\[10pt]
 {\it Let   $ \g_o \=  \g^{( U_o)}$ be  a fixed $ \cK$-controlled curve and  for any  $ \cK$-controlled variation $ F$ with $ F(t,0) = \g_o$
consider the real function  $\cW^{ F}: [0,1] \to \bR$ defined by
\beq  \label{formulacor-bis}
 \cW^{ F}(\d) \=  \int_0^T \left( \int_0^\d    \left(Y^a \frac{\p \cP_{\s^{(s)}_t} }{\p u^a} \Bigg|_{u^{(s)}(t)} -  \frac{\p^2  \mu' }{\p t\,\p s}\bigg|_{s}\right)ds \right) dt 
\ .\eeq
Then  $ \g_o \=  \g^{( U_o)}$ is a solution to the   Mayer  problem  determined by $( \cK,  L,  C)$ only if for any $ \cK$-variation as above 
\beq   \label{formulacor-bis1}  \cW^F(\d) \leq 0\qquad \text{for each}\ \d \in [0,1]\ .\eeq
}
\\
It should be pointed out that  an explicit check of \eqref{formulacor-bis1}  is expected to  be  quite hard in  generic situations:  it demands the  study  of the sign 
behaviour of the functions $\cW^F$  for {\it any}   $ \cK$-controlled variation $F$  of a candidate $ \g_o$.  Nonetheless it has    
a number of  consequences,  the most elementary one  represented  by an infinitesimal version of Theorem  \ref{corollone},   the concluding result of this section. \par
Let $ \g_o \=  \g^{( U_o)}$ be a $ \cK$-controlled curve associated with the pair $ U_o = (u_o(t),  \s_o)$  and denote by 
$\cF$   the full collection   of   the $ \cK$-controlled variations of $ \g_o$.  Let also denote by $\cJ ac^{(n)}$ the class of  all  vectors fields $V$ in $T J^n(\cQ|\bR)$, defined just at the points of the trace of $ \g^{(n)}_o(t)$, 
 of the form 
$V|_{ \g^{(n)}_o(t)} \= Y_{ F}|_{ \g^{(n)}_o(t)}$ for some $F \in \cF$.   In other words,  $\cJ ac^{(n)}$ is the family of all  Jacobi vector fields  of   $ \g^{(n)}_o(t) $, which are  determined by the variations of curves of jets determined by  the $ \cK$-controlled variations of $ \g_o$. 
\par
\smallskip
\begin{theo}  
 A  $ \cK$-controlled curve   $ \g_o \=  \g^{( U_o)}$ is    a solution to the  Mayer  problem  of  $( \cK,  L,  C )$ only if the following two conditions hold for any $V \in  \cJ ac^{(n)}$: 
 \begin{itemize}[itemsep = 10pt]
 \item[(1)] $dC \left(V|_{ \g^{(n)}_o(T)}\right) \geq  0$; 
 \item[(2)] $ \displaystyle  \int_0^T   \left(Y^a \frac{\p \cP_{\s^{(s)}_t} }{\p u^a} \Bigg|_{u_o(t)} -  \frac{\p^2   \mu' }{\p t\,\p s}\bigg|_{s = 0}\right) dt  \leq 0$. 
 \end{itemize}
\end{theo}
\begin{pf}  (2) is an immediate consequence of     \eqref{formulacor-bis1}.  For (1), one can obtain it by just  observing the function $\cW^F$ coincides with  the map 
$\d\overset{\cW^F} \longmapsto C|_{F^{(n)}(T, 0)} -  C|_{F^{(n)}(T, \d)}$. The inequality is  then obtained  taking the derivative of this expression at $\d = 0$. 
\end{pf}
\par
\medskip
  \section{The generalised   Pontryagin Maximum Principle}
  \label{Sect7}
  In this section we want to show  that  the Principle of Minimal Labour yields  a  strict  analogue   of the classical PMP   for a wide class of  generalised Mayer problems of  higher order differential constraints of variational type  and without the usual restriction of  fixed initial values.  The class we  consider  is characterised by  {\it   differential constraints of normal type}, which are   defined in the next subsection. Immediately after this, we introduce a  generalised version  of the classical Pontryagin needle variations and we finally  prove the advertised result. 
  \par
  \smallskip  
\subsection{Differential constraints of normal type}  \hfill\par
\label{defnormal}
As usual,  we consider a Mayer problem determined by a triple $(\cK, L, C)$ with controlled Lagrangian $L$ of actual order $r\geq 1$ and associated   controlled Euler-Lagrange equations  of order at most $2 r$.  Let us    indicate such a system of  differential constraints  as  
\beq \label{first-1} \cE^j\left(t, q^i, \frac{d q^i}{dt}, \ldots, \frac{d^{2r} q^i}{dt^{2r}}, u^a(t)\right) = 0\ ,\qquad 1 \leq i \leq N\ .\eeq
By considering  appropriate auxiliary variables,   say 
$$y^1(t) = q^1(t) \ ,\ \ y^2(t) = q^2(t)\ ,\ \  \ldots\ \ , y^{N+1}(t) = \frac{d q^1}{dt}\ ,\ \  y^{N+2}(t) = \frac{d q^2}{dt}\ ,\ \ \ldots\ ,$$
the original  system  \eqref{first-1} can be always transformed  into an equivalent system on the new functions $y^A(t)$, which consists only of {\it first} order differential  equations    of the form 
\beq \label{first-2}  \cG^B\left(t, y^1, y^2,  \ldots, y^{\wt N},  \frac{d y^1}{dt}, \frac{d y^2}{dt}, \ldots,  \frac{d y^{\wt N}}{dt},  u^a(t)\right) = 0\ ,\qquad  1 \leq B \leq \wt N' \eeq
for  some appropriate    $\wt N$-tuples of variables  and $\wt N'$-tuples of equations, with  both $\wt N$ and $\wt N'$ greater than or equal to $N$.  In general, there is not a unique way to reduce the constraints \eqref{first-1} into a form \eqref{first-2}. For instance,  the second order equation on curves $q(t) \in \cQ = (0, + \infty)$ 
\beq\label{ex} q \frac{d^2 q}{dt^2} + \left(\frac{d q}{dt} \right)^2 - u(t) = 0\eeq can be reduced not only  to the equivalent first order system 
\beq\label{bad} y^1 \frac{d y^2}{dt}  + (y^2)^2 - u = 0\ ,\qquad \frac{d y^1}{dt} - y^2 = 0\ ,\eeq 
but also to the  system (via the change of variable $\wt y^1 = q^2$)
\beq \label{good} \frac{d \wt y^2}{dt} - 2 u = 0\ ,\qquad  \frac{d \wt y^1}{dt} - \wt y^2 = 0\ .\eeq 
We say that {\it the differential constraints \eqref{first-1} are  of normal type} if they  are equivalent to at least one   first order differential system  \eqref{first-2}, 
consisting only of  first order equations in normal form, that is having the form
\beq \label{first-3}  \frac{d y^A}{d t} = g^A(t, y^B,  u^a(t)) \ , \eeq
{\it with some smooth function $g = (g^A): \bR \times \bR^{\wt N} \times K \subset \bR^{ \wt N + M +1} \to  \bR^{\wt N}$}. 
Note that, in the previous  example, the system \eqref{bad} does not satisfy such a regularity  assumption, but \eqref{good} does. This is good enough for us to consider the  \eqref{ex} as a constraint of normal type. \par
\smallskip
 In case it is possible to transform the original constraints into some of the form \eqref{first-3} in which the $g^A$ 
 and the partial derivatives $\frac{\p g^A}{\p y^B}\big|_{(t, y^B, u)}$ have  
uniform    bounds   we say that they are 
of {\it bounded normal type}. 
\par
\begin{rem} \label{cut-off}
In dealing with   localised properties of  solutions of  a set of differential constraints,  
using cut-off functions  it is in general possible to  replace some (possibly unbounded) smooth functions  $g^A$ by  other 
 functions  $g'{}^A$, 
coinciding with the $g^A$ on an appropriate compact set and satisfying global bounds for their values and  the partial derivatives $\frac{\p g^A}{\p y^B}$. 
This idea  is a common tool   for  extending  estimates  on  solutions of  constraints of {\it bounded} normal type to   those of  more general  constraints of normal type.
\end{rem}
Consider now  the distance function on the family of control curves  $u, u': [0,T] \to K \subset \bR^M$ defined by 
$$\operatorname{dist}(u, u') \= \text{measure}\left\{\ t \in [0, T]\ :\ u(t) \neq u'(t)\  \right\}$$
and denote by $\r$ a fixed metric on the jet space  $J^{n-1}(\cQ|\bR)|_{t =0}$ of initial values, which generates the standard topology of such a space.
The next lemma is an immediate consequence of a classical fact on  systems of controlled first order differential equations (see e.g. \cite[Prop. 3.2.2]{BP}). 
It gives  the main motivation for considering the class of  differential constraints  of normal type.  
\begin{lem} \label{normallemma} Let $(\cK, L, C)$ be a defining triple of a Mayer problem with differential constraints of order $2r$  of  bounded normal type and $K \subset J^{n-1}(\cQ|\bR)|_{t = 0}$ a compact subset in the  space of  the $(n-1)$-jets at $t = 0$.   There exists 
a constant ${\mathfrak c}$,  depending only on $L$, such that for any two  $\cK$-controlled curves
$$\g^{(U)}, \g^{(U')} : [0, T] \longrightarrow [0,T] \times \cQ \ ,$$ determined by pairs  $U  = (u(\cdot), \s)$, $U'= (u'(\cdot), \s') $  with  $\s$, $\s'$ in $K$, one has 
\beq \| \g^{(U)} - \g^{(U')}\|_{\cC^{2r-1}} \leq {\mathfrak c}\, \left(\operatorname{dist}(u, u') + \r(\s, \s')\right)\ .\eeq
\end{lem}  
\par
\smallskip
 \subsection{Generalised needle variations}\hfill\par
 \label{generalisedneedle}
Let  $ U_o  = (u_o(\cdot),  \s_o)$ be a fixed element in  $ \cK$ and   $  \g^{( U_o)}: [0,T] \to [0,T] \times  \cQ $ the uniquely corresponding $ \cK$-controlled curve. 
Pick also a strictly positive  time $\t \in (0, T)$, a point $\o \in K$ and two real numbers  $0 < k < < 1$ and $0 < \ve_o $  small enough so that  $[\t - \ve_o - k \ve_o^2, \t + k \ve_o^2] \subset (0,T)$.  After fixing $k$ and the triple $(\t, \o, \ve_o)$,
we may consider  the piecewise continuous map  
\beq \label{uk}
u^{(\t, \o, \ve_o)}: [0, T] \to K\ ,\qquad u^{(\t, \o, \ve_o)}(t) \= \left\{\begin{array}{ll} u_o(t) & \text{if} \ t \in \big[0, \tau- {\ve_o} \big),\\[4pt]
 \omega & \text{if}\ t \in \big[\tau - {\ve_o},  \tau\big),\\[4pt]
 u_o(t) & \text{if} \ t \in \big[\tau, T\big]\ \\
 \end{array}\right.
 \eeq 
 and  an appropriate associated {\it smooth} map  $\check u^{(\t, \o, \ve_o)}:  [0, T] \to \check K$ with values in the convex hull $\check  K \subset \bR^M$ of a slightly larger  open neighbourhood  of $K$, described as follows. 
 We assume that $\check u^{(\t, \o, \ve_o)}(t)$  is equal to   $u^{(\t, \o, \ve_o)}(t)$ for  all  points $t$   in $[0,T]$ with the only exception  of  those in two  small intervals of the form  $[\t - \ve_o - k \ve_o^2,  \t - \ve_o]$ and $[ \t, \t + k \ve_o^2]$ with $k<<1$, in which  the function $\check u^{(\t, \o, \ve_o)}(t)$  is required just to take values in $\check K$, with no   further restrictions.  We call $u^{(\t, \o, \ve_o)}$  the {\it needle modification of $u_o$ at $t = \t$ of ceiling value $\o$ and  width $\ve_o$}. Any associated smooth approximation  $\check u^{(\t, \o, \ve_o)}$ will be  called {\it smoothed needle modification} (see Fig. 2 and Fig. 3).  
There are of course several ways to build up a smoothed needle modification for  a given discontinuous one. Nonetheless,   we assume that  a fixed algorithm  has been chosen and that each discontinuous needle modification has a uniquely associated  smoothed one.   \\
 \medskip
 \centerline{
\begin{tikzpicture}
%\draw[help lines] (0,0) grid (13,4);
\draw[<->, line width = 1] (1,3.8) to (1,0.5) to (6.5,0.5);
\draw[<->, line width = 1] (7,3.8) to (7,0.5) to (12.5,0.5);
\draw[fill]  (3, 0.5) circle [radius = 0.06];
\node at  (3, 0.2) {\tiny $\t - \ve_o$};
\draw[fill]  (4.5, 0.5) circle [radius = 0.06];
\node at  (4.5, 0.2) {\tiny $\t$};
\draw[fill]  (6.2, 0.5) circle [radius = 0.06];
\node at  (6.2, 0.2) {\small $T$};
\draw[fill]  (9, 0.5) circle [radius = 0.05];
\draw[fill,purple]  (8.5, 0.5) circle [radius = 0.06];
\node[purple] at  (8.1, 0.2) {\tiny $\t {-} \ve_o{-} k\ve_o^2$};
%\draw[fill,purple]  (9.3, 0.5) circle [radius = 0.06];
\node[purple] at  (11.2, 0.2) {\tiny $\t {+} k\ve_o^2 $};
\draw[fill]  (10.5, 0.5) circle [radius = 0.05];
%\draw[fill, purple]  (10.2, 0.5) circle [radius = 0.06];
\draw[fill,purple] (11, 0.5) circle [radius = 0.06];
\draw[fill]  (12.2, 0.5) circle [radius = 0.06];
\node at  (12.2, 0.2) {\small $T$};
\draw[fill]  (1, 3) circle [radius = 0.06];
\node at  (0.8, 3) {\small $\o$};
\draw[fill]  (7, 3) circle [radius = 0.06];
\node at  (6.8, 3) {\small $\o$};
\node[blue] at (2.3,1.7) {\tiny$u_o(t)$};
\node[blue] at (5.3,1.9) {\tiny$u_o(t)$};
\node[black] at (3.7, 2.5) {$u^{(\t, \o, \ve_o)}$};
\draw [line width = 0.7, blue] (1, 1) to [out=-20, in=200] (3,1.7)   ; 
\draw [line width = 0.7, blue, dashed](3,1.7)  to (3, 3)   ; 
\draw [line width = 0.7, blue] (3, 3) to  (4.5,3)   ; 
\draw [line width = 0.7, blue, dashed] (4.5,3)  to(4.5, 1.7)   ; 
\draw [line width = 0.7, blue] (4.5, 1.7) to [out=-20, in=150] (6.2,1.5)   ; 
\draw [line width = 0.7, blue] (7, 1) to [out=-20, in=200] (9,1.7)   ; 
\draw [line width = 0.7, blue, dashed](9,1.7)  to (9, 3)   ; 
\draw [line width = 0.7, blue] (9, 3) to  (10.5,3)   ; 
\draw [line width = 0.7, blue, dashed] (10.5,3)  to(10.5, 1.7)   ; 
\draw [line width = 0.7, blue] (10.5, 1.7) to [out=-20, in=150] (12.2,1.5)   ; 
\draw [line width = 1, purple](8.5,1.46)  to  [out=20, in=180] (9.1, 3) ; 
\draw [line width = 1, purple] (10.4,3)  to  [out=0, in=180] (11, 1.64)   ; 
\node[blue] at (8.3,1.7) {\tiny$u_o(t)$};
\node[blue] at (11.3,1.9) {\tiny$u_o(t)$};
\node[black] at (9.7, 2.5) {$\check u^{(\t, \o, \ve_o)}$};
 \end{tikzpicture}
 }
   \centerline{\tiny\hskip 1 cm \bf Fig. 2 Needle modification \hskip 2 cm \bf Fig. 3 Smoothed needle modification}
\ \\
 \smallskip
From now on, we also assume the following convenient assumption:  \\[2pt]
 {\it The differential problems that determine the $ \cK$-controlled 
curves are well defined and with unique solutions  also for the pairs $( u(t), \s)$, in which the initial datum  $\s$  is as usual and   $ u(t)$  takes values in  the convex hull    $\check K$ of some open neighbourhood  of  $K$}.\\[4pt]
 On the basis of this assumption, we 
may speak of a  $ \cK$-controlled curve $ \g^{( U)}$,  $ U = ( u(t),  \s)$, not only when $u([0,T])$ is entirely included  in $K$,  but also when it is  within an appropriate convex set $\check K \supset K$.\par
\smallskip
In the next definition, the function $\Sigma = \Sigma(\ve, s)$ is a continuous two-parameters family of initial data in $J^{n-1}(\cQ|\bR)|_{t = 0}$ satisfying the condition   $\Sigma(\ve, 0) = \s_o$ for each $\ve \in [0, \ve_o]$. 
\par
\begin{definition}
The {\it (generalised) needle variation of  the $ \cK$-controlled curve $\g^{( U_o)}$, corresponding to the triple $(\t, \o, \ve_o)$ and led by the family  $\Sigma = \Sigma(\ve, s)$} is the  one-parameter family of $ \cK$-controlled variations 
\beq\label{Cardino} \needle^{(\t, \o, \ve_o, \Sigma)} (\g^{(U_o)})\= \{\  F^{(\t, \o, \ve, \Sigma)}:[0,T] \times [0,1] \to [0,T] \times  \cQ\ , \ \ve \in (0, \ve_o]\ \}\ ,\eeq  in which the  variations $ F^{(\t, \o, \ve, \Sigma)}$ 
are  determined as follows.   For each $\ve \in (0,\ve_o]$,   let $ U^{(\ve, s)} = (u^{(\ve, s)}(t),  \s^{(\ve, s)} =   \s_o)$ be the curve in $ \cK$, in which:
\begin{itemize}[leftmargin = 20pt]
\item the initial    data $ \s^{(\ve, s)}$ are given by   $ \s^{(\ve, s)} =   \Sigma(\ve, s)$; 
\item  the one-parameter family of curves $u^{(\ve, s)}(t)$ takes values in the convex set  $\check K$ and is defined    by 
\beq  u^{(\ve, s)}(t) = (1-s) u_o(t) + s \check u^{(\t, \o, \ve)}(t)\ ,\qquad s \in [0,1]\ .\eeq
\end{itemize}
Then, we define  $ F^{(\t, \o, \ve, \Sigma)}$  as the  $ \cK$-controlled variation 
$$ F^{(\t, \o, \ve, \Sigma)}(t, s) \=  \g^{( U^{(\ve, s)})}(t)\ .$$
\end{definition}
In the following, given a needle variation $\needle^{(\t, \o, \ve_o, \Sigma)}$,  for each $\ve \in (0, \ve_o]$ we  denote by $ \mu'{}^{(\t, \o, \ve, \Sigma)}, \mu'{}^{(\t, \o, \ve, \Sigma)}:[0,T] \times [0,1] \to \bR$ the   functions defined in  \eqref{7.17},  \eqref{7.17'},  by means of  the homotopy $F^{(\t, \o, \ve, \Sigma)}$. \par
 \smallskip
\subsection{The generalised   Pontryagin Maximum Principle}\hfill\par
We are now ready to prove our announced analogue of the PMP for the class of generalised Mayer problems with differential constraints of normal type. As in all previous sections,  we  still consider  a fixed generalised Mayer problem  determined by one of the triples $(\cK, L, C)$ described in Sect.\ \ref{definingtriples}.
We start with a preparatory lemma. \par
\begin{lemma}   \label{propnormal}
Let   $ \g_o \=  \g^{( U_o)}$ be a  $ \cK$-controlled curve  and  $\needle^{(\t, \o, \ve_o, \Sigma)} (\g_o)= \{ F^{(\t, \o, \ve, \Sigma)}, 0 < \ve \leq \ve_o\}$ a needle variation of $\g_o$ 
 with associated  function  $\mu'{}^{(\t, \o, \ve, \Sigma)}: [0,T] \times [0,1] \to \bR$ as defined in  \S\ \ref{generalisedneedle}.
If  the Mayer problem has differential constraints of normal type,  the  limit
$\lim_{\ve  \to 0}
 \cP_{j^{r}_{\t} ( \g^{( U^{(\ve,s)})})}(\o) $
  exist and is  equal to 
\beq  \label{firstlimit-bis}  \lim_{\ve  \to 0} \cP_{j^{r}_{\t} ( \g^{( U^{(\ve, s)})})}(\o) = \cP_{j^{r}_{\t} ( \g^{( U_o)})}(\o) 
 \ .\eeq
\end{lemma}
\begin{pf}  First of all, we recall  that, for any $\ve \in (0, \ve_o]$, the curve of controls $U^{(\ve, s)} = (u^{(\ve, s)}(t),  \s^{(\ve, s)} =   \Sigma(\ve, s))$, $s \in [0, 1]$, have all  initial data $\s^{(\ve, s)}$ in the compact set $\Sigma([0, \ve_o] \times [0,1])$  and  all control curves  $u^{(\ve, s)}$  differ  from the control curve $u_o(t)$ only at the points of an interval of measure $\ve$ (more precisely, of measure  $\ve + o(\ve)$  since we are considering  smoothed needle variations). 
Thus, 
 since all of  the  $ \g^{( U^{(\ve,s)})}$ and their jets lie in a compact subset of $J^{2r+1}(\cQ|\bR)$,  by   the standard use of cut-off functions  described  Remark \ref{cut-off} and Lemma \ref{normallemma}, we have that  $\| \g^{(U(\ve, s))} - \g^{(U_o)}\|_{\cC^{2r-1}} $ goes   to $0$ for $\ve \to 0$. From this  the function  $\cP_{j^{r}_{\t} ( \g^{( U^{(\ve,s)})})}: K \to \bR$ tends  
 to the function $\cP_{j^{r}_{\t} ( \g^{( U_o)})}: K \to \bR$. 
\end{pf}
  \begin{theo}[Generalised Pontryagin Maximum Principle] \label{theorem51}\label{POMPone} \label{generalizedPOMP}
Let   $ \g_o \=  \g^{( U_o)}$ be a  $ \cK$-controlled curve for a generalised Mayer problem  with  differential constraints of normal type determined by a triple $(\cK, L, C)$. A necessary condition for  
$ \g_o$ to be  a solution   to the  Mayer  problem   is that for any    $\needle^{(\t, \o, \ve_o, \Sigma)} (\g_o)= \{ F^{(\t, \o, \ve, \Sigma)}, 0 < \ve \leq \ve_o\}$ 
  \beq  \label{712}  \cP_{j^{r}_{\t} ( \g_o)}(\o) -  \liminf_{\ve \to 0^+}  \frac{\mu'{}^{(\t, \o, \ve, \Sigma)}(T, 1) - \mu'{}^{(\t, \o, \ve, \Sigma)}(T,0) }{\ve}   \leq \cP_{j^{n}_{\t} ( \g_o)} (u_o(\t))\ .\eeq
\end{theo}
\begin{pf} Let $ \needle^{(\t, \o, \ve_o, \Sigma)}(\g_o)$  be a fixed  needle variation
and $\cZ:(0, \ve_o] \to \bR$  the function   defined by 
\beq  \cZ(\ve) \=  \int_0^T \left( \int_0^1    \left(Y^{(\ve)a} \frac{\p \cP_{\s^{(s)}_t} }{\p u^a} \Bigg|_{u^{(\ve, s)}(t)} -  \frac{\p^2  \mu'{}^{(\t, \o, \ve, \Sigma)}}{\p t\,\p s}\bigg|_{(t,s)}\right)ds \right) dt\ , \eeq
where, as usual, $\s^{(s)}_t $ denotes the curve of jets  $\s^{(s)}_t  = j^n_t(\g^{(U(\ve,s))})$ and the $Y^{(\ve)a} $ are the components of the vector field $Y_{ F^{(\t, \o, \ve, \Sigma)}}$ in the directions of the $u^a$-axes.  
We recall that, 
due to the particular construction of the homotopies  $ F^{(\t, \o, \ve, \Sigma)}$,  
\beq \label{limited} Y^{(\ve)a}(t,s) = \frac{\p  u^{(\ve, s)a}(t)}{\p s}  = \frac{\p \left((1-s) u^a_o(t) + s \check u^{(\t, \o, \ve)a}(t)\right)}{\p s} = \check u^{(\t, \o, \ve)a}(t) - u^a_o(t)\eeq
and $\left(\check u^{(\t, \o, \ve, \Sigma)}-  u_o\right)\big|_{[0,T] \setminus [\t - \ve - k \ve^2, \t + k \ve^2]} = 0 $. Thus 
the  functions  $Y^{(\ve)}{}^a|_{{u^{(s, \ve)}(t)}}$  are  equal to $0$ at the points  outside of the rectangle  $ [\t - \ve - \frac{k\ve^2}{2}, \t + \frac{k \ve^2}{2}] \times [0,1]$. 
By Theorem \ref{corollone},  a necessary condition for $ \g_o$ to be a solution to the Mayer problem is that $\cZ(\ve) \leq 0$ for any $\ve \in (0, \ve_o]$ . Hence we have that for any such $\ve$ 
\begin{align}
\label{step1} 0 & \geq  \frac{1}{\ve} \cZ(\ve) =  \frac{1}{\ve} \Bigg( \iint_{[0, T] \times [0,1]}    Y^{(\ve)a} \frac{\p \cP_{\s^{(s)}_t} }{\p u^a} \Bigg|_{u^{(\ve, s)}(t)} dt ds -\\
\nonumber & \hskip 5 cm   - \iint_{[0, T] \times [0,1]}   \frac{\p^2  \mu'{}^{(\t, \o, \ve, \Sigma)}}{\p t\,\p s}\bigg|_{(t,s)} dt ds \Bigg) \ .
\end{align}
From \eqref{limited} the first double integral  in \eqref{step1}  reduces to the sum
\begin{multline*}
 \int_0^1 \int_{\t - \ve - k \ve^2}^{\t- \ve} \!\!   Y^{(\ve)a} \frac{\p \cP_{\s^{(s)}_t} }{\p u^a} \Bigg|_{u^{( \ve, s)}(t)} dt ds +  \int_0^1 \int_{\t - \ve}^{\t }    Y^{(\ve)a} \frac{\p \cP_{\s^{(s)}_t} }{\p u^a} \Bigg|_{u^{( \ve, s)}(t)} dt ds {+} \\
 {+} \int_0^1 \int_{\t}^{\t + k \ve^2 }  \!\!  Y^{(\ve)a} \frac{\p \cP_{\s^{(s)}_t} }{\p u^a} \Bigg|_{u^{( \ve, s)}(t)} dt ds\,  {=} \int_0^1 \int_{\t - \ve}^{\t }    (\o^a - u^a_o(t)) \frac{\p \cP_{\s^{(s)}_t} }{\p u^a} \Bigg|_{u^{( \ve, s)}(t)} dt ds + O(\ve^2)
 \ .
 \end{multline*}
From this, we obtain 
\begin{align}
\nonumber 0 \geq &  \frac{1}{\ve} \cZ(\ve) = \int_0^1 \frac{1}{\ve} \left(  \int_{\t - \ve}^\t \left( \o^a - u^a_o(t)\right) \frac{\p \cP_{\s^{(s)}_t} }{\p u^a} \Bigg|_{u^{(\ve, s)}(t)} dt   \right) ds    -\\
\nonumber & \hskip 6 cm  - \frac{\mu'{}^{(\t, \o, \ve, \Sigma)}(T, 1) - \mu'{}^{(\t, \o, \ve, \Sigma)}(T,0)  }{\ve}+   O( \ve)  = \\
\nonumber &  {=} \int_0^1  \left( \o^a - u^a_o(\t )\right) \frac{\p \cP_{\s^{(s)}_{\t }} }{\p u^a} \Bigg|_{u^{( \ve, s)}(\t )}  \hskip - 0.7cm  ds -\frac{\mu'{}^{(\t, \o, \ve, \Sigma)}(T, 1) - \mu'{}^{(\t, \o, \ve, \Sigma)}(T,0)  }{\ve}+   O( \ve) =
\end{align}
\begin{align}
\nonumber & = \int_0^1 \frac{\p u^{(\ve, s)a}}{\p s}  \bigg|_{t = \t } \frac{\p \cP_{\s^{(s)}_{\t }} }{\p u^a}\Bigg|_{u^{( \ve, s)}(\t )}  \hskip - 0.7cm  ds  - \frac{\mu'{}^{(\t, \o, \ve, \Sigma)}(T, 1) - \mu'{}^{(\t, \o, \ve, \Sigma)}(T,0)  }{\ve}+   O( \ve) = \\
  \label{derivative}   
  & =  \cP_{j^{n}_{\t} ( \g^{( U^{(\t, \o,  \ve)})})}(\o) {-}   \cP_{j^{n}_{\t} ( \g^{( U_o)})}(u_o(\t)) -\frac{\mu'{}^{(\t, \o, \ve, \Sigma)}(T, 1) - \mu'{}^{(\t, \o, \ve, \Sigma)}(T,0)  }{\ve}+   O( \ve) 
  \ .
\end{align}
From this and \eqref{firstlimit-bis}  the claim follows. 
\end{pf}
\par
In order to have a  truly helpful theorem,   the previous result should be combined with some  efficient way to determine  the sign of the  corrective term  $\liminf_{\ve \to 0^+} \frac{\mu'{}^{(\t, \o, \ve, \Sigma)}(T, 1) - \mu'{}^{(\t, \o, \ve, \Sigma)}(T,0)  }{\ve}$ appearing in \eqref{712}.  This can be reached exploiting  the next technical lemma. \par
\begin{lem} \label{lemmarisolu} Let  the triple $(\cK, L, C)$ and the curve   $ \g_o \=  \g^{( U_o)}$ be  as  in  Theorem \ref{POMPone}. 
For any  $\needle^{(\t, \o, \ve_o, \Sigma)} (\g_o)= \{ F^{(\t, \o, \ve, \Sigma)}, 0 < \ve \leq \ve_o\}$ of  $\g_o $, one has that 
\begin{multline} \label{717}\mu'{}^{(\t, \o, \ve, \Sigma)}(T, 1) - \mu'{}^{(\t, \o, \ve, \Sigma)}(T,0)  =   C\big|_{{F}^{ (\t, \o, \ve, \Sigma)(n)}(T,1)} -  C\big|_{\g^{(n)}_o(T)} - \\
 - \int_0^T L\big|_{F^{  (\t, \o, \ve, \Sigma)(n)}(t,1)} dt  +  \int_0^T L\big|_{\g^{(n)}_o(t)} dt + \\
 + \int_0^1  \sum_{\d = 1}^r \sum_{\ve = 0}^{\d-1} (-1)^{\ve}  \frac{d^\ve}{dt^\ve}  \left(\frac{\p L  }{\p q^i_{(\d)}} \right) 
(Y_{F})^{i}_{{(\d-(\ve+1))}}  \big|_{{F}^{ (\t, \o, \ve, \Sigma)(n)}(T,s)} ds - \\
-  \int_0^1  \sum_{\d = 1}^r \sum_{\ve = 0}^{\d-1} (-1)^{\ve}  \frac{d^\ve}{dt^\ve}  \left(\frac{\p L  }{\p q^i_{(\d)}} \right) 
(Y_{F})^{i}_{{(\d-(\ve+1))}}  \big|_{{F}^{ (\t, \o, \ve, \Sigma)(n)}(0,s)} ds
\ . \end{multline}
\end{lem}
\begin{pf}
 Let $\Pi$ and $\Pi^{PC}$ be the $1$-forms on $J^n(\wh \cQ|\bR)$ defined by
 \begin{align}
\nonumber & \Pi \=  \bigg\{\frac{1}{2}\!\!\!\sum_{\smallmatrix 1 \leq i \leq N\\
   0 \leq  \b \leq r-1
   \endsmallmatrix} \!\!\! \bigg( ( \pa^i_{\b(1)})^2 {-} ( \pa'{}^i_{\b(2)})^2 {-}  ( \pa''{}^i_{\b(2)})^2 
    \bigg) + \\
& \hskip 3 cm   + \!\!\! \sum_{\smallmatrix 1 \leq i \leq N\\
   0 \leq  \b \leq r-1
   \endsmallmatrix} \!\!\! \bigg( \frac{1}{2} ( \pa^i_{\b})^2+ \frac{\pi^4}{32 T^4} ( \pa'{}^i_{\b})^2 +  \frac{\pi^4}{32 T^4}( \pa''{}^i_{\b})^2  \bigg)   \bigg\} dt\ ,\\
\nonumber &\Pi^{PC} :=    \Pi 
+   \l\hskip-10pt  \sum_{\smallmatrix 1 \leq i \leq N\\
   0 \leq  \b \leq r-1
   \endsmallmatrix} \bigg( \pa^i_{\b(1)} \varpi^i_{\b(0)}  -  \pa'{}^i_{\b(2)} \varpi'{}^i_{\b(1)} -  \pa''{}^i_{\b(2)} \varpi''{}^i_{\b(1)} + \\
   & \hskip 7 cm +     \pa'{}^i_{\b(3)} \varpi'{}^i_{\b(0)} +  \pa''{}^i_{\b(3)} \varpi''{}^i_{\b(0)} \bigg)
\end{align}
   (for the definitions of the $1$-forms $\varpi_{\b(\d)}^i$ etc., see \eqref{holonomic}). Consider  the   variation   $\wh F = \wh F^{(\t, \o, \ve, \Sigma)}$ 
   in the extended space $\bR \times \wh \cQ$,  which is uniquely associated with one of the    $F = F^{(\t, \o, \ve, \Sigma)}$ belonging 
   to the considered needle variation, as we indicated at the beginning of Sect.\ \ref{section7.2}. Let also denote by  $X_{\wh{F}}$ and $Y_{\wh{F}}$, 
   the corresponding vector fields in $\wh \cS = \wh F([0,T] \times [0,1])$,  as defined in \eqref{vectorX-A} and \eqref{vectorY-A}. Since  
 $\Pi$ and $\Pi^{PC}$ are variationally equivalent (see definition in Sect.\ \ref{variationallyequiv}) and, at each point of $\wh \cS$,  the vector field $X_{\wh{F}}$  is   tangent to  curves of jets of  one of the curves in $\wh \cQ$ determined by the homotopy $\wh F$,  we have   
\begin{multline} \label{720}
- \int_0^T  \Pi^{PC} \left(X_{\wh{F}}\right)|_{{\wh{F}}^{(n)}(t, 1)}  dt  + \int_0^T  \Pi^{PC} \left(X_{\wh{F}}\right)|_{{\wh{F}}^{(n)}(t, 0)}  dt  =\\
=  - \int_0^T  \Pi \left(X_{\wh{F}}\right)|_{{\wh{F}}^{(n)}(t, 1)}  dt  + 
+  \int_0^T  \Pi \left(X_{\wh{F}}\right)|_{{\wh{F}}^{(n)}(t, 0)}  dt  = \\
  =  \mu^{(\t, \o, \ve, \Sigma)}(T, 1) - \mu^{(\t, \o, \ve, \Sigma)}(T,0)    + \int_0^T L\big|_{\wh F^{(n)}(t,1)} dt  -  \int_0^T L\big|_{\g^{(n)}_o(t)} dt  \ .
\end{multline}
On the other hand, by the Stokes Theorem 
\begin{multline} \label{pollini-bisbis} 
- \int_0^T  \Pi^{PC} \left(X_{\wh{F}}\right)|_{{\wh{F}}^{(n)}(t, 1)}  dt  + \int_0^T  \Pi^{PC} \left(X_{\wh{F}}\right)|_{{\wh{F}}^{(n)}(t, 0)}  dt  
= \int_0^1 \Pi^{PC} \left( Y_{\wh F}\right)\big|_{\wh F^{(n)} (0,s)} ds - \\
- \int_0^1 \Pi^{PC} \left(Y_{\wh F}\right) \big|_{\wh F^{(n)} (T,s)}ds 
+\iint_{[0,T] \times [0,1]}d \Pi^{PC}(X_{\wh{F}}, Y_{\wh{F}})\big|_{{\wh{F}}^{(n)}(t,s)} dt ds \end{multline}
We  claim that the third summand  in the right hand side of  \eqref{pollini-bisbis} is $0$. Indeed,  by the same arguments of the proof of Lemma \ref{pcprop},  the differential $d \Pi^{PC}$ is equal to a sum of $2$-forms that are (a) either identically vanishing on the vector field $X_{\wh F}$    or (b)
have coefficients that vanish identically along the solutions of the Euler-Lagrange equations determined by  
 \begin{multline*}L' =  \bigg\{\frac{1}{2}\!\!\!\sum_{\smallmatrix 1 \leq i \leq N\\
   0 \leq  \b \leq r-1
   \endsmallmatrix} \!\!\! \bigg( ( \pa^i_{\b(1)})^2 {-} ( \pa'{}^i_{\b(2)})^2 {-}  ( \pa''{}^i_{\b(2)})^2 
    \bigg) + \\
  + \!\!\! \sum_{\smallmatrix 1 \leq i \leq N\\
   0 \leq  \b \leq r-1
   \endsmallmatrix} \!\!\! \bigg( \frac{1}{2} ( \pa^i_{\b})^2+ \frac{\pi^4}{32 T^4} ( \pa'{}^i_{\b})^2 +  \frac{\pi^4}{32 T^4}( \pa''{}^i_{\b})^2  \bigg)   \bigg\}\ .
   \end{multline*}
 Since we are integrating along the points of the surface $\wh \cS$ (whose components $\pa^i_{\b}(t)$, $\pa'{}^i_{\b}(t)$, $\pa''{}^i_\b$ are   solutions  precisely to such Euler-Lagrange equations),  the claim follows.  From this,  \eqref{720} and \eqref{pollini-bisbis},   we  obtain 
\beq   \label{pollini-bis} 
\begin{split} \mu'{}^{(\t, \o, \ve, \Sigma)}(T, 1) & - \mu'{}^{(\t, \o, \ve, \Sigma)}(T,0)    +  \int_0^T L\big|_{\wh F^{(n)}(t,1)} dt  -  \int_0^T L\big|_{\g^{(n)}_o(t)} dt = \\
&= \int_0^1 \Pi^{PC} \left(  Y_{\wh F} \right) \big|_{{\wh{F}}^{(n)}(0,s)} ds - \int_0^1\Pi^{PC} \left(  Y_{\wh F} \right) \big|_{{\wh{F}}^{(n)}(T,s)} ds + \\
& \hskip 1 cm +  \sum_{\b = 1}^{r-1}  \sum_{i= 1}^N  \int_0^1 \left(    \pa'{}^i_{\b(3)} Y'{}^i_{\b(0)} + \pa''{}^i_{\b(3)} Y''{}^i_{\b(0)}  \right)\bigg|_{{\wh{F}}^{(n)}(T, v)} d v\ .
\end{split}
\eeq
We now recall  that the initial and terminal conditions on the functions   $\pa^i_{\b}(t)$, $\pa'{}^i_{\b}(t)$, $\pa''{}^i_\b$  have been selected  such a way  that 
the  three terms  of the right hand side of \eqref{pollini-bis}  are  equal to minus  the corresponding integrals along the two  ``vertical sides'' of  $\p \wh \cS$ of  the $1$-form
\beq 
 \sum_{\d = 1}^r \sum_{\z = 0}^{\d-1} (-1)^{\z}  \frac{d^\z}{dt^\z}  \left(\frac{\p  \left(L + \frac{ dC}{dt}\right) }{\p q^i_{(\d)}} \right) \o^i_{(\d-(\z+1))} 
\eeq
(see  Remarks \ref{remark52} and \ref{remark53}).  From this and  the fact that the cost function $C$ vanishes identically on $J^n(\cQ|\bR)_{ t = 0}$, 
it follows that 
\beq   \label{pollini-terter} 
\begin{split} \mu'{}^{(\t, \o, \ve, \Sigma)}(T, 1)  &- \mu'{}^{(\t, \o, \ve, \Sigma)}(T,0)   = - \int_0^T L\big|_{\wh F^{(n)}(t,1)} dt  +  \int_0^T L\big|_{\g^{(n)}_o(t)} dt + \\
&\hskip 0.5 cm  +  \int_0^1  \sum_{\d = 1}^r \sum_{\z = 0}^{\d-1} (-1)^{\z}  \frac{d^\z}{dt^\z}  \left(\frac{\p L  }{\p q^i_{(\d)}} \right) 
Y_{\wh F}{}^{i}_{{(\d-(\z+1))}}  \big|_{{\wh{F}}^{(n)}(T,s)} ds - \\
 &\hskip 0.5 cm -  \int_0^1  \sum_{\d = 1}^r \sum_{\z = 0}^{\d-1} (-1)^{\z}  \frac{d^\z}{dt^\z}  \left(\frac{\p L  }{\p q^i_{(\d)}} \right) 
Y_{\wh F}{}^{i}_{{(\d-(\z+1))}}  \big|_{{\wh{F}}^{(n)}(0,s)} ds + \\
& \hskip 0.5 cm + \int_0^1  \sum_{\d = 1}^r \sum_{\z = 0}^{\d-1} (-1)^{\z}  \frac{d^\z}{dt^\z} \frac{\p  }{\p q^i_{(\d)}}\left( \frac{ dC}{dt} \right)
Y_{\wh F}{}^{i}_{{(\d-(\z+1))}}  \big|_{{\wh{F}}^{(n)}(T,s)} ds. 
\end{split}
\eeq
We now observe that the fifth summand in \eqref{pollini-terter}  is equal to  the sum along  the two ``vertical sides'' of $\partial \wh \cS$ of the $1$-form 
$$\b^{PC} \= \frac{d C}{dt} dt + \sum_{\d = 1}^r \sum_{\z = 0}^{\d-1} (-1)^{\z}  \frac{d^\z}{dt^\z} \frac{\p  }{\p q^i_{(\d)}}\left( \frac{ dC}{dt} \right)\o^i_{(\d-(\z+1))} \ . $$
This $1$-form is variationally equivalent to the $1$-form $\b \= \frac{d C}{dt} dt $.  Hence, by Stokes Theorem and the properties of the variationally equivalent $1$-forms
\begin{multline} \label{ultimino} \int_0^1  \sum_{\d = 1}^r \sum_{\z = 0}^{\d-1} (-1)^{\z}  \frac{d^\z}{dt^\z} \frac{\p  }{\p q^i_{(\d)}}\left( \frac{ dC}{dt} \right)
Y_{\wh F}{}^{i}_{{(\d-(\z+1))}}  \big|_{{\wh{F}}^{(n)}(T,s)} ds = \\
=  \int_0^T  \left( \frac{d C}{dt}  \bigg|_{{\wh{F}}^{(n)}(t,1)} -  \frac{d C}{dt}  \bigg|_{{\wh{F}}^{(n)}(t,0)}\right) dt +  \iint_{[0,T] \times [0,1]} d \b^{PC} (X_{\wh F}, Y_{\wh F})\bigg|_{\wh F^{(n)}(t,s)} dt ds = 
\\
= C\big|_{{\wh{F}}^{(n)}(T,1)} -  C\big|_{\g^{(n)}_o(T)}  +   \iint_{[0,T] \times [0,1]} d \b^{PC} (X_{\wh F}, Y_{\wh F})\bigg|_{\wh F^{(n)}(t,s)} dt\, ds  \ .
\end{multline}
Using once again  the  proof of Lemma \ref{pcprop}, and that  the curves ${\wh{F}}^{(n)}(\cdot, s)$, $s \in [0,1]$, have tangent vectors on which the holonomic $1$-forms vanish identically, we obtain that the double integral in \eqref{ultimino} reduces to the integral of a linear combination of $2$-forms with coefficients given by 
the Euler-Lagrange operator applied to the Lagrangian $\frac{d C}{dt}$. By the well-known property that  an  Euler-Lagrange operator  on a    total differential  gives an identically vanishing function,  we  conclude that  $  \iint_{[0,T] \times [0,1]} d \b^{PC} (X_{\wh F}, Y_{\wh F})\big|_{\wh F^{(n)}(t,s)} dt\, ds = 0$.  From \eqref{ultimino} and \eqref{pollini-terter} the lemma follows.
 \end{pf}
\begin{cor}[Generalised Pontryagin Maximum Principle -- II Version] \label{Corollone}  Let   $(\cK, L, C)$  be as in  Theorem \ref{POMPone} and for any $\cK$-controlled curve   $ \g \=  \g^{( U)}$    denote by $\Gneedles(\g)$   the  class of   needle  variations  $\needle^{(\t, \o, \ve_o, \Sigma)} (\g_o)= \{ F^{(\t, \o, \ve, \Sigma)}, 0 < \ve \leq \ve_o\}$ of  $\g_o $, 
which satisfy  the following inequality  for any $\ve \in (0, \ve_o]$
\begin{multline} \label{additional}     \int_0^T \Bigg(L\big|_{F^{  (\t, \o, \ve, \Sigma)(n)}(t,1)}   -   L\big|_{\g^{(n)}_o(t)} \Bigg)dt  
 + \int_0^1 \Bigg(- \frac{\p C}{\p q^i_{(\b)}} (Y_F)^i_{(\b)}\bigg|_{{F}^{ (\t, \o, \ve, \Sigma)(n)}(T,s)} -  \\
-  \sum_{\d = 1}^r \sum_{\ve = 0}^{\d-1} (-1)^{\ve}  \frac{d^\ve}{dt^\ve}  \left(\frac{\p L  }{\p q^i_{(\d)}} \right) 
(Y_{F})^{i}_{{(\d-(\ve+1))}}  \big|_{F^{ (\t, \o, \ve, \Sigma)(n)}(T,s)} \Bigg) ds + \\
 + \int_0^1  \sum_{\d = 1}^r \sum_{\ve = 0}^{\d-1} (-1)^{\ve}  \frac{d^\ve}{dt^\ve}  \left(\frac{\p L  }{\p q^i_{(\d)}} \right) 
(Y_{F})^{i}_{{(\d-(\ve+1))}}  \big|_{F^{ (\t, \o, \ve, \Sigma)(n)}(0,s)} ds \geq 0
\ . \end{multline}
A $\cK$-controlled curve 
$ \g_o = \g^{(U_o)}$ is   a solution   to the  Mayer  problem  only if  for any    needle variation in  $\Gneedles(\g_o)$
  \beq  \label{712-bis}  \cP_{j^{r}_{\t} ( \g_o)}(\o)   \leq \cP_{j^{n}_{\t} ( \g_o)} (u_o(\t))\ .\eeq
\end{cor}
\begin{pf} From Lemma \ref{lemmarisolu},  if  a needle variation satisfies \eqref{additional}, then  the expression $- \frac{\mu'{}^{(\t, \o, \ve, \Sigma)}(T, 1) - \mu'{}^{(\t, \o, \ve, \Sigma)}(T,0)  }{\ve}$ is non-negative for any $\ve \in (0, \ve_o]$.  From this and   Theorem \ref{POMPone}, the necessary condition \eqref{712-bis} holds. 
\end{pf}
\par
\smallskip
\subsection{The classical  Pontryagin Maximum Principle}\hfill\par
\label{classicalPMP}
Consider now a classical Mayer problem, i.e. a problem as described in Sect.\ \ref{section211}.  Let us  represent the controlled evolutions of the system by curves 
$x(t) = (x^1(t), \ldots x^{N'}(t))$, $t \in [0, T]$,  in some $\bR^{N'}$. They are   constrained by the  conditions: 
\begin{itemize}[leftmargin = 25pt]
\item[(a)] $x(0) = x_o$ for a fixed initial value $x_o$;  
\item[(b)] they satisfy the differential constraints $\frac{d x^i}{dt}  = f^i(t, x(t), u(t))$.
\end{itemize}
  As    pointed out in Sect.\ \ref{Cardinuccio},   if we   add the  auxiliary variables 
$p = (p_1, \ldots, p_{N'})$,   impose that they  are solutions to the   equations 
$ \frac{d p_i}{dt} = - p_\ell  \frac{\p f^\ell}{\p x^i}(t, x(t), u(t))$ 
and   set 
$$q^1 \= x^1\ ,\ \ldots\ ,\ q^{N'} \= x^{N'}\\  , \ q^{N' + 1} \= p_1, \ \ldots\ , \ q^{2N'} \= p_{N'} \ ,$$
 such a classical  problem can be considered as a generalised  Mayer problem on $J^3(\bR^N|\bR)$, $N = 2 N'$ (\footnote{We consider a jet bundle of order $3$ just to be sure that \eqref{crucialineq} is satisfied by the actual order $r = 1$ of the controlled Lagrangian  defined below.})   with  the  defining triple given by: 
 \begin{itemize}[itemsep=5pt, leftmargin=10pt]
\item   the set $\cK$  of  the  pairs $U = (u(t), \s)$,  in which $u(t)$ is a  
smooth curve  $u: [0,T]   \to K\subset \bR^M$  and   $\s = (A^i =  x^i(0), B_\ell = p_\ell (0))$  is a $0$-th order jet 
where $x(0)= x_o$ and   $p(0)$ is (provisionally) unconstrained. 
\item  the  controlled Lagrangian  
$$L(t, q^j_{(\b)}, u^a) \=   p_i  \left(x^i_{(1)}  - f^i(t, x^i, u^a) \right)\ . $$ 
\item a  cost function  $C: J^3(\bR^N|\bR) \to \bR$    which is of actual order $r = 0$   and coincides on $J^3(\bR^N|\bR)|_{t = T}$  with a classical  terminal cost function,  depending just on the coordinates $x^i$. With no loss of generality, we  assume that $C$ depends just  on the  $x^i$ at all points. 
\end{itemize}
Such  a (generalised) Mayer problem is manifestly of normal type, as defined in Sect.\ \ref{defnormal},  and Corollary  \ref{Corollone}  applies.  Let us therefore determine what are the needle variations of the class $\Gneedles$  for this setting.  First of all, we observe that along any solution of the controlled Euler-Lagrange equations of this problem, the function  $L$ vanishes identically. Hence, since the actual order $r$ of the Lagrangian is $r = 1$,  the characterising inequality \eqref{additional} reduces to  
\beq  \label{additional-bis}  C\big|_{{F}^{ (\t, \o, \ve, \Sigma)(n)}(T,1)} -  C\big|_{\g^{(n)}_o(T)} \leq \\
 - \int_0^1  p_i 
\frac{\p x^i}{\p s} \bigg|_{F^{(\t, \o, \ve, \Sigma)}(T,s)} ds +  \int_0^1  p_i 
\frac{\p x^i}{\p s} \bigg|_{F^{(\t, \o, \ve, \Sigma)}(0,s)} ds
\ . 
\eeq
Here, we denoted by $x^i = x^i(t, s)$ the component  in the  $x^i$-direction of the $\cK$-controlled variation $F(t,s) = F^{(\t, \o, \ve, \Sigma)}(t,s)$ of a considered generalised needle variation.  We now observe that 
$$ C\big|_{{F}^{ (\t, \o, \ve, \Sigma)(n)}(T,1)} -  C\big|_{\g^{(n)}_o(T)} = \int_0^1   \frac{\p C}{\p x^i}
\frac{\p x^i}{\p s} \bigg|_{F^{(\t, \o, \ve, \Sigma)}(T,s)} ds$$
 and  $x^i(0, s) \equiv x^i_o$  so that $\frac{\p x^i}{\p s} \big|_{F^{(\t, \o, \ve, \Sigma)}(0,s)} = 0$.
Thus, \eqref{additional-bis} is equivalent to 
\beq  \label{additional-bisbis}  -  \int_0^1  \left( \frac{\p C}{\p x^i} + p_i\right) \frac{\p x^i}{\p s} \bigg|_{F^{(\t, \o, \ve, \Sigma)}(T,s)}  ds \geq 0
\ . 
\eeq
We now recall that  we are free to impose any initial  condition on  the auxiliary variables $p_j$. Furthermore, by the particular form  of the differential constraints on the curves $p_j(t)$,  
it is certainly possible to determine  a family $\Sigma = \Sigma(\ve, s)$ of initial data for a needle variation, with the property  that  the corresponding functions $p_j(t, s)$ of the $\cK$-controlled curves of 
a homotopy $F^{(\t, \o, \ve, \Sigma)}$  satisfy the terminal conditions  
\beq \label{questionable-2} p_i(T, s) =  -  \frac{\p C}{\p x^i}\bigg|_{x(T,s)}\ , \eeq
provided,  of course,  that such terminal conditions are satisfied in the first place  by the components  $p_{io}(T)  = p_i(T, 0)$  of $\g_o(T)$. 
 (see also \cite[Sect.\ 5.1]{CS}). 
This  and \eqref{additional-bisbis}    has  the following  crucial consequence: {\it  if $\g_o$ satisfies $ p_i(T) =  -  \frac{\p C}{\p x^i}\big|_{x(T)}$,  all of the  needle variations of  $\g_o$, which are led by a $\Sigma$  forcing \eqref{questionable-2}, are  in the class $\Gneedles(\g_o)$ described in  Corollary \ref{Corollone}}. In particular, there is a needle variation $\needle^{(\t, \o, \ve_o, \Sigma)} (\g_o)$ in the   class $\Gneedles(\g_o)$  for any  choice of $\t \in (0, T]$, $\o \in K$ and $\ve_o$ sufficiently small.  Thus,  by Corollary \ref{Corollone},  a necessary condition for such a curve $\g_o = \g^{(U_o)}$ to be a solution of the Mayer problem is that  the inequality \eqref{712-bis} holds for {\it any} needle variation as above. \par
\smallskip
We finally observe that,  for the classical Mayer problem considered in this section, we have 
$ \cP_{j^{1}_{t} ( \g_o)}(u^a) = \left(p_i f^i(t, x^i, u^a) - p_i x^i_{(1)}\right)\big|_{j^1_t(\g_o)}$.
Hence, if for each $(t, x^i,  p_i) \in \bR \times \bR^{2 N'}$,  we   denote by  $\cH(t, x^i,  p_i): K \to \bR$  the  classical  Pontryagin function
$$\cH(t, x^i,  p_i)(u^a) \= \sum_{i = 1}^{N'} p_i f^i(t, x^i, u^a)\ , $$ 
for any needle variation we have  
$ \cP_{j^{1}_{t} ( \g_o)}(\o^a)    = \cH(t, x^i,  p_i)|_{ \g_o(t)} (\o^a)- \left( p_i x^i_{(1)}\right)\big|_{j^1_{t}(\g_o)}$.
 From this,  we immediately derive   the  following version of the classical  Pontryagin Maximum Principle. 
\begin{cor}[Pontryagin Maximum Principle]
 A  $ \wt \cK$-controlled curve   $ \g_o \=  \g^{( U_o)}$,   with components $p_j(t)$ satisfying  $ p_i(T) =  -  \frac{\p C}{\p x^i}\big|_{x(T)}$,  is    a solution to the      Mayer  problem  determined by the above described triple  $( \wt \cK,  L, C)$ only if  for any $\t \in (0, T]$ and  $\o \in K \subset \bR^M$  the following inequality holds:  
 \beq   \label{PMPclassical} \cH|_{\g_o(\t)}(\o) \leq \cH|_{ \g_o(\t)}(u_o)  \ . \eeq
\end{cor}
\par
\medskip
  \section{A discussion of  a  basic  example and  some concluding remarks}\label{Sect8}
  \ \\[-1.5 cm]
 \subsection{Comparison of  different approaches to an elementary problem: \\ \phantom{aaa} the controlled linearised pendulum}\hfill\par
\label{sect8.1}
Let $x(t) \in \bR$, $t \in [0,T]$, be the coordinate which describes the position in time of 
a linearised pendulum controlled by a force $u(t)$, that is a dynamical system   subjected to the differential constraint $\ddot x(t)= - x(t) + u(t)$. 
Assume also that  the force $u(t)$ is bound  to take values in  $K = [-1, 1] \subset \bR$ and that  the  initial position for the pendulum  is constrained to be  $x(0)= 0$ and  the initial velocity  $\dot x(0)$ is   bound to be in $[-v_{\text{max}}, v_{\text{max}}]$ for some  $v_{\text{max}} > 0$.  
We  want to  discuss the Mayer problem  corresponding to finding  a time-dependent force $u(t)$, under which 
 the position $x(t)$ at  $t = T$ is maximised or, equivalently,  the terminal  cost  $C(x(T)) \= - x(T)$ is minimised. \par
\smallskip
The classical approach to such  a problem is the following.  First, let us reduce the differential constraint to a system of first order. This can be done by introducing an auxiliary variable, i.e. by representing the dynamical system with  curves  $(x^1(t), x^2(t)) \in \bR^2$ with    $x^1(t) = x(t)$ and   $x^2(t) = \dot x^1(t)$. In this way the evolution of the system is described by curves 
$ (t, x^1(t), x^2(t))$,  subjected to the differential constraints  
\beq \label{firsthalf} \begin{array}{l} \dot x^1(t) = x^2(t)\ ,\\[10pt]  \dot x^2(t) = - x^1(t) + u(t)\ , \end{array} \qquad x^1(0) = 0 \ .
\eeq
 In these coordinates the terminal cost   is  determined  by the function $C(x^1, x^2) = - x^1$. \par
\smallskip
Second, the Pontryagin auxiliary variables $p_1(t)$ and $p_2(t)$ are introduced and the  evolutions of the system is now described by curves $\g(t) = (t, x^1(t), x^2(t), p_1(t), p_2(t))$  in $\bR \times \cQ$,  $\cQ = \bR^4$, constrained by the  \eqref{firsthalf} and, at the same time,  by
\beq \label{secondhalf} \begin{array}{l} \dot p_1(t) = p_2(t)\ ,\\[10pt] \dot p_2(t) = - p_1(t)\ ,\end{array} \ \qquad \begin{array}{l} p_1(T) = - \frac{\p C}{\p x^1}\big|_{x^1(T)} = 1\ ,
\\[10pt] p_2(T) = - \frac{\p C}{\p x^2}\big|_{x^1(T)}  = 0\ .\end{array}
\eeq
The  \eqref{secondhalf}  and  the  \eqref{firsthalf}  are uncoupled. This  allows to determine explicitly the components $p_i(t)$ for each curve $\g(t) = (t, x^1(t), x^2(t), p_1(t), p_2(t))$. 
They are 
$$p_1(t) = \cos(T - t) \ ,\qquad p_2(T) = \sin(T - t)\ .$$
Due to this, for each given  $\g_o(t) = (t, x^1_o(t), x^2_o(t), p_{o1}(t), p_{o2}(t))$  satisfying   the above  constraints, the associated Pontryagin function $\cH: K \to \bR$ takes the form 
$$\cH(\o) = \cos(T-t) x^2_o(t) + \sin(T-t) (- x^1_o(t) + \o)\ .$$
Then, the usual  PMP implies that an  optimal control curve $u_o(t)$  must  satisfy  the following conditions:
\beq \begin{array}{ll} u_o(t) = +1  & \text{when}\ \sin(T-t) > 0\ ,\\[15pt]
u_o(t) = - 1 & \text{when}\ \sin(T - t) < 0 \ .
\end{array}
\label{oc}
\eeq
 %These constraints are so strong that basically fix the  values for the optimal control curve $u_o(t)$ at almost all points, independently on the choice of the initial speed $\dot x(0)$. 
 The meaning of this result  is the following: {\it for any  initial velocity $ v \in [-v_{\text{max}}, v_{\text{max}}]$,  if $u(t)$ is different from \eqref{oc},  there exists at least one needle modification $\wt u(t)$ of $u(t)$,  such
 that the controlled curve  determined by  the pair  $\wt U = \big(\wt u(t), \s = (x(0) = 0, \dot x(0) = v)\big)$ has a   cost which is smaller
than the cost of the curve determined by  $U = \big(u(t), \s = (x(0) = 0, \dot x(0) = v)\big)$}.  Thus,  if  an optimal control $(u_o(t), \s)$ exists,  it  must have $u_o(t)$ as in \eqref{oc}. After establishing this, by varying only the initial datum $\s$ (i.e. the initial velocity $v_o$), it is quite straightforward to find that an  optimal control exists and it is  $U = \big(u_o(t), \s = (x(0) = 0, \dot x(0) = v_{\text{max}})\big)$.
  \par
 \smallskip 
Let us now see how our results offer an alternative way to solve the same  problem. 
 Let us consider just one auxiliary variable, say $p$, and the controlled  Lagrangian on $J^5(\cQ|\bR) \times K$,  $\cQ = \bR$,  $K = [-1,1]$, of actual order $r= 2$: 
\beq L(t, x, \dot x, \ddot x, p,  u) \=  p (\ddot x + x - u)\ .\eeq
The corresponding controlled Euler-Lagrange equations  are  only  two and are 
\beq
\ddot x + x - u = 0\ ,\qquad \ddot p + p  = 0\ .\eeq
Our  problem can be now recognised as the generalised Mayer problem determined by the triple $(\cK, L, C)$, where $C: J^5(\cQ|\bR) \to \bR$ is the smooth function $C(j^5_t(x)) = - x(t)$ and $\cK$ is the set of 
pairs $U = (u(t), \s)$, given by a  control curve taking values in $K = [-1, 1] $ and  initial conditions  $\s = j^5_{t = 0}(\g) =  (x(0), \dot x(0), \ddot x(0), \ldots, p(0), \dot p(0), \ddot p(0), \ldots)$  subjected to no restrictions  except for the value $x(0) = 0$.
\par
\smallskip
In this setting, given a controlled curve $\g_o = \g^{(U_o)}$,    the condition \eqref{additional},  which characterises the needle variations  $\needle^{(\t, \o, \ve_o, \Sigma)} (\g_o)= \{ F^{(\t, \o, \ve, \Sigma)}, 0 < \ve \leq \ve_o\}$ in the class $\Gneedles(\g_o)$ is very simple and we give it in the next formula \eqref{86},  where we denote by $(x(t,s), p(t,s), \dot x(t, s), \dot p(t,s), \ldots)$ the components of  the jets homotopy $F^{(\t, \o, \ve, \Sigma)(5)}(t,s)$,  associated with  the curve  $U(s) \in \cK$
$$U(s) = \big(u(t,s), (x(0, s) = 0, p(0, s), \dot x(0, s) = 0, \dot p(0, s) , \ldots)\big)  \ . $$
It is 
\beq\label{86} \int_0^1 \bigg( \frac{\p x}{\p s}\bigg|_{(t = T, s)}  - p(T, s) \frac{\p \dot x}{\p s} \bigg|_{(t = T, s)} + \dot p(T, s) \frac{\p x}{\p s} \bigg|_{(t = T, s)} \bigg) ds \geq 0\ .\eeq
Since we are free to choose any initial condition for the variable $p$, we may  consider  a family $\Sigma = \Sigma(\ve, s)$ of initial data  so that  the corresponding functions $p(t, s)$ satisfies the terminal conditions  
\beq \label{questionable-2bis} p(T, s) = 0\ ,\qquad \dot p(T, s) = -1 .\eeq
For such $\Sigma$, the condition \eqref{86} is automatically satisfied. This means that, for any choice of $\t$, $\o$ and $\ve_o$, the class $\Gneedles(\g_o)$ is not empty and contains all needle variations satisfying \eqref{questionable-2bis}. Note also that \eqref{questionable-2bis} together with the differential constraint $\ddot p + p = 0$ completely determines the function $t \mapsto p(t,s)$  which is 
\beq \label{questionable-ter} p(t, s) = \sin(T - t)\qquad \text{for any}\ s\ .\eeq
On the other hand,  given a controlled curve $\g_o(t) = (t,  x_o(t), p_o(t))$  with $p_o(t)$ as in  \eqref{questionable-ter},  the corresponding function $\cP_{j^2(\g_o)}: K \to \bR$ is 
\beq \label{eqsol} \cP_{j^2(\g_o)}(\o) =  - \sin(T-t) \left(\ddot x + x - \o\right)\ .\eeq
From  the above observations and our Generalised Pontryagin Maximum Principle  (Corollary \ref{Corollone}), it follows that a curve $\g_o(t)$, determined by a control curve $u_o(t)$  is a solution to our Mayer problem only if it is as in \eqref{oc}, 
 as prescribed by the classical approach a la Pontryagin. As before, this  condition completely determines a (non-smooth) optimal control curve $u_o(t)$ and, by the same argument of before,  it follows that there is an optimal control and it is the pair $(u_o(t), \s_o)$ with $\s_o \= (x(0) = 0, \dot x(0) = v_{\text{max}})$.
   \begin{rem} \label{remark81} The control curve  $u_o(t)$ determined via \eqref{eqsol}  is clearly non smooth and hence not fitting the   simplifying  assumptions  adopted throughout this paper.  However, using the results of this paper and  standard approximation techniques, one can  rigorously prove that   the  above (discontinuous)  
 curve $u_o(t)$ is the only measurable control curve,  for which there exists  no needle modification corresponding to a  curve with   a strictly smaller terminal cost.  A detailed proof of this claim is given  in \cite{CGS1}. 
 \end{rem}
  We stress the fact that this second (new) approach to the addressed Mayer  problem  involves {\it  just one} auxiliary variable, instead of the three used in the classical approach. In fact, the same circle of 
ideas can easily find   solutions to the large class of similar Mayer problems with  one control variable $u \in [-1,1]$,  one  dependent variable $x(t)$, $t \in [0,T]$,  constrained by a differential problem of order $m$ of the form 
$$\sum_{\ell = 0}^m a_\ell \frac{d^\ell x}{d t^\ell} = u\ ,\qquad x(0) = \frac{d x}{dt}\bigg|_{t = 0} = \ldots =  \frac{d^{m-1} x}{dt^{m-1}}\bigg|_{t = 0} = 0$$
with constant coefficients $a_\ell$,  and the cost function $C = - x$. By considering just one auxiliary variable $p$ and  the controlled Lagrangian 
$L (t, x_{(\ell)}, p) \= p \big(\sum_{\ell = 0}^m a_\ell x_{(\ell)}  - u\big)$, 
one can find  the optimal control $u_o(t)$, $t \in [0,T]$,  of such a  Mayer problem with the same arguments of before. It is 
$u_o(t) = \text{sign}(p_o(t))$, 
where   $p_o(t)$  is the unique solution to the  differential problem 
\begin{multline}\sum_{\ell = 0}^m (-1)^\ell a_\ell \frac{d^\ell p}{d t^\ell} = 0\qquad \text{in}\ \ \ [0,T]\ \\
\text{with terminal conditions} \  p(T) = \frac{d p}{dt}\bigg|_{t = T} = \ldots =  \frac{d^{N-2} p}{dt^{N-2}}\bigg|_{t = T} = 0\ \ , \quad \frac{d^{N-1} p}{dt^{N-1}}\bigg|_{t = T} = -1\ .
\end{multline}
Of course, the same solution can be easily found also using the classical approach and the classical PMP,  provided that, instead of the above single auxiliary variable $p$, one introduces and handles $2m-1$ auxiliary variables: In fact, one needs $m-1$ auxiliary variables  to reduce the  constraint to a system of  first order   equations, and  the $m$ Pontryagin auxiliary variables $p_j$. 
\par
\medskip
We conclude showing that there exists also  a third different approach to  the  above  linearised pendulum  problem.  In fact, with no need of  auxiliary variables, this problem
can be immediately seen to be  a generalised Mayer problem  (in the sense of this paper) with differential constraints given  by the  controlled  Lagrangian 
\beq L: J^3(\bR|\bR)\times K \to \bR\ ,\qquad L(t, x, \dot x, u) \= \frac{1}{2} \dot x^2 - \frac{1}{2} x^2 + u x\ ,\qquad K = [-1,1]\ .\eeq
We now  restrict to the cases in which  $T \neq \pi k$ for any  $k \in \bN$. Under this assumption,  a direct inspection of the general solutions of the differential constraints
 shows that, for a given control  map $u:[0, T] \to K$, the  map
\beq \label{surjective} \f: \bR \to \bR\ ,\qquad \f(v) \= x^{(U)}(T)\ ,\qquad U = \big(u(t), \s = (x(0) = 0, \dot x(0) = v)\big)\ ,\eeq
is surjective (here, as usual, we denote by  $x^{(U)}(t)$ the unique solution of the differential constraints with initial datum $\s$ and control curve $u(t)$). 
Due to the particularly simple form of the cost function, this information immediately leads to the conclusion   that {\it if there exists an optimal control $U_o = (u_o(t), \s_o)$, then the  initial velocity must take an  extremal value}. However,  there is also 
an enlightening (but  longer) argument based on the generalised PMP  yielding to the same conclusion. Using the surjectivity of the map 
\eqref{surjective},  one can  see that for   $v \neq  \pm v_{\text{max}}$ and $u(t)$ arbitrary,  there exist {\it good} needle variation for any pair $(\t, \o) \in  (0, T) \times K$, provided that $\o$ is  sufficiently close to $u(\t)$. 
From this and the explicit expression of the map $\cP_{j^5_\t(\g^{(U)})}: K \to \bR$, it follows that  if $v \in (0, v_{\text{max}})$ and  if  $u(\t) \neq - \text{sign}(x(\t))$ at some $\t \in (0, T)$, then  there is  {\it an appropriate needle modification  and  an appropriate change of the  initial  velocity}, which gives  a curve with a strictly lower  terminal cost. By studying the terminal costs  determined by the control curves of the form $u(t) = -\text{sign}(x(t))$, one concludes (once again) that the desired optimal control $U_o = (u_o(t), \s_o)$ must have initial velocity $v_o = \pm  v_{\text{max}}$.\par
\smallskip
It now remains to  look for an optimal control $U_o$ in the restricted family of  pairs $U = (u(t), \s_o)$ in which  the  initial value is completely fixed $\s_o = (x(0) = 0, \dot x(0) = \pm v_{\text{max}})$. In this subclass,  
$\dot x(0) = v_{\text{max}}$ is  no longer modifiable and there is no way of constructing {\it good} needle variations.  Our Corollary \ref{Corollone}  gives no exploitable condition in this situation.  Thus, the   easiest way to conclude  is now to go back to one of the previous approaches and  determine explicitly  the   control curve $u_o(t)$ using  auxiliary variables.  Nonetheless, with a not  so  big effort (the    details are omitted just for not  weighing  the reader down), all terms     in \eqref{712}  can be made fully explicit for any  needle variation and,   from Theorem \ref{POMPone},   the  previously determined optimal control is once again derived. 
\par
\smallskip
 \subsection{Concluding remarks}\hfill\par
In Sect.\ \ref{sect8.1} we showed how a very simple case of a control problem with a differential constraint of  higher order can be   studied 
with   no need of  reducing  it to an equivalent  system of   first order constraints. The same approach can be adopted  in much  less trivial situations. 
Just to fix the ideas, consider a  control problem determined by a  differential constraint  of the third order
\beq \label{esempio-0} \dddot x = f(t, x, \dot x, \ddot x, u) \eeq
with  $u  = (u^a)$ varying in a compact set $K \subset \bR^m$ and terminal cost function $C = C(x(T), \dot x(T), \ddot x(T))$. By introducing  a dual variable $p$, the constraint \eqref{esempio-0} can be identified  as one of the Euler-Lagrange equations of
 the controlled Lagrangian 
\beq \label{esempio-1} L(t, x, \dot x, \ddot x, p, u) \= p\big( \dddot x - f(t, x, \dot x, \ddot x, u)\big)\ . \eeq 
Indeed,   the Euler-Lagrange equation corresponding to $p$ is precisely \eqref{esempio-0}, while  the Euler-Lagrange equation  corresponding to  $x$ is
\begin{multline}  \label{esempio-2} \frac{\p L}{\p x} - \frac{d}{dt} \left(\frac{\p L}{\p \dot x}\right) +  \frac{d^2}{dt^2} \left(\frac{\p L}{\p \ddot x}\right) 
-  \frac{d^3}{dt^3} \left(\frac{\p L}{\p \dddot x}\right) =\\
=  p \left(- \frac{\p f}{\p x}  +  \frac{d}{dt} \left( \frac{\p f}{\p \dot x}\right) -\frac{d^2}{dt^2} \left( \frac{\p f}{\p \ddot x}\right)\right)+ \dot p\left( \frac{\p f}{\p \dot x} - 2\frac{d}{dt} \left( \frac{\p f}{\p \ddot x}\right) \right) 
- \ddot p \frac{\p f}{\p \ddot x}   
- \dddot p = 0\ . 
\end{multline}
Given a control curve $u_o(t)$ and an initial condition $\s_o$, there is a uniquely associated controlled  curve $\g(t) =  (t, x(t), p(t))$, in which  $x(t)$ satisfies \eqref{esempio-0} and the  initial condition $  j^2_{t = 0}(x) = \s_o$, and   $p(t)$ satisfies \eqref{esempio-2} and the  terminal conditions
\beq 
\begin{split}
 & p(T) =\left. \frac{\p C}{\p \ddot x}\right|_{
j^{2}_{T}(x)}\ , \hskip 2 cm 
 \dot p(T)=  \left.-\frac{\p C}{\p \dot x}\right|_{
j^{2}_{T}(x)}
- \left. p(T)\frac{\p f}{\p \ddot x} \right|_{
(j^{2}_{T}(x), u(T))}\ ,\\
 \label{esempio-3} & \ddot p(T) =  \left.\frac{\p C}{\p x}\right|_{j^2_T(x)}  \hskip -0.4 cm +  \left. p(T)\frac{\p f}{\p \dot x} \right|_{
(j^{2}_{T}(x), u(T))}  \hskip -0.5 cm -   \left.  \dot p(T) \frac{\p f}{\p \ddot x} \right|_{
(j^{2}_{T}(x), u(T))}  \hskip -0.5 cm -   \left.p(T)  \frac{d}{dt}\left(\frac{\p f}{\p \ddot x} \right)\right|_{
(j^{2}_{T}(x), u(T))}\hskip - 0.4 cm .
\end{split}
\eeq
These conditions are the higher order analogues of the terminal condition \eqref{questionable-2} used in the classical case of first order differential constraints. In fact, they are specially designed  to  make the left hand side of \eqref{additional} equal to  $0$ for any generalised needle variation that consists of homotopies  of  curves  satisfying the above  two-point boundary problem (see \cite{CGS1}). This means that any such generalised needle variation is in   $\Gneedles (\g)$ and that the generalised PMP given by   Corollary \ref{Corollone} holds at all points of  a solution  $\g(t) = (t, x(t), p(t))$  of the above problem, in perfect analogy with  the  case of differential  constraints of  first order.   \par
We stress the fact  that,   for each control curve $u(t)$ and  associated solution $x(t)$ to  \eqref{esempio-1}, the corresponding curve   $p(t)$ is  completely determined by the terminal conditions \eqref{esempio-3} and  the {\it linear} equation \eqref{esempio-2} with time-dependent coefficients. Such  linearity of the equation  for  $p$ is  a general and quite useful fact. Indeed, it  holds {\it  for any control problem determined by an  $r$-th order  differential  system of  $n$ equations
 the form }
 \beq \label{esempio-4} \frac{d^r x^i}{d t^r} = f^i\left(t,x^j(t), \ldots, \frac{d^{r-1} x^j}{d t^{r-1}},u^a(t)\right)\ , \eeq 
provided that  such constraints are coupled with all other  Euler-Lagrange equations of the controlled Lagrangian 
 \beq \label{esempio-5} L\left(t,x^j(t), \ldots, \frac{d^{r-1} x^j}{d t^{r-1}}, p_j, u^a\right) \=  p_i\left( x^i_{(r)}  - f^i\left(t,x^j(t), \ldots, \frac{d^{r-1} x^j}{d t^{r-1}}, u^a\right)\right)\ . \eeq  
On the other hand, as  we mentioned in Remark \ref{remark81},  standard approximation techniques yield to versions of the above generalised PMP 
that hold under much weaker regularity assumptions  (\cite{CGS1}).   We   expect that   the above described  two-point boundary problems together with   such  extended versions of the PMP  
 can be    helpful also for the  direct computations of the optimal controls  under higher order constraints.  Investigations on  this aspect are left to future work. 
%%%%%%%%%%%%%%%%%%%%%%%%%%%%%%%%%%%%%%%%%%%%%%%%%%%%%%%

%%%%%%%%%%%%%%%%%%%%%%%%%%%%%%%%%%%%%%%%%%%%
\bigskip
\bigskip
\font\smallsmc = cmcsc8
\font\smalltt = cmtt8
\font\smallit = cmti8
\hbox{\parindent=0pt\parskip=0pt
\vbox{\baselineskip 9.5 pt \hsize=3.1truein
\obeylines
{\smallsmc 
Franco Cardin
Dipartimento di Matematica % Pura ed Applicata
Tullio  Levi-Civita
Universit\`a degli Studi di Padova,
Via Trieste, 63, 
I-35121 Padua
ITALY
\ 
\ 
\ 
}\medskip
{\smallit E-mail}\/: {\smalltt cardin@math.unipd.it 
\ 
}
}
\hskip 0.0truecm
\vbox{\baselineskip 9.5 pt \hsize=3.7truein
\obeylines
{\smallsmc
Cristina Giannotti \& Andrea Spiro
Scuola di Scienze e Tecnologie
Universit\`a di Camerino
Via Madonna delle Carceri
I-62032 Camerino (Macerata)
ITALY
\ 
}\medskip
{\smallit E-mail}\/: {\smalltt cristina.giannotti@unicam.it
\smallit E-mail}\/: {\smalltt andrea.spiro@unicam.it}
}
}

\begin{thebibliography}{30}


 \bibitem{AS} A. A. Agrachev and Yu. L. Sachkov, {\it Control theory from the geometric viewpoint} in  ``Encyclopaedia of Mathematical Sciences, 87. Control Theory and Optimization, II.'',  Springer-Verlag, Berlin, 2004. %xiv+412 pp.

% \bibitem{Br2} Aldo Bressan, {\it On the application of control theory to certain problems for Lagrangian systems, and hyper-impulsive motion for these. I. Some general mathematical considerations on controllizable parameters}, 
% Atti Accad. Naz. Lincei Rend. Cl. Sci. Fis. Mat. Natur., {\bf 8}(1) (1988), 91--105.

% \bibitem{Br1} Aldo Bressan, {\it On control theory and its applications to certain problems for Lagrangian systems. On hyper-impulsive motions for these. II. Some purely mathematical considerations for hyper-impulsive motions. Applications to Lagrangian systems}, Atti Accad. Naz. Lincei Rend. Cl. Sci. Fis. Mat. Natur., {\bf 8}(1) (1988),  107--118.
% 
%\bibitem{Br3} Aldo Bressan,  {\it On control theory and its applications to certain problems for Lagrangian systems. On hyper-impulsive motions for these. III. 
% Strengthening of the characterizations performed in Parts I and II, for Lagrangian systems. An invariance property}, 
%  Atti Accad. Naz. Lincei Rend. Cl. Sci. Fis. Mat. Natur., {\bf 8}(3) (1990), 461--471.
%
%
%\bibitem{Br} Alberto Bressan,  {\it Impulsive Control of Lagrangian Systems and locomotion in fluids}, Discrete Contin. Dyn. Syst., {\bf 20}(1), (2008), 1--35. 
% 
%  
\bibitem{BP}A. Bressan and B. Piccoli, Introduction to the mathematical theory of control,  {\it American Institute of Mathematical Sciences (AIMS), Springfield, MO}, 2007. %xiv+312 pp.
%
% \bibitem{Ca} F. Cardin,  {\it Elementary symplectic topology and mechanics},  Lecture Notes of the Unione Matematica Italiana, vol. 16. Springer,  2015. 
%
%\bibitem{C1} F. Cardin, {\it Franco's notes on Pontrjagin maximum principle}, unpublished (2015).
%
% \bibitem{CF1}  F. Cardin and M. Favretti, {\it On nonholonomic and vakonomic dynamics of mechanical systems with nonintegrable constraints},  Journal of Geometry and Physics, {\bf 18}(4)  (1996), 295--325.
% 
% \bibitem{CF2}  F. Cardin and  M. Favretti, {\it Hyper-impulsive motion on manifolds},  Dynam. Contin. Discrete Impuls. Systems, {\bf 4}(1) (1998), 1-21.
% 
  \bibitem{CS}  F. Cardin and  A. Spiro,  {\it Pontryagin maximum principle and {S}tokes theorem},
J. Geom. Phys.  {\bf 142},  (2019), 274--286.

  \bibitem{CGS1}  F. Cardin, C. Giannotti and  A. Spiro,  {\it On the Pontryagin Maximum Principle under differential constraints of higher order},
in preparation.


\bibitem{Ce} L. Cesari,
   Optimization---theory and applications, {\it Springer-Verlag, New York}, 1983.
% \bibitem{CGS} F. Cardin, C. Giannotti and A. Spiro, {\it Examples}, in preparation.
%
%\bibitem{Ch} B. C. Chachuat,  Nonlinear and dynamic optimization: from theory to Practice, Lecture Notes for IC-32: Winter Semester,  
%{\it \`Ecole Polytechnique F\'ed\'erale de Lausanne}, 2006.
%
%\bibitem{Ev} L. C. Evans, Partial Differential Equations - $2^{\text{nd}}$ edition, American Mathematical Society, Providence, 2010.
%
%
\bibitem{FS} E. Fiorani and A. Spiro
{\it Lie algebras of conservation laws of variational ordinary
              differential equations},
 {J. Geom. Phys.},
  {\bf 88}
  (2015),
56--75.

\bibitem{Ju} V.  Jurdjevic,  Geometric control theory,  {\it Cambridge University Press, Cambridge}, 1997.
%
%\bibitem{Le} U. Ledzewicz and H. Sch\"attler, {\it An extended maximum principle}, 
%Nonlinear Anal. {\bf 29} (1997), 159--183.
%
%\bibitem{Le1} U. Ledzewicz and H. Sch\"attler, {\it A high-order generalized local maximum principle},
%SIAM J. Control Optim.,  {\bf 38}, (2000), 823--854.
%
% \bibitem{Ma1} C.-M. Marle, {\it Sur la g\'eom\'etrie des syst\`emes m\'ecaniques \`a liaisons actives},  C. R. Acad. Sci. Paris S\'er. I Math.,  {\bf 311} (1990), 839--845. 
%
% \bibitem{Ma2} C.-M. Marle,  {\it G\'eom\'etrie des syst\`emes m\'ecaniques \`a liaisons actives} in   ``Symplectic geometry and mathematical physics (Aix-en-Provence, 1990),'',  {\it Birkh\"auser Boston, Boston, MA}, 1991.
% 
% \bibitem{MBLP}  E.  Massa, D. Bruno, G. Luria and E. Pagani, {\it Geometric constrained variational calculus. I: Piecewise smooth extremals}, 
% Int. J. Geom. Methods Mod. Phys. {\bf 12} (2015), 1550061, 42 pp. 
% 
%  \bibitem{MBLP1} E.  Massa, D. Bruno, G. Luria and E. Pagani, {\it Geometric constrained variational calculus. II: The second variation (I Part)}, 
% Int. J. Geom. Methods Mod. Phys. 13 (2016), no. 1, 1550132, 31 pp.
%
% \bibitem{MBLP2}  E.  Massa, D. Bruno, G. Luria and E. Pagani, {\it Geometric constrained variational calculus. II: The second variation (II Part)}, 
%  Int. J. Geom. Methods Mod. Phys. 13 (2016), no. 4, 1650038, 39 pp. 
%
% \bibitem{MT} M. R. Menzio and  
%W. M. Tulczyjew, 
%{\it Infinitesimal symplectic relations and generalized Hamiltonian dynamics},
%Ann. Inst. Henri Poincar\'e {\bf XXVIII} (1978), 349--367.
%
 \bibitem{PBGM} L. S. Pontryagin, V. G. Boltyanskii,  R. V.  Gamkrelidze and E. F.  Mishchenko, The mathematical theory of optimal processes (translated from Russian by D. E. Brown), {\it A Pergamon Press Book, The Macmillan Co., New York}, 1964.
% 
% 
%  \bibitem{PR}P.-D. Prieto-Martinez, N. Roman-Roy,
%Lagrangian-Hamiltonian unified formalism for autonomous higher order dynamical systems. 
%J. Phys. A 44 (2011), no. 38,  35 pp. 
% 
% 
% 
% \bibitem{Ra1} F. Rampazzo, {\it On Lagrangian systems with some coordinates as controls}, Atti Accad. Naz. Lincei, Classe di Scienze Mat. Fis. Nat. Serie 8, {\bf 82} (1988), 685--695.
% 
% \bibitem{Ra2} F. Rampazzo, {\it On the Riemannian structure of a Lagrangean system and the problem of adding time-dependent coordinates as controls}, European J. Mechanics A/Solids, {\bf 10} (1991), 405--431.
% 
\bibitem{Sa} D.J. Saunders,
 {\rm The geometry of jet bundles},
 {\it Cambridge University Press, Cambridge}, 1989.
 
 \bibitem{Sp} A. Spiro, {\it Cohomology of Lagrange complexes invariant under pseudogroups of local transformations},
Int. J. Geom. Methods Mod. Phys., {\bf 4} (2007), 669--705.
%



\bibitem{Su3} 
H.J. Sussmann and  J.C. Willems, Three Centuries of Curve Minimization: From the
Brachistochrone to Modern Optimal Control Theory, 2003  {\it (a book downloadable from 
www.math.rutgers.edu/{\tiny$\sim$}sussmann/papers/main-draft.ps.gz)}.
%\bibitem{Su} H. J. Sussmann, {\it A Local Second- and Third-order Maximum Principle} in ``Proceeding of the American Control Conference, 
%Anchorage, AK, May 8-10, 2002'', American Automatic Control Council, Evanston, IL, 2002.
 
%  \bibitem[]{}
%  
%    \bibitem[]{}
%  
%  ==========================
%  
%    \bibitem[]{}
%    
%      \bibitem[]{}
%
%\bibitem[BT]{BT} S. Benenti and W. M. Tulczyjew, {\it Symplectic linear relations}, 
% Mem. Accad. Sci. Torino Cl. Sci. Fis. Mat. Natur. (5) {\bf 5} (1981),  71--140. 

%
%
%\bibitem[GG]{GG} G. Gallavotti, {\it 
%The elements of mechanics.} Texts and Monographs in Physics. Springer-Verlag, New York, 1983. xiv+575.
%
%\bibitem[LC]{LC} T. Levi-Civita, 
%
%\bibitem[Ja]{Ja} R. D. James, {\it Co-existent Phases in the One-Dimensional Static Theory of Elastic Bars}, 
%Arch.   Rat. Mech. and Analysis, {\bf 72} (1979), 99 -- 140.
%
%
%
%
%
%\bibitem[Si]{Si} A. Signorini,  
%
%
%\bibitem[ST]{ST} A. Sniatycki and W. M. Tulczyjew, 
%{\it Generating forms of Lagrangian submanifolds}, 
%Indiana Math. J.  {\bf 22} (1972/73),  267--275.
%
%
%\bibitem[Tu1]{Tu1} W. M. Tulczyjew, 
%{\it Hamiltonian systems, Lagrangian systems and the Legendre transformation}, 
%Istituto Nazionale di Alta Matematica, Symposia Matematica, Vol. XIV (1974), 247--258.
%
%\bibitem[Tu2]{Tu2} W. M. Tulczyjew, 
%{\it Gauge independent formulation of dynamics of charged extended objects},
%Ann. Inst. Henri Poincar\'e {\bf XXXIV} (1981), 25--43.
%
%\bibitem[Tu3] {Tu3} W. M. Tulczyjew, 
%{\it The origin of variational principles}, 
%in ``Classical and Quantum Integrability'', Banach Center Publications, Vol. 59, Institute of Mathematics, 
%Polish Academy of Sciences, Warsawa, 2003. 
%
%\bibitem[Tu4]{Tu4} W. M. Tulczyjew, 
%{\it Variational formulations I: Statics of mechanical systems},
%Comm. in Mathematics {\bf 19} (2011), 179--206.
%
%
%\bibitem[TZ]{TZ} W. M.  Tulczyjew and S. Zakrzewski,  {\it  Partial categories} 
%Atti Accad. Sci. Torino Cl. Sci. Fis. Mat. Natur. {\bf 121} (1987), 46--52.
%
%\bibitem[TZ1]{TZ1} W. M.  Tulczyjew and S. Zakrzewski, {\it Composition classes and partial categories},  Atti Accad. Sci. Torino Cl. Sci. Fis. Mat. Natur. {\bf 122} (1988),  148--154.
%
%
%\bibitem[TZ2]{TZ2} W. M.  Tulczyjew and S. Zakrzewski, {\it Partial categories of differentiable relations}. Atti Accad. Sci. Torino Cl. Sci. Fis. Mat. Natur. {\bf 124}  (1990),   161--168.


%\bibitem[Vu]{Vu} S. Vukmirovi\v c, {Para-quaternionic reduction}, posted in ArXiv:math/0304424v1


\end{thebibliography}
\end{document}